\documentclass[12pt]{amsart}
\usepackage[usenames,dvipsnames]{xcolor}
\usepackage{amsmath,amssymb,amsthm,graphicx,mathrsfs,url,latexsym,enumerate}
\usepackage{a4wide}
\usepackage{ifthen}
\usepackage{verbatim}
\usepackage{fancyhdr}
\usepackage{url}
\usepackage{bbm}
\usepackage[sans]{dsfont}
\usepackage{transparent}
\usepackage{amsbsy}
\usepackage[colorlinks=true,linkcolor=Red,citecolor=Green]{hyperref}
\usepackage{wrapfig}
\usepackage{tikz}

\definecolor{bleu}{RGB}{27,88,145}
\definecolor{mauve}{RGB}{138,20,79}

\renewcommand{\Im}{\operatorname{Im}}
\renewcommand{\Re}{\operatorname{Re}}

\newcommand{\Id}{\operatorname{Id}}

\newcommand{\E}{\mathbb E}

\newcommand{\supp}{\operatorname{supp}}

\newcommand{\diag}{\operatorname{diag}}
\renewcommand{\span}{\operatorname{span}}
\newcommand{\C}{\mathbb C}

\newcommand{\N}{\mathbb N}
\newcommand{\R}{\mathbb R}

\newcommand{\Hess}{\operatorname{Hess}}

\newcommand{\Ran}{\operatorname{Ran}}

\newcommand{\bsigma}{\boldsymbol{\sigma}}

\def\<{\langle}
\def\>{\rangle}
\newcommand{\bp}{{\it Proof. }}
\newcommand{\ep}{\hfill $\square$}
\newcommand{\be}{\begin{equation}}
\newcommand{\ee}{\end{equation}}
\newcommand{\bes}{\begin{equation*}}
\newcommand{\ees}{\end{equation*}}
\renewcommand{\d}{\operatorname{d}}

\numberwithin{equation}{section}
\numberwithin{figure}{section}

\def\m{\mathbf{m}}
\def\s{\mathbf{s}}
\def\j{\mathbf{j}}
\def\k{\mathbf{k}}
\def\um{\underline{\m}}

\everymath{\displaystyle}

\DeclareMathOperator{\eps}{\varepsilon}
\DeclareMathOperator{\phii}{\varphi}

\newcommand{\vvert}[1]{\left\lVert#1\right\rVert}
\newcommand{\lb}{\operatorname{[\negthinspace [}}
\newcommand{\rb}{\operatorname{]\negthinspace ]}}

\newtheorem{theorem}{Theorem}
\newtheorem{defin}{Definition}[section]

\newtheorem{lemma}[defin]{Lemma}

\newtheorem{proposition}[defin]{Proposition}
\newtheorem{remark}[defin]{Remark}

\newtheorem{assumption}{Assumption}

\def\bbb{{\mathcal B}}\def\ccc{{\mathcal C}}
\def\fff{{\mathcal F}} 
 \def\lll{{\mathcal L}}
\def\mmm{{\mathcal M}}\def\nnn{{\mathcal N}} \def\ppp{{\mathcal P}}

\def\uuu{{\mathcal U}}\def\vvv{{\mathcal V}}\def\www{{\mathcal W}}

\def\Dr{{\mathscr D}}
\def\Er{{\mathscr E}}

\def\Mr{{\mathscr M}}

\begin{document}
\title{Sharp spectral gap for some degenerate Witten Laplacians}
\maketitle

\begin{center}\textsc{Loïs Delande\footnote{CERMICS, École des Ponts Champs-sur-Marne, France}\footnote{MATHERIALS, Inria Paris, Paris, France }}\end{center}

\begin{abstract}
    We consider Witten Laplacians associated to some non-Morse potentials. We prove Eyring–Kramers formulas for the bottom of the spectrum of these operators in the semiclassical regime and quantify the spectral gap separating these eigenvalues from the rest of the spectrum. The main ingredient is the construction of sharp degenerate Gaussian quasimodes through an adaptation of the WKB method.
\end{abstract}

\section{Introduction}
\subsection{Motivations}
The Witten Laplacian associated to a smooth function $V:M\to\R$, where $M$ is $\R^d$ or a compact manifold of dimension $d$ was introduced by Witten \cite{Wi82} in order to give an analytical proof of the Morse inequalities. For $h>0$, this semiclassical operator acting on functions is defined as
\be\label{eq:WittenLap}
\Delta_V = -h^2\Delta + |\nabla V|^2 - h\Delta V.
\ee
Here we will just focus on the Witten Laplacian over $\R^d$, although due to the very local properties of the objects we consider, one should be able to generalize the results of this paper to compact smooth manifolds. The study of this operator quickly extended beyond the field of differential geometry. In stochastic calculus for example we can consider the overdamped Langevin dynamic with drift $2\nabla V$
\be\label{eq:EDS}
dX_t = -2\nabla V(X_t)dt + \sqrt{2h}dB_t
\ee
with $B_t$ a Brownian motion. It appears that its generator, the Kramers-Smoluchowski operator
\bes
\lll = -h\Delta + 2\nabla V\cdot \nabla
\ees
is linked to the Witten Laplacian through the conjugation
\be\label{eq:conjug}
e^{-V/h}\circ h\lll\circ e^{V/h} = \Delta_V.
\ee
The overdamped Langevin dynamic originally arose in the study of molecular systems, and in this context, $h$ is proportional to the system's temperature \cite{LeRoSt10_01}. This operator also arises in control theory \cite{LaLe22} or in dynamical system \cite{DaRi21}. The need of precision in the laws describing the eigenvalues is important for numerical simulation and computational physics. We refer to \cite{LeRoSt10_01} for further information regarding the numerical aspects.

Many points of view can be adopted to study \eqref{eq:EDS}, for example let us consider a bounded open set $\Omega$ in $\R^d$ and $x_0\in\Omega$. To $X_t^{x_0}$ the solution of \eqref{eq:EDS} which starts at $X_0^{x_0} = x_0$, we associate an exit time
\bes
\tau = \inf\{t>0,\ X_t^{x_0}\notin\Omega\},
\ees
The question of the expected exit time is an important one and has been extensively studied in the low temperature regime $h\to0$. The intuition was that the mean exit time is the inverse of the first non zero eigenvalue of the generator $\lll$ with Dirichlet boundary conditions. It was proved in some restricted cases in \cite{FrWe84}, then improved by \cite{Da83} (for example when the deterministic differential equation has only one stable critical point), both methods rely on the large deviations principle. Then this result has been extended to a wider class of SDE \cite{DiGLeLePNe19}, \cite{Ne21}, \cite{LeLePNe22}, \cite{LeMiNe23}. We can also mention \cite{Ne20}, in which the author considered a function that is not necessarily Morse on the boundary of $\Omega$. These later works use the notion of quasi-stationary distribution and their first result was to compute the expected exit time when the starting point is given by this distribution, it has then been extended to any starting point in the domain.

We also have that $\lll$ generates the following PDE
\bes
\left\{
\begin{aligned}
    &\partial_t u + \lll u = 0,\\
    &u_{|t=0}=u_0,
\end{aligned}
\right.
\ees
whose solution is $u(t,x)=\mathbb{E}(u_0(X_t^{x}))$. This gives a motivation to have a sharp description of the spectrum of $\Delta_V$. Estimates regarding the bottom of its spectrum in the semiclassical limit $h\to0$ has been done shortly after Witten's work, in \cite{HeSj84_01}, the authors obtained rather precise results when $V$ is a Morse function. Namely, they showed that it has $n$ exponentially (with respect to $h^{-1}$) small eigenvalues, $n$ being the number of minima of $V$, and they also obtained the existence of a spectral gap after the $n$-th eigenvalue of size $h$. Several years later, the estimates have been refined to Eyring-Kramers laws,
\bes
\lambda_j = a_jhe^{-2S_j/h}(1+o(1))
\ees
for $j\in\lb1,n\rb$, with explicit $a_j$ and $S_j$. This result was obtained by probabilistic methods \cite{BoGaKl05_01} and semiclassical analysis \cite{HeKlNi04_01} under a generic assumption (we discuss this assumption in the paragraphs following \eqref{ass.gener}). These results continued to be investigated and very similar results have been shown for more general Fokker-Planck operators which are not self-adjoint anymore adding difficulties finding good resolvent estimates \cite{HeHiSj08-2}, \cite{BoLePMi22}.

Because the Witten Laplacian is self-adjoint, these Eyring-Kramers laws and spectral gap immediately lead to a metastable behavior of the solution of \eqref{eq:EDS}. Which means that during a series of exponentially (with respect to $h^{-1}$) long time intervals, the process will concentrate around one eigenstate associated to a small eigenvalue in each interval before finally converging to the stable state. In general we need to apply a quantitative Gearhart-Prüss theorem \cite{HeSj10_01}, \cite{HeSj21_01}, \cite{HeSj23} to obtain it \cite[Corollary 1.5]{BoLePMi22}.

More recently, \cite{Mi19} for the Witten Laplacian without the generic assumption and \cite{BoLePMi22} for Fokker-Planck type operators gave a very sharp description of the small eigenvalues of these operators under not so strong assumptions, though still assuming that the potential $V$ is Morse. They described and showed how we can compute the semiclassical expansion of the $\lambda_j$.

Now we can wonder, what if $V$ is not necessarily Morse? Several works studied that question, we can mention for example \cite{LePNiVi20} where the authors prove an Arrhenius law under very weak assumptions, namely requiring that $V$ is Lipschitz with a finite number of singularities and critical values.

In \cite{HoSt88}, \cite{HoKuSt89} and \cite{Mi95}, considering a smooth potential with small gradient on a compact manifold, the authors showed a law in-between Arrhenius and Eyring-Kramers. More precisely they obtained the following bounds on the eigenvalues
\be\label{eq:ArrEK}
a_jh^{5d} e^{-2S_j/h} \leq \lambda_j \leq b_jh^{-3d} e^{-2S_j/h}
\ee
for some $a_j,b_j > 0$ and explicit $S_j$.

In \cite{LeNi15}, the authors built a potential in 1D and 2D where $\supp\nabla V$ has a non-empty interior, i.e. $V$ is flat on open domains of $\R^d$. Assuming the flatness only appears at the minimum of $V$ and not saddle points, they managed to prove an Eyring-Kramers law for the first eigenvalue of the Witten Laplacian in a bounded domain with Dirichlet boundary condition and no critical point on the boundary. A hard generalization is to consider $V$ to also be flat at saddle points. We will see in this article that the semiclassical constructions heavily depend on the saddle points and rather little on the minima, hence the difficulty to flatten the potential at the saddle points.

Another work under a slightly weaker assumption than Morse is \cite{AsBoMi22} where they consider Morse-Bott function, i.e. the set of critical points of $f$ is a disjoint union of boundaryless connected submanifolds $\Gamma$ of $\R^d$ such that the transversal Hessian of $f$ at any point of $\Gamma$ is non degenerate. Under this hypothesis, the authors obtained Eyring-Kramers laws similar to \cite{BoLePMi22}'s.

Our goal here is to present a result similar to \cite{BoLePMi22}'s for Witten Laplacian with degenerate potential. Degenerate in the sense that the potential vanishes at higher orders at its critical points. It will be particularly useful in a forthcoming paper concerning degenerate Fokker–Planck operators.

\subsection{Framework and results}
Let us consider a potential $V\in\ccc^\infty(\R^d,\R)$. Denoting
\be\label{eq:dv}
\d_V = e^{-V/h}\circ h\nabla\circ e^{V/h} = h\nabla + \nabla V,
\ee
one can show that the Witten Laplacian $\Delta_V$ \eqref{eq:WittenLap} defined on $\ccc^\infty_c(\R^d)$, satisfies $\Delta_V = \d_V^*\d_V$ and is therefore a non-negative operator. We thus have that it is essentially self-adjoint (see for example \cite[Theorem 9.15]{He13}), which self-adjoint extension we still denote $\Delta_V$.

\begin{assumption}\label{ass.confin}
There exist $C>0$ and a compact  set $K\subset\R^d$ such that
\bes
V(x)\ \geq\  -C,\;\;\; \vert\nabla V(x)\vert \geq \frac1C\;\;\;\text{and}\;\;\; \vert\Delta V(x)\vert \leq C \vert\nabla V(x)\vert^2\,.
\ees
for all $x\in\R^d\setminus K$.
\end{assumption}

Under this assumption, its essential spectrum is contained in $[C_0,+\infty)$ for a certain $C_0>0$ independent of $h\in(0,h_0]$ with $h_0$ small (see subsection \ref{ssec:essspec}). 
\begin{lemma}\label{lem:soussolaffine}
    Let $V$ satisfying Assumption \ref{ass.confin}, then there exists $b\in\R$ such that
    \bes
    \forall x\in\R,\ \ V(x) \geq \frac1C|x| + b
    \ees
    with $C>0$ given by Assumption \ref{ass.confin}. 
\end{lemma}
With this lemma (whose proof is postponed to the appendix), $e^{-V/h}\in L^2(\R^d)$, hence $\ker \Delta_V = \C e^{-V/h}$ and thus $0\in\sigma(\Delta_V)$.
\begin{remark}
    In \cite{BoGaKl05_01}, the confining assumption is slightly different, they assume $|\Delta V| \leq C |\nabla V|$ without the square on the gradient. But using that $|\nabla V|\geq1/C$, if $|\Delta V| \leq C |\nabla V|^\alpha$ for any $\alpha \leq 2$ then $|\Delta V| \leq C^{3-\alpha} |\nabla V|^2$. And an assumption of the form $|\Delta V| \leq C |\nabla V|^\alpha$ with $\alpha>2$ would be a problem in order to prove the lower bound on $\Delta_V$'s essential spectrum when $\nabla V$ is not bounded. Hence, the hypothesis $|\Delta V| \leq C |\nabla V|^2$ is the more general we can consider.
\end{remark}
We denote by $\uuu$ the set of critical points of $V$.
\begin{assumption}\label{ass.morse}
    For any critical point $x^*\in\uuu$, there exist a neighborhood $\vvv\ni x^*$, $(t_i^{x^*})_{1\leq i\leq d}\subset \R^{*}$, $(\nu_i^{x^*})_{1\leq i\leq d}\subset \N\setminus\{0,1\}$, a $\ccc^\infty$ change of variable $U^{x^*}$ defined on $\vvv$ such that $U^{x^*}(x^*)=x^*$, $U^{x^*}$ and $d_{x^*}U^{x^*}$ are invertible and
    \be\label{eq:morse}
    \forall x\in\vvv,\ \ V\circ U^{x^*}(x) - V(x^*) = \sum_{i=1}^dt_i^{x^*}(x_i-x_i^*)^{\nu_i^{x^*}}.
    \ee
\end{assumption}
\begin{remark}
    We notice from \eqref{eq:morse} that any $x^*\in\uuu$ is an isolated critical point, but Assumption \ref{ass.confin} implies that $\uuu\subset K$ which we recall is compact, therefore the set $\uuu$ is finite. Furthermore, if $V$ is a Morse function, $V$ satisfies Assumption \ref{ass.morse} through the Morse Lemma with $\nu_i^{x^*} = 2$ for all $i$.
\end{remark}
\begin{remark}\label{rem:relax}
    We can relax \eqref{eq:morse} to
    \be\label{eq:morsealt}
    \forall x\in\vvv,\ \ V\circ U^{x^*}(x) - V(x^*) = \sum_{i=1}^dt_i^{x^*}(x_i-x_i^*)^{\nu_i^{x^*}}(1+r_i(x-x^*)),
    \ee
    with $r_i(x) = O(x)$ smooth.
    To see this, consider the map
    \bes
    \phi : (x_1,\ldots,x_d)\to(x_1(1+r_1(x))^{1/\nu_1^{x^*}},\ldots,x_d(1+r_d(x))^{1/\nu_d^{x^*}}),
    \ees
    then $d_{0}\phi = \Id$, hence it is a $\ccc^\infty$-diffeomorphism in a neighborhood of $0$. Considering now $U^{x^*}$ that satisfies \eqref{eq:morsealt}, then $\hat U^{x^*} = U^{x^*}\circ\tau_{-x^*}\circ\phi^{-1}\circ\tau_{x^*}$ satisfies \eqref{eq:morse}, where $\tau_a(x) = x-a$ and moreover we have $d_{x^*}\hat U^{x^*} = d_{x^*}U^{x^*}$.
\end{remark}
Most of the time when the reference is clear, we will just write $U,t_i$ and $\nu_i$ instead of $U^{x^*},t_i^{x^*}$ and $\nu_i^{x^*}$ in order to simplify the notation.

We denote $\uuu$ the set of critical points of $V$, and we consider the partition $\uuu = \uuu^{odd}\sqcup\uuu^{even}$ where $x^*\in\uuu^{even} \iff \forall i,\ \nu_i^{x^*}\in 2\N$. Now we shall say that $x^*\in\uuu^{even}$ is of index $j\in\lb0,d\rb$ if $\sharp\{i\ |\ t_i^{x^*}<0\}=j$ and therefore $\uuu^{even}=\bigsqcup_{j=0}^d\uuu^{(j)}$ where $\uuu^{(j)}$ is the set of critical points of $V$ with even order in each direction and of indhat they have no conflict of interest.ex $j$, we also denote $n_0=\sharp\uuu^{(0)}$ the number of minima of $V$.

We shall also set a few notations concerning these degeneracy orders. For a set of critical points $A\subset \uuu$ we will need
\be\label{eq:notationnu}
\begin{split}
    &\overline \nu^A = \sup_{x^*\in A}\sup_i\nu_i^{x^*},\\
    &\underline \nu^A = \sup_{x^*\in A}\inf_i\nu_i^{x^*},\\
    &\underline \nu^A_* = \sup_{x^*\in A}\inf\{\nu_i^{x^*}\ |\ t_i > 0\},\\
    &\hat \nu^A = \sup_{x^*\in A}\inf\{\nu_i^{x^*}\ |\ \nu_i^{x^*} > \underline \nu^{\{x^*\}} \}.
\end{split}
\ee
In order to simplify the notations, if $A=\uuu$, we will forget the superscript and just write $\overline \nu$ for example, and if $A$ contains a single element, we will write $\overline \nu^{x^*}$ instead of $\overline \nu^{\{x^*\}}$. If $A$ contains a maximum or a homogeneous point, $\underline \nu^A_*$ or $\hat \nu^A$ can be infinite with the convention $\inf\emptyset = +\infty$, and whenever this occurs, we will discuss it to avoid any ambiguity.

With these assumptions, we can state a first result concerning the rough localization of the spectrum of $\Delta_V$.

\begin{theorem}\label{thm:1}
    Suppose that Assumptions \ref{ass.confin} and \ref{ass.morse} are satisfied. There exist $h_0,\eps > 0$ such that for all $h\in(0,h_0]$, $\sigma(\Delta_V)\cap\{\Re z<\eps h^{2-\frac2{\overline{\nu}}}\}$ is made of $n_0$ eigenvalues counted with multiplicities. Hence we can denote them $\{\lambda(\m,h), \m\in\uuu^{(0)}\}$. Moreover there exists $c>0$ independent of $h$ such that
    \bes
    \forall\m\in\uuu^{(0)},\ \lambda(\m,h) \leq e^{-c/h}.
    \ees
\end{theorem}

Assumption \ref{ass.morse} is enough to prove Theorem \ref{thm:1} and the work done in section \ref{sec:rough}, but from section \ref{sec:sharp} onward and for the following theorem on a more precise description of the low-lying eigenvalues, we need to have a slightly stronger assumption:
\begin{assumption}\label{ass.morsebis}
    The map $d_{x^*}U^{x^*}$ in Assumption \ref{ass.morse} is unitary instead of merely invertible.
\end{assumption}

For the rough localization of Theorem \ref{thm:1}, we only need to build cutoffs around the critical points of $V$ but we do not need to be very precise about these cutoffs. In order to obtain sharp asymptotics, we have to refine our constructions. It requires the introduction of the following topological definitions we recall for the reader's convenience from \cite[Definition 1.3]{BoLePMi22} and what follows. But before stating this definition we need to recall a result which is well known in the Morse case, but for which we did not find references in the degenerate case of this paper.

\begin{proposition}\label{prop:saddle}
    Let $\s\in\uuu$, for $r>0$ small enough, the set $B(\s,r)\cap\{V < V(\s)\}$ has exactly two connected components if $\s\in\uuu^{(1)}$ and one otherwise (where $B(\s,r)$ denotes the open ball of center $\s$ and radius $r$).
\end{proposition}

We postpone the proof of this proposition to the Appendix.

\begin{defin}\label{def:labeling}
    We say that $\s\in\uuu$ is a separating saddle point if, for every $r>0$ small enough, $\{x\in B(\s,r),\ V(x)<V(\s)\}$ is composed of two connected components that are contained in two different connected components of $\{x\in\R^d,\ V(x)<V(\s)\}$. Hence, from Proposition \ref{prop:saddle}, the set of these points is a subset of $\uuu^{(1)}$, we will then denote it by $\vvv^{(1)}$, we also denote by separating saddle value of $V$ a point in $V(\vvv^{(1)})$.

    We say that $E\subset\R^d$ is a critical component of $V$ if there exists $\sigma\in V(\vvv^{(1}))$ such that $E$ is a connected component of $\{V\leq\sigma\}$ and $\partial E\cap \vvv^{(1)}\neq\emptyset$.
\end{defin}

Let us now recall the labeling procedure. Under the assumptions \ref{ass.confin} and \ref{ass.morse}, $V(\vvv^{(1}))$ is finite. We denote $N = \sharp V(\vvv^{(1})) + 1$ and by $\sigma_2 > \sigma_3 > \ldots > \sigma_N$ the elements of $V(\vvv^{(1}))$, for convenience, we also introduce an artificial infinite saddle value $\sigma_1 = +\infty$. Starting from $\sigma_1$, we recursively associate to each $\sigma_i$ a finite family of local minima $(\m_{i,j})_j$ and a finite family of critical component $(E_{i,j})_j$:
\begin{itemize}
    \item Let $X_{\sigma_1} = \{x\in\R^d,\ V(x) < \sigma_1\} = \R^d$. We let $\m_{1,1}$ be any global minimum of $V$ and $E_{1,1} = \R^d$. In the following we will write $\um = \m_{1,1}$.
    \item Next, we consider $X_{\sigma_2} = \{x\in \R^d,\ V(x) < \sigma_2\}$. This is the union of its finitely many connected components. Exactly one contains $\um$ and the other components are denoted by $E_{2,1},\ldots,E_{2,N_2}$. They are all critical and, in each component $E_{2,j}$ we pick up a point $\m_{2,j}$ which is a global minimum of $V_{|E_{2,j}}$.
    \item Suppose now that the families $(\m_{k,j})_j$ and $(E_{k,j})_j$ have been constructed until rank $k = i-1$. The set $X_{\sigma_i} = \{x\in\R^d,\ V(x) < \sigma_i\}$ has again finitely many connected components and we label $E_{i,j}$, $j\in\lb1,N_i\rb$, those of these components that do not contain any $\m_{k,j}$ for $k<i$. They are all critical and, in each $E_{i,j}$, we pick up a point $\m_{i,j}$ which is a global minimum of $V_{|E_{i,j}}$.
\end{itemize}
At the end of this procedure, all the minima have been labeled. Throughout, we denote by $\s_1$ a fictive saddle point such that $V(\s_1) = \sigma_1 = \infty$ and, for any set $A$, $\ppp(A)$ denotes the power set of $A$. From the above labeling, we define two mappings
\bes
E:\uuu^{(0)}\to\ppp(\R^d)\ \ \mbox{ and }\ \ \j:\uuu^{(0)}\to\ppp(\vvv^{(1)}\cup\{\s_1\})
\ees
as follows: for every $i\in\lb1,N\rb$ and every $j\in\lb1,N_i\rb$,
\bes
E(\m_{i,j}) = E_{i,j},
\ees
and
\bes
\j(\um) = \{\s_1\}\ \ \mbox{ and }\ \ \j(\m_{i,j}) = \partial E_{i,j}\cap \vvv^{(1)}\ \mbox{ for }i\geq2.
\ees
In particular, we have $E(\um) = \R^d$ and, for all $i\in\lb1,N\rb$, $j\in\lb1,N_i\rb$, one has $\emptyset\neq\j(\m_{i,j})\subset\{V=\sigma_i\}$. We then define the mappings
\bes
\bsigma:\uuu^{(0)}\to V(\vvv^{(1)})\cup\{\sigma_1\}\ \ \mbox{ and }\ \ S:\uuu^{(0)}\to(0,+\infty],
\ees
by
\be\label{eq:defS}
\forall \m\in\uuu^{(0)},\ \ \bsigma(\m) = V(\j(\m))\ \ \mbox{ and }\ \ S(\m) = \bsigma(\m) - V(\m),
\ee
where, with a slight abuse of notation, we have identified $V(\j(\m))$ with its unique element. Note that $S(\m) = \infty$ if and only if $\m = \um$. 

We now consider the following generic assumption in order to lighten the result and the proof.
\be\label{ass.gener}\tag{Gener}
\begin{array}{l}
    (\ast)\mbox{ for any }\m\in\uuu^{(0)},\m \mbox{ is the unique global minimum of } V_{|E(\m)}\\
    (\ast)\mbox{ for all }\m\neq\m'\in\uuu^{(0)},\j(\m)\cap\j(\m')=\emptyset.
\end{array}
\ee
In particular, \eqref{ass.gener} implies that $V$ uniquely attains its global minimum at $\um$. This assumption allows us to avoid some heavy constructions regarding the set $\uuu$ and lighten the definition \ref{def:quasimodes} of the quasimodes, see \cite{Mi19} for the general setting. But it seems that its not a true obstruction and that we can pursue the computations without this assumption as described in \cite[Section 6]{BoLePMi22}, \cite{No24}, in the spirit of \cite{Mi19} and we show an example of what can happen when \eqref{ass.gener} is not satisfied in Section \ref{sec:triple}.

This assumption comes from \cite[(1.7)]{BoGaKl05_01} and \cite[Assumption 3.8]{HeKlNi04_01} where they were hard to handle properly. Then it changed to become \cite[Hypothesis 5.1]{HeHiSj11_01} when finally reaching the form of \cite[Assumption 4]{LePMi20}.

One can show that \eqref{ass.gener} is weaker than \cite{BoGaKl05_01}'s, \cite{HeKlNi04_01}'s and \cite{HeHiSj11_01}'s assumptions. More precisely, they supposed that $\j(\m)$ is a singleton while here we have no restriction on the size of $\j(\m)$.

For the proof of Theorem \ref{thm:2} below, while computing the eigenvalues, a term we would want to be an error one appears, it is of order $h^{\alpha_0}$ with $\alpha_0$ defined in \eqref{eq:alpha0}. The following assumption is to ensure this term is indeed a sub-principal one.
\begin{assumption}\label{ass.alpha}
    We consider $\alpha_0$ as in \eqref{eq:alpha0}. We have $\alpha_0 > 0$.
\end{assumption}

This assumption is necessary to complete the calculus and does not appear to be merely a computational barrier although we did not try to prove it. In Section \ref{sec:example} we go through a quite large list of potentials that satisfy this assumption. That being said, we have all the ingredients needed to state the sharp description theorem.

\begin{theorem}\label{thm:2}
    Suppose that Assumptions \ref{ass.confin}, \ref{ass.morse}, \ref{ass.morsebis} and \ref{ass.alpha} are satisfied, assume also \eqref{ass.gener} hold true. There exist $h_0,\beta > 0$ such that for all $h\in(0,h_0]$, one has $\lambda(\um,h) = 0$ and for all $\m\in\uuu^{(0)}\setminus\{\um\}$, $\lambda(\m,h)$ satisfies the following Eyring–Kramers type formula
    \bes
    \lambda(\m,h) = v(\m) h^{\mu(\m)}e^{-2S(\m)/h}(1+O(h^\beta)),
    \ees
    where $v(\m)$ is defined in \eqref{eq:defvmu} and $\mu(\m)$ is defined in \eqref{eq:defmu}. They depend explicitly on $V$: both depend on the $\nu_i^{x^*}$ and only $v(\m)$ depends on the $t_i^{x^*}$'s. Finally, $S$ is the standard height function defined by \eqref{eq:defS}.
\end{theorem}

\begin{remark}
    If we assume that $V$ is a Morse function, we retrieve the now well-known result about small eigenvalues of the Witten Laplacian \cite{BoGaKl05_01}\cite{HeKlNi04_01}.
\end{remark}

From \eqref{eq:defvmu}, notice that $\mu(\m)\in\big(\frac{3-d}{2},\frac{3+d}{2}\big)$, thus in high dimension ($d\geq4$) it is possible to have $\mu(\m)<0$. We also observe that the estimate \eqref{eq:ArrEK} is still valid in our settings since $\big(\frac{3-d}{2},\frac{3+d}{2}\big)\subset (-3d,5d)$.

We can see here a link between our result and some work done on Schrödinger operators. Let us denote
\bes
P_W = -h^2\Delta + W
\ees
the Schrödinger operator with potential $W$. In \cite{BoPo19}, the authors showed that for our potential $V$, the first eigenvalue $\lambda_1$ of $P_{|\nabla V|^2}$ goes to $0$ slower than $h^2$ : $\frac{\lambda_1}{h^2}\xrightarrow[h\to0]{} +\infty$. Considering that $h\Delta V$ shifts some eigenvalues near $0$, we can confirm \cite{BoPo19}'s result with our bound on the spectral gap, being of size $h^{2-\frac{2}{\overline{\nu}}}$. Their work is also very interesting concerning infinitely flat potential like $e^{-|x|^\alpha}$. In this paper we do not consider such degenerate potential, but in this context the result of \cite{BoPo19} still hold in dimension $d=1$ and can be extended to all the eigenvalues.

In \cite{BeGe10}, the authors considered two-wells potentials satisfying an assumption similar to Assumption \ref{ass.morse}. Under such a hypothesis, they managed to compute the mean time for a solution of \eqref{eq:EDS} to reach one well of $V$ starting from the other. We know (see \cite{BoGaKl05_01}) that under some generic assumptions (such that the potential is Morse, which is not the case here), this exit time is the inverse of the eigenvalue associated with that well of the Witten Laplacian over $\R^d$ up to a factor of the form $1 + O(e^{-c/h})$. Hence, assuming the link between the eigenvalue and the mean time is still valid in this case (which is far from being proven), we can compare the results obtained by \cite{BeGe10} and ours. They mainly stated their results with some specific potentials, and in each cases, they obtained Eyring-Kramers law of the form
\bes
\E(\tau) = ah^\alpha e^{S/h}(1+o(1))
\ees
where we denoted $\E(\tau)$ the expected reaching time, and $a$, $\alpha$, $S$ depend explicitly on the potential (in fact the $o(1)$ is more precise in their results but we are not very interested in this comparison here). Since the $S$ in the exponential is the height of the well as usual and as it has been proven under weaker assumptions, we will not linger on that matter. Moreover, one can check that the prefactors $a$ are the same for both our works. So let us pay attention to the power of $h$ befront the exponential : $\alpha$. The potentials they studied have two non-degenerate minima and a unique saddle point, which is degenerate, we give the expression of the potential near that saddle point they considered and the power they obtained,
\begin{itemize}
    \item $V_1 = -x_1^4 + \sum_{j\geq2}x_j^2$, with $\alpha = -\frac14$;
    \item $V_2 = -x_1^{2p} + \sum_{j\geq2}x_j^2$ with $\alpha = -\frac{p-1}{2p}$;
    \item $V_3 = -x_1^2 + x_2^4 + \sum_{j\geq3}x_j^2$ with $\alpha = \frac14$; and
    \item $V_4 = -x_1^2 + x_2^4 + x_3^4 + \sum_{j\geq4}x_j^2$ with $\alpha = \frac12$.
\end{itemize}
Looking at Section \ref{sec:example} we see that those potentials satisfy Assumption \ref{ass.alpha}, and we would like to observe that $\alpha = -\mu$ with $\mu$ defined in \eqref{eq:defvmu}. We actually observe that $\alpha = 1-\mu$ which is explained by a multiplication by $h$ when conjugating the generator of the SDE to obtain the the Witten Laplacian \eqref{eq:conjug}. Therefore, we see that the powers match perfectly giving hope for the link between mean reaching time and first non-zero eigenvalue to also be true for non-Morse potentials.

\subsection*{Acknowledgements}
The author is grateful to Laurent Michel for his advice through this work and to Jean-François Bony for helpful discussions.

\subsection*{Declarations}
This work is supported by the ANR project QuAMProcs 19-CE40-0010-01. The author declare that he has no conflict of interest. Due to the subject of this paper and its production, no data were used, and no human or animal participants were involved.

The rest of the paper is organized as follows. In the next section, we prove a hypocoercivity result regarding the operator $P$, leading to the proof of Theorem \ref{thm:1}. It gives a rough localization of the eigenvalues of $P$ as well as a resolvent estimate which is crucial for sharp spectral estimates. In Section \ref{sec:sharp}, we detail the precise construction of the quasimodes around a given saddle point of the potential $V$. Section \ref{sec:geo} focuses on the global definition of these quasimodes. We also develop the study of the interaction matrix resulting in the proof of Theorem \ref{thm:2}, giving the Eyring-Kramers law for the small eigenvalues of $P$. We continue in Section \ref{sec:example} by a non-exhaustive list of potentials that satisfy the assumptions made in Theorem \ref{thm:2}. At last, in Section \ref{sec:triple}, we give a hint about arising behaviors when Assumption \eqref{ass.gener} is no longer satisfied.

\section{Resolvent and rough spectrum estimates}\label{sec:rough}

Let $0\leq\chi_{x^*}\leq1$, $x^*\in \uuu$ be a family of cutoffs in $\ccc_c^\infty(\R^d)$ such that $\chi_{x^*}$ is supported in $B(x^*,r)$ for some $r>0$ to be chosen small enough so that the supports are disjoint, $U^{x^*}$ is defined on $\supp\chi_{x^*}$, and $\chi_{x^*}=1$ near $x^*$.
Then, for $\m\in\uuu^{(0)}$, we introduce the quasimodes
\bes
f_\m(x)=\chi_\m(x) e^{-(V(x)-V(\m))/h}|\det d_x(U^{\m})^{-1}|,
\ees
and hence the vector space 
\bes
F_h=\span\{f_\m,\;\m\in\uuu^{(0)}\}
\ees
has dimension $n_0$.

\begin{proposition}\label{prop:gapWitten}
    There exists $h_0>0$ such that for all $h\in(0,h_0]$,
    \bes
    \exists C>0,\ \forall u\in F_h^\bot\cap D(\Delta_V)\ \ \<\Delta_Vu,u\> \geq Ch^{2-\frac2{\overline{\nu}}}\vvert{u}^2,
    \ees
    recalling $\overline{\nu}=\max_{i,x^*}\nu_i^{x^*}$.
\end{proposition}

\bp Let us consider an extra non-negative smooth cutoff $\chi_\infty$ such that $\chi_\infty \equiv 1$ at infinity, $\chi_\infty\equiv0$ around each minimum $\m\in\uuu^{(0)}$, and $\sum_{a\in\uuu^{(0)}\cup\{\infty\}}\chi_a^2$ never vanishes (therefore, it is uniformly bounded from below). Consider then the family $\hat\chi_a = \frac{\chi_a}{\sqrt{\sum_{b\in\uuu^{(0)}\cup\{\infty\}}\chi_b^2}}$ for $a\in\uuu^{(0)}\cup\{\infty\}$. This way, the family $(\hat\chi_a)_{a\in\uuu\cup\{\infty\}}$ forms a smooth partition of unity in the sense of \cite[Definition~3.1]{CyFrKiSi87_01} and we have $\hat\chi_{x^*} = \chi_{x^*}$ around $x^*$. In the following we will omit the hat in order to simplify the notation.

Applying the IMS localization formula \cite[Theorem 3.2]{CyFrKiSi87_01}, let $u\in F_h^\bot\cap D(\Delta_V)$,
\be\label{eq:IMS}
\<\Delta_V u,u\> = \sum_{a\in\uuu\cup\{\infty\}} (\vvert{\d_V(\chi_au)}^2 - h^2\vvert{\nabla\chi_au}^2).
\ee

For $a\in\uuu$, let $U$ be the substitution $U^a$ of Assumption \ref{ass.morse} and $p=\sum_{i=1}^dt_i(x_i-a_i)^{\nu_i}$ the polynomial defined on $\R^d$ which coincides with $V\circ U$ where defined. Let $v\in D(\Delta_V)$,
\bes
\begin{split}
    (\d_Vv)\circ U &= (h\nabla v + \nabla Vv)\circ U = h(d_xU^T)^{-1}\nabla(v\circ U) +  (d_xU^T)^{-1}\nabla(V\circ U)v\circ U\\
    &= (d_xU^T)^{-1}(h\nabla\tilde v + \nabla \tilde V\tilde v) = (d_xU^T)^{-1}\d_{p}\tilde v,
\end{split}
\ees
where we denoted $\tilde f = f\circ U$, and $\d_p$ is the twisted derivative associated with $p$ \eqref{eq:dv}. Then, since $d_aU$ is invertible, there exist $C>0$ and $\vvv$ a neighborhood of $a$ such that
\be\label{eq:jacob}
\forall x\in\vvv,\quad\frac1C \leq |\det d_{U^{-1}(x)}U|\leq C \mbox{ and } (d_xU^Td_xU)^{-1} \geq \frac1C.
\ee
Assuming $r>0$ is small enough to have $\supp\tilde\chi_a\subset \vvv$, we obtain for $u\in F_h^\bot\cap D(\Delta_V)$
\bes
\begin{split}
    \vvert{\d_V(\chi_au)}^2 &= \int_{\R^d} |\d_V(\chi_au)|^2 \geq C^{-1}\int_{U(\vvv)} |\d_V(\chi_au)|^2|\det d_{U^{-1}(x)}U|^{-1}\\
    &\geq C^{-1}\int_\vvv |\d_V(\chi_au)|^2\circ U = C^{-1}\int_\vvv \<(d_xU^Td_xU)^{-1}\d_p(\widetilde{\chi_au}),\d_p(\widetilde{\chi_au})\>\\
    &\geq C^{-2}\int_\vvv |\d_p(\widetilde{\chi_au})|^2 = C^{-2}\<\Delta_p \widetilde{\chi_au},\widetilde{\chi_au}\>,
\end{split}
\ees
where we denote by $\Delta_p$ the Witten Laplacian associated with $p$.

For $b\in\R^d$ and $\lambda\in\R^*$, denoting $T(b)f(x) = f(x-b)$ unitary over $L^2(\R^d)$ and $D(\lambda)f(x) = \sqrt{\lambda}f(\lambda x)$ unitary over $L^2(\R)$,
\be\label{eq:DK}
\Delta_p = T(a)\sum_{i=1}^dD(h^{-\frac{1}{\nu_i}})h^{2-\frac{2}{\nu_i}}K_i^aD(h^{\frac{1}{\nu_i}})T(-a)
\ee
where
\bes
K_i^a = -\partial_i^2 + (t_i\nu_i)^2x_i^{2\nu_i-2}-\nu_i(\nu_i-1)t_ix_i^{\nu_i-2},
\ees
which is the non-semiclassical Witten Laplacian associated with the one-dimensional potential $t_ix_i^{\nu_i}$. Thus, $\sigma(K_i^a)\subset\R_+$ and
\bes
0\in\sigma(K_i^a) \iff e^{-t_ix_i^{\nu_i}}\in L^2(\R) \iff t_i>0\mbox{ and }\nu_i\in2\N,
\ees
hence
\bes
0\in\sigma(\Delta_p) \iff \forall i,\ t_i>0\mbox{ and }\nu_i\in2\N \iff a\in\uuu^{(0)}.
\ees
Moreover, in this case, we know that its associated eigenfunction is $e^{-p/h}$, and the next smallest eigenvalue is of the form $ch^{2-\min_i\frac{2}{\nu_i}}$ using \eqref{eq:DK}.

Thus, there exists $c'>0$ such that, if $a$ is not a minimum,
\bes
\forall u\in F_h^\bot\cap D(\Delta_V),\ \ \vvert{\d_V(\chi_au)}^2 \geq c'h^{2-\frac{2}{\overline\nu^a}}\vvert{\widetilde{\chi_au}}^2
\ees
using the notations \eqref{eq:notationnu}. If $a$ is a minimum, using that for $u\in F_h^\bot$, we have $u\bot\chi_ae^{-V/h}|\det dU^{-1}|$ which implies that $\widetilde{\chi_au} \bot e^{-\tilde V/h}$ thus it leads to the same result
\bes
\forall u\in F_h^\bot\cap D(\Delta_V),\ \ \vvert{\d_V(\chi_au)}^2 \geq c'h^{2-\frac{2}{\overline\nu^a}}\vvert{\widetilde{\chi_au}}^2.
\ees
Using \eqref{eq:jacob}, on $\vvv$, $|\det d_{x}U^{-1}|\geq \frac1C$, and thus $\vvert{\widetilde{\chi_au}}^2\geq \frac1C\vvert{\chi_au}^2$. Therefore
\be\label{eq:chia}
\exists c>0,\ \forall a\in\uuu,\ \forall u\in F_h^\bot\cap D(\Delta_V),\ \ \vvert{\d_V(\chi_au)}^2 \geq ch^{2-\frac{2}{\overline{\nu}}}\vvert{\chi_au}^2.
\ee
For $\chi_\infty$, let us write for $u\in D(\Delta_V)$
\bes
\begin{split}
    \vvert{\d_V(\chi_\infty u)}^2 &= \vvert{h\nabla(\chi_\infty u)}^2 + \vvert{\nabla V\chi_\infty u}^2 - h\<\Delta V\chi_\infty u,\chi_\infty u\>\\
    &\geq \vvert{\nabla V\chi_\infty u}^2 - h\sqrt C\vvert{\nabla V\chi_\infty u}^2
\end{split}
\ees
for some $C>0$, noticing that $\chi_\infty$ localizes outside the critical points of $V$, and using Assumption \ref{ass.confin}. Then for $h$ small enough, there exists $c>0$,
\be\label{eq:chiinf}
\forall u\in F_h^\bot\cap D(\Delta_V),\ \ \vvert{\d_V(\chi_\infty u)}^2 \geq \frac 1{2C}\vvert{\chi_\infty u}^2 \geq ch^{2-\frac{2}{\overline{\nu}}}\vvert{\chi_\infty u}^2.
\ee

Combining \eqref{eq:IMS}, \eqref{eq:chia}, and \eqref{eq:chiinf} we have
\bes
\exists c>0,\ \forall u\in F_h^\bot\cap D(\Delta_V),\ \ \<\Delta_V u,u\> \geq ch^{2-\frac{2}{\overline{\nu}}}\vvert{u}^2 - h^2 \sum_{a\in\uuu\cup\{\infty\}}\vvert{\nabla\chi_au}^2
\ees
using that $\sum_{a\in\uuu\cup\{\infty\}}\chi_a^2 = 1$. Taking $h\in(0,h_0]$ with $h_0>0$ small enough concludes the proof.

\ep

\subsection{Proof of Theorem \ref{thm:1}}
Let $G_h = \span\{g_\m,\ \m\in\uuu^{(0)}\}$ with
\bes
g_\m = \chi_\m e^{-(V-V(\m))/h}.
\ees
Observe that $\dim(G_h) = \dim(F_h) = n_0$. Using \eqref{eq:dv}, we have 
\bes
\d_Vg_\m = h\nabla\chi_\m e^{-(V-V(\m))/h}.
\ees
Since $\chi_\m = 1$ in a neighborhood of $\m$ implying that $\nabla\chi_\m$ is supported away from $\m$, there exists $c'>0$ such that for $h$ small enough
\bes
\<\Delta_V g_\m,g_\m\> = \vvert{\d_Vg_\m}^2 \leq e^{-c'/h}.
\ees
Using a degenerate Laplace method \eqref{eq:Laplace2}, there exists $c''>0$ such that, for $h$ small enough,
\bes
\vvert{g_\m}^2 \geq c''h^{\sum_{i=1}^d\frac{1}{\nu_i^\m}}
\ees
which leads to
\bes
\exists c,h_0\>0,\ \forall h\in(0,h_0],\ \forall u\in G_h,\ \ \<\Delta_V u,u\> \leq e^{-c/h}\vvert{u}^2.
\ees
Thus, using Proposition \ref{prop:gapWitten} and the Min-Max Theorems we obtain the spectral estimate, recalling that
\bes
\sigma_{ess}(\Delta_V)\subset[C_0,+\infty)
\ees
for some $C_0>0$ and all $h\in(0,h_0]$.

\section{Sharp quasimodes}\label{sec:sharp}

We now want to have a better view of the small eigenvalues of $\Delta_V$. To this end, we construct sharp quasimodes, following the approach of \cite[section 3 and 4]{BoLePMi22}. Their theorem does not apply here because, for a general potential $V$ satisfying Assumption \ref{ass.morse}, $\Delta_V$ does not satisfy condition (1.9), (Harmo), or (Morse) from their paper. Our goal is therefore to obtain a similar result without relying on these assumptions. As in \cite{BoLePMi22}, given $\s\in\vvv^{(1)}$, we look for an approximate solution to the equation $\Delta_V u = 0$ in a neighborhood $\www$ of $\s$.

Under Assumption \ref{ass.morse}, and up to a translation, there exist a neighborhood $\www$ of $\s$, degeneracies $(\nu_i)_i,(t_i)_i$ and a map $U$ such that $U(0) = \s$ and
\bes
\forall x\in U^{-1}(\www),\quad V\circ U(x) = V(\s) + \sum_{i=1}^dt_ix_i^{\nu_i}.
\ees
Since $\s$ has index 1, we may assume without loss of generality that $t_i > 0$ for $i \geq 2$ and $t_1 < 0$. From this point on, we work under Assumptions \ref{ass.morse} and \ref{ass.morsebis} instead of merely Assumption \ref{ass.morse}.

Due to the locality of the constructions we are about to make, the only degeneracy orders that will appear are those of $\s$. This is why throughout this section (and this section only), we will always omit the superscript $\s$ in the notations of \eqref{eq:notationnu}, hence, here, $\overline \nu$ will denote $\sup_i \nu_i^\s$, $\underline \nu$ will denote $\inf_i \nu_i^\s$, etc.

We choose an approximate solution of $\Delta_V u = 0$ of the form
\bes
u=\chi e^{-(V-V(\m))/h},
\ees
and we set
\bes
\chi(x)=\int_0^{\ell(x,h)}\zeta(s/\tau)e^{-\frac{s^{\nu_1}}{\nu_1h}}ds
\ees
where the function $\ell\in\ccc^\infty(\www)$ has a formal classical expansion $\ell\sim\sum_{j\geq0}h^j\ell_j$. Here, $\zeta$ denotes a fixed smooth even function equal to $1$ on $[-1,1]$ and supported in $[-2,2]$, and $\tau > 0$ is a small parameter which will be fixed later. Due to the degeneracy of $V$, the construction used in \cite{BoLePMi22}, which involves the standard Gaussian $e^{-s^2/(2h)}$, is not suitable here. Using that function would lead to more complex computations and less satisfactory results.

The goal of this section is to construct the function $\ell$. We adapt the construction from \cite[Section 3]{BoLePMi22}.

Since $\Delta_V(e^{-V/h}) = 0$, $\Delta_V(\chi e^{-V/h}) = [\Delta_V,\chi](e^{-V/h})$, but
\bes
\begin{split}
    [\Delta_V,\chi] &= -h^2[\Delta,\chi] = -h^2(\Delta \chi + 2\nabla \chi\cdot\nabla),\\
    \nabla \chi &= \nabla\ell\zeta(\ell/\tau)e^{-\frac{\ell^{\nu_1}}{\nu_1h}},\\
    \Delta \chi &= \Big(\Delta\ell\zeta(\ell/\tau) + \frac1\tau|\nabla\ell|^2\zeta'(\ell/\tau)-\zeta(\ell/\tau)|\nabla\ell|^2\frac{\ell^{\nu_1-1}}{h}\Big)e^{-\frac{\ell^{\nu_1}}{\nu_1h}}.
\end{split}
\ees
Therefore, if we set $\ell_0(\s) = 0$, there exists a smooth function $r$ such that $r\equiv0$ near $\s$, and both $r$ and its derivatives are locally uniformly bounded with respect to $h$. Moreover
\be\label{eq:wplusr}
\Delta_V(\chi e^{-V/h}) = h(w+r)e^{-(V+\frac{\ell^{\nu_1}}{\nu_1})/h},
\ee
where
\bes
w = 2\nabla V\cdot\nabla\ell + \ell^{\nu_1-1}|\nabla\ell|^2 - h\Delta\ell.
\ees
As with $\ell$, the function $w$ admits a formal classical expansion $w \sim \sum_{j\geq0} h^jw_j$. Formally solving $w=0$ and identifying the powers of $h$ leads to a system of equations
\be\label{eq:eik}\tag{Eik}
2\nabla V\cdot\nabla\ell_0 + \ell_0^{\nu_1-1}|\nabla\ell_0|^2 = 0
\ee
and for $j\geq1$,
\be\label{eq:transp}\tag{T$_j$}
\big((2\nabla V + 2\ell_0^{\nu_1-1}\nabla\ell_0)\cdot\nabla + (\nu_1-1)\ell_0^{\nu_1-2}|\nabla\ell_0|^2\big) \ell_j + R^j = 0
\ee
where $R^j$ is a smooth polynomial of the $\partial^\alpha\ell_k$ for $|\alpha|\leq2$ and $k<j$. As in \cite{BoLePMi22} and by analogy with the WKB method, we refer to the first equation as the eikonal equation and the subsequent ones as transport equations. Due to some technical barriers, we aim to solve these equations only approximately, to all orders in $h$; that is
\be\label{eq:eikapprox}\tag{Eik$^\infty$}
2\nabla V\cdot\nabla\ell_0 + \ell_0^{\nu_1-1}|\nabla\ell_0|^2 = O(x^\infty)
\ee
and
\be\label{eq:transpapprox}\tag{T$_j^\infty$}
\big((2\nabla V + 2\ell_0^{\nu_1-1}\nabla\ell_0)\cdot\nabla + (\nu_1-1)\ell_0^{\nu_1-2}|\nabla\ell_0|^2\big) \ell_j + R^j = O(x^\infty).
\ee

In the following Lemma \ref{lem:ell}, various local results are given under different assumptions and accordingly different levels of accuracy. Separating saddle points will be classified exhaustively later, according to the different cases. Except Assumption \ref{ass.case4}, these assumptions are obviously local ones, in a neighborhood of a critical point $\s$.

\begin{assumption}\label{ass.case1}$U$ is linear, so that $x \mapsto d_x U$ is constant. Under Assumption \ref{ass.morsebis}, this implies that for all $x$, we have $d_xU^T d_xU = d_0U^T d_0U = \Id$; in other words, $U$ is a unitary matrix.
\end{assumption}

\begin{assumption}\label{ass.case4}
    We are in a one-dimensional case.
\end{assumption}

\begin{assumption}\label{ass.case2}$V\circ U$ is exactly quadratic in the drop direction, that is, $\nu_1=2$.
\end{assumption}

\begin{assumption}\label{ass.case3}The potential has quadratic directions and the drop one is not among them. In other words, $\nu_1 > \underline \nu = 2$.
\end{assumption}

\begin{assumption}\label{ass.case5}
    There is a unique minimal degeneracy order and the drop direction is the most degenerate one. This means that there exists a unique $i\in\lb1,d\rb$ such that $\nu_i=\underline\nu$ with $i>1$, and $\nu_1 = \overline\nu$.
\end{assumption}

The apparently missing case when $\nu_1,\underline\nu>2$ and $\nu_1<\overline\nu$ will be later handled by using the rough statement of Lemma \ref{lem:ell} $i)$ below, which does not give a full asymptotic expansion.

\begin{lemma}\label{lem:ell}~

\begin{itemize}
    \item[$i)$] Considering $\ell(x) = \ell_0(x) = (2\nu_1|t_1|)^{\frac{1}{\nu_1}}U^{-1}(x)_1$, we have
    \be\label{eq:PchiCase0}
    \Delta_V(\chi e^{-V/h}) = h\Big(\sum_{j\geq2}\omega_{1,j}\circ U^{-1}\tilde t_j|U^{-1}(x)_j|^{\nu_j-1} + O(h)\Big) e^{-\big(V+\frac{\ell^{\nu_1}}{\nu_1}\big)/h}
    \ee
    with $\tilde t_j = 2(2\nu_1|t_1|)^\frac{1}{\nu_1}\nu_jt_j$ and $\omega_{1,j} = O(x)$ for $j\geq2$, and where the $O(h)$ term is uniform in $x$ over any compact set containing $\s$ in its interior.
    \item[$ii)$] Under Assumption \ref{ass.case1}, the function $\ell$ defined in $i)$ exactly solves \eqref{eq:eik} and all \eqref{eq:transp} for $j\geq1$, resulting in a significantly improved approximation:
    \be\label{eq:PchiCase1}
    \Delta_V(\chi e^{-V/h}) = O(h^\infty)e^{-\big(V+\frac{\ell^{\nu_1}}{\nu_1}\big)/h}
    \ee
    uniformly around $\s$.
    \item[$iii)$] Under Assumption \ref{ass.case4}, the function $\ell$ defined in $i)$ exactly solves \eqref{eq:eik} resulting in the following estimate
    \be\label{eq:PchiCase4}
    \Delta_V(\chi e^{-V/h}) = O(h^2)e^{-\big(V+\frac{\ell^{\nu_1}}{\nu_1}\big)/h}
    \ee
    uniformly around $\s$.
    \item[$iv)$] Under Assumption \ref{ass.case2}, there exists a sequence $(\ell_j)_{j\geq0}\subset\ccc^\infty(\www)$ such that $\ell \sim \sum_{j\geq0}\ell_jh^j$ in $\ccc^\infty(\www)$,  which solve \eqref{eq:eikapprox} and \eqref{eq:transpapprox} resulting in
    \be\label{eq:PchiCase2}
    \Delta_V(\chi e^{-V/h}) = O(x^\infty)e^{-\big(V+\frac{\ell^{2}}{2}\big)/h}
    \ee
    uniformly around $\s$ and with respect to $h$.
    \item[$v)$] Under Assumption \ref{ass.case3}, we can solve higher orders of the eikonal equation, more precisely, we can construct an $h$-independent function $\ell$ (equal to $\ell_0$ from $i)$ at principal order) such that, ordering the variables such that $\nu_i = 2 \iff i\in\lb2,J\rb$ for a certain $J\in\lb2,d\rb$ and recalling $\hat\nu = \min_{\nu_i\neq2}\nu_i$ from \eqref{eq:notationnu},
    \be\label{eq:PchiCase3}
    \begin{split}
        \Delta_V(\chi e^{-V/h}) &= hO\bigg(|U^{-1}(x)|\sum_{j\notin\lb2,J\rb}|U^{-1}(x)_j|^{\nu_j-1}\\&\phantom{*******} + |U^{-1}(x)|^{\hat\nu-1}\sum_{j\in\lb2,J\rb}|U^{-1}(x)_j| + h\bigg)e^{-\big(V+\frac{\ell^{\nu_1}}{\nu_1}\big)/h}
    \end{split}
    \ee
    uniformly around $\s$.
    \item[$vi)$] Under Assumption \ref{ass.case5}, choosing the last variable to be the least degenerate, the hypothesis can be stated as $\forall i<d,\ \nu_i>\nu_d$ and $\nu_1=\overline\nu$. Then we obtain an estimate similar to that of the previous case, namely
    \be\label{eq:PchiCase5}
    \begin{split}
        \Delta_V(\chi e^{-V/h}) &= hO\bigg(|U^{-1}(x)|\sum_{j<d}|U^{-1}(x)_j|^{\nu_j-1}\\&\phantom{*******} + |U^{-1}(x)|^{\hat\nu-\underline\nu+1}|U^{-1}(x)_d|^{\nu_d-1} + h\bigg)e^{-\big(V+\frac{\ell^{\nu_1}}{\nu_1}\big)/h}
    \end{split}
    \ee
    uniformly around $\s$.
    \item[$vii)$] In all cases, the chosen function $\ell$ elliptizes $V$ around $\s$ according to
    \bes
    \big(V+\frac{\ell_0^{\nu_1}}{\nu_1}\big)\circ U(x) = V(\s) + \sum_{i=1}^d|t_i|x_i^{\nu_i}(1+O(x)).
    \ees
\end{itemize}
\end{lemma}

By $\ell \sim \sum_{j\geq0}\ell_jh^j$ in $\ccc^\infty(\www)$, we mean that for all $n\in\N$, there exists $C_n$ such that
\be\label{eq:classicexpansion}
\vvert{\ell - \sum_{j=0}^n\ell_jh^j}_{\ccc^\infty(\www)} \leq C_n h^{n+1}.
\ee

\begin{remark}
    In the usual Morse case, \eqref{eq:eik} is solved through topological and geometric constructions, see \cite[Lemma 3.2]{BoLePMi22}. They use the non-degeneracy of the potential and its coefficients to extract a Lagrangian manifold which projects nicely onto the $x$-space, its generating function serves as the starting point for constructing $\ell_0$. Here, it is far from clear whether such a manifold exists, and even if it does, Remark \ref{rem:counterexample} suggests that the standard method likely cannot be carried through to completion. Therefore, we must pursue an alternative approach to construct $\ell_0$. Our method consists in mimicking the resolution of \eqref{eq:transp} using homogeneous polynomials.
\end{remark}

\bp In this proof, we denote $\tilde f = f\circ U$ for any function $f$ for which this composition makes sense. Noticing that for any $f\in\ccc^1$, $\nabla\tilde f = d_xU^T\nabla f\circ U$, we have
\bes
\tilde w = 2\Omega\nabla\tilde V\cdot\nabla\tilde\ell + \tilde\ell^{\nu_1-1}\Omega\nabla\tilde\ell\cdot\nabla\tilde\ell - h\Delta\ell\circ U
\ees
where $\Omega = (d_xU^Td_xU)^{-1}$, thus $\Omega = \Id + O(x)$ from Assumption \ref{ass.morsebis}. Notice that $\tilde \ell$ and $\tilde w$ have a decomposition similar to $\ell$ and $w$ : $\tilde \ell \sim \sum_{j\geq0}\tilde\ell_jh^j$ and $\tilde w \sim \sum_{j\geq0}\tilde w_jh^j$. Instead of working with the $w_j$, it will be handier to work with the $\tilde w_j$ (and note that $(w\circ U)_j = w_j\circ U$, so there is no ambiguity in the notation of $\tilde w_j$ and $\tilde \ell_j$). Through this change of variables, the eikonal equation becomes
\be\label{eq:eiktilde}\tag{$\widetilde{\text{Eik}}$}
2\Omega\nabla\tilde V\cdot\nabla\tilde\ell_0 + \tilde\ell_0^{\nu_1-1}\Omega\nabla\tilde\ell_0\cdot\nabla\tilde\ell_0 = 0
\ee
and the transport ones become
\be\label{eq:transptilde}\tag{$\tilde T_j$}
\lll \tilde\ell_j + \tilde R^j = 0
\ee
where
\bes
\lll = 2\Omega(\nabla\tilde V + \tilde\ell_0^{\nu_1-1}\nabla\tilde\ell_0)\cdot\nabla + (\nu_1-1)\tilde\ell_0^{\nu_1-2}\Omega\nabla\tilde\ell_0\cdot\nabla\tilde\ell_0,
\ees
with $R^j$ given in \eqref{eq:transp} is a smooth function of the $\partial^\alpha\ell_k$ for $|\alpha|\leq2$ and $k<j$.

In order to solve the eikonal equation, we will assume that $\tilde\ell_0$ admits a formal expansion in powers of $x$, $\tilde\ell_0 \sim \sum_{j\geq0}\tilde\ell_{0,j}$ where $\tilde\ell_{0,j}\in\ppp^j_{hom}$ the set of homogeneous polynomials of degree $j$. Consequently, $\tilde w_0$, which is the left hand side of \eqref{eq:eiktilde}, has a similar development, $\tilde w_0 \sim \sum_{j\geq0}\tilde w_{0,j}$, with $\tilde w_{0,j}\in\ppp^j_{hom}$ satisfying
\be\label{eq:decomposition}
\tilde w_{0,j} = 2\sum_{p+q+r = j+2}\Omega_p\nabla\tilde V_q\cdot\nabla\tilde\ell_{0,r} + \sum_{\sum i_k +p+q+r = j+2}\tilde\ell_{0,i_1}\cdots\tilde\ell_{0,i_{\nu_1-1}}\Omega_p\nabla\tilde\ell_{0,q}\cdot\nabla\tilde\ell_{0,r},
\ee
using the Taylor expansions $\Omega = \sum_{k\geq0}\Omega_k$ and $\tilde V = \sum_{k\geq 0} \tilde V_k$, with $\Omega_k,\tilde V_k \in \ppp^k_{hom}$. We note that $\Omega_0 = \Id$ and $\tilde V_k \neq 0 \iff \exists \nu_i=k$ thus in particular $\tilde V_0 = \tilde V_1 = 0$.

Solving \eqref{eq:eiktilde} leads to solve $\tilde w_{0,j} = 0$ for all $j\in\N$. Using that we set $\ell_0(\s) = 0$, we thus have $\tilde\ell_{0,0} = 0$ which leads to $\tilde w_{0,0} = 0$.

{\it Proof of $i)$.} Then, we start by looking for a linear solution of \eqref{eq:eiktilde}, which means $\tilde\ell_0(x) = \tilde \ell_{0,1}(x) = \xi\cdot x$ for a certain $\xi\in\R^d$ to be determined. Therefore we can write
\bes
\tilde w_0 = 2\Gamma\xi\cdot (x_i^{\nu_i-1})_i + (\xi\cdot x)^{\nu_1-1}|\xi|^2 + O\big(\sum_{i,j=1}^d|x_i||x_j|^{\underline\nu-1}\big),
\ees
where we define $\Gamma = \diag(\nu_it_i)_i$ and we recall $\underline \nu = \inf\nu_i$. Seeing that the first term has separated variables, the second need that too and therefore $\xi$ is colinear to a vector of the canonical basis $(e_i)_{1\leq i\leq d}$ of $\R^d$ (note that when $\nu_1=2$ there is a priori no reason to take such vector, but the result we obtain is still valid in this case). If $\xi = ce_i$ with $i\geq 2$ then
\bes
2\Gamma\xi\cdot (x_i^{\nu_i-1})_i + (\xi\cdot x)^{\nu_1-1}|\xi|^2 = c(2\nu_it_ix_i^{\nu_i-1} + c^{\nu_1}x_i^{\nu_1-1}).
\ees
Even if $\nu_i = \nu_1$, this cannot vanish (because $t_i>0$ and the $\nu_i$ are even), thus we need to take $\xi = ce_1$ which leads to $c = (-2\nu_1t_1)^\frac{1}{\nu_1}$. Plugging $\tilde\ell_0(x) = (2\nu_1|t_1|)^\frac{1}{\nu_1}x_1$ into \eqref{eq:eiktilde} and denoting $\Omega = (\omega_{i,j})_{i,j}$, we obtain
\be\label{eq:w0Case1}
\tilde w_0 = 2c\omega_{1,1}\nu_1t_1x_1^{\nu_1-1} + 2c\sum_{j\geq2}\omega_{1,j}\nu_{j}t_{j}x_{j}^{\nu_{j}-1} + c^{\nu_1+1}\omega_{1,1}x_1^{\nu_1-1} = 2c\sum_{j\geq2}\omega_{1,j}\nu_{j}t_{j}x_{j}^{\nu_{j}-1}
\ee
because we set $c$ so that the first and last term compensate exactly. Considering $\tilde\ell = \tilde\ell_{0}$, noticing that $\omega_{1,j} = O(x)$ for $j\geq2$, and having that $\Delta\ell$ is uniformly bounded around $\s$, we obtain the first result of the lemma.

{\it Discussion on the resolution of \eqref{eq:eiktilde}.} Now let us plug back $\tilde\ell_0 = \sum_{j\geq1}\tilde\ell_{0,j}$, with $\tilde\ell_{0,1}(x) = (2\nu_1|t_1|)^{\frac{1}{\nu_1}}x_1$ in \eqref{eq:eiktilde}. Using \eqref{eq:decomposition} we obtain the following system of equations
\bes
\lll_0 \tilde\ell_{0,j} + R^0_{j+\underline\nu-2} = 0
\ees
for all $j\geq \underline\nu$, where $R^0_{j+\underline\nu-2} \in \ppp^{j+\underline\nu-2}_{hom}$ is a smooth polynomial of the $\tilde\ell_{0,k},\nabla\tilde\ell_{0,k}$ for $k<j$ and $\lll_0$ has the form
\be\label{eq:lll0-1}
\lll_0 = 2\nabla\tilde V_{\nu_1}\cdot\nabla + 2\tilde\ell_{0,1}^{\nu_1-1}\nabla\tilde\ell_{0,1}\cdot\nabla + (\nu_1-1)\tilde\ell_{0,1}^{\nu_1-2}|\nabla\tilde\ell_{0,1}|^2
\ee
if $\underline \nu = \nu_1$ and
\be\label{eq:lll0-2}
\lll_0 = 2\nabla\tilde V_{\underline\nu}\cdot\nabla
\ee
otherwise. When $\underline \nu = \nu_1$, the explicit expression of $\tilde\ell_{0,1}$ implies
\bes
\lll_0 = 2\Lambda(x_i^{\nu_i-1})_i\cdot\nabla + 2\nu_1(\nu_1-1)|t_1|x_1^{\nu_1-2},
\ees
with $\Lambda = \diag(\nu_i|t_i|\delta_{\nu_i,\nu_1})_i$ and $\delta$ the Kronecker symbol. When $\underline \nu < \nu_1$, \eqref{eq:lll0-2} implies $\lll_0(x_1^j) = 0$. In both cases, $\lll_0$ is a linear map from $\ppp^{j}_{hom}$ to $\ppp^{j+\underline\nu-2}_{hom}$. The transport operator $\lll_0 : \ppp^{j}_{hom} \to \ppp^{j+\underline\nu-2}_{hom}$ is thus not injective in all cases. It is actually invertible when $\underline \nu = \nu_1 = 2$ (see Case $iii)$). For the existence of a solution, we are rather concerned with the surjectivity of $\lll_0$ and, in the cases $iv)$ and $v)$, where it is not surjective, we will check that the remainder term $R^0_{j+\underline\nu-2}$ actually belongs to the range of $\lll_0$. In Remark \ref{rem:counterexample}, an example is given for which it is not true and for which the eikonal equation \eqref{eq:eiktilde} cannot be solved with an arbitrary accuracy. All of this shows that without any further assumptions, we cannot solve \eqref{eq:eiktilde} any better with this method.

{\it Proof of $ii)$.} Under Assumption \ref{ass.case1}, $\Omega = \Id$, therefore by \eqref{eq:w0Case1}, we know that taking $\tilde\ell = \tilde\ell_{0,1} = (2\nu_1|t_1|)^{\frac{1}{\nu_1}}x_1$, we have $\tilde w_0 = 0$ ($\Omega = \Id$ leading to $\omega_{1,j} = 0$ for $j>1$). Moreover with that same $\ell$, $\tilde w = \tilde w_0 + h\tilde w_1$ and $\tilde w_1 = -\Delta\ell\circ U$. Using that $U$ is a unitary matrix, $\Delta\ell\circ U = \Delta(\ell\circ U) = \Delta\tilde\ell = 0$ because $\tilde\ell$ is linear.
All this leads to
\bes
\Delta_V(\chi e^{-V/h}) = O(h^\infty)e^{-\big(V+\frac{\ell^{\nu_1}}{\nu_1}\big)/h},
\ees
uniformly with respect to $x$ in a neighborhood of $\s$ because we still have the $r = O(h^\infty)$ uniformly in \eqref{eq:wplusr}.

{\it Proof of $iii)$.} Consider now that Assumption \ref{ass.case4} is satisfied. When $d=1$, any $\s\in\uuu^{(1)}$ is in fact a maximum, it only has a drop direction and no other, therefore the sum in the result of $i)$ is empty and thus the linear function $\tilde\ell_0 = (2\nu|t|)^{\frac1\nu}x$ is not only the solution to the first order but an exact solution of the eikonal equation.

For the transport equations, assume that, for all $j\geq1$, $\tilde\ell_j$ admits a formal expansion in powers of $x$, $\tilde\ell_j \sim \sum_{k\geq0}\tilde\ell_{j,k}$ where $\tilde\ell_{j,k}\in\ppp^k_{hom}$. Consequently, $\tilde w_j$, which is the left hand side of \eqref{eq:transptilde}, has a similar development, $\tilde w_j \sim \sum_{k\geq0}\tilde w_{j,k}$. For $k\geq 0$, $w_{j,k} = 0$ can be written
\bes
\lll_0\tilde\ell_{j,k-\nu+2} + \hat R^j_{k} = 0,
\ees
with $\hat R^j_{k}$ a smooth function of the $\partial^\alpha\tilde\ell_i$ and $\partial^\beta\tilde\ell_{j,r}$ for $i<j$, $|\alpha|\leq2$, $r<k$, and $|\beta|\leq1$. Here we have the following expression for $\lll_0$
\bes
\begin{array}{c|ccc}
    \lll_0:&\ppp^k_{hom}(x)&\to&\ppp^{k+\nu-2}_{hom}(x)\\
    &p&\mapsto&2\nu|t|x^{\nu-1}\cdot\partial_xp + 2\nu(\nu-1)|t|x^{\nu-2}p.
\end{array}
\ees
Since both spaces are one-dimensional and $\lll_0$ is non-zero, it is invertible. But the problem arises beforehand. If $k<\nu-2$, then we only have
\bes
\hat R^j_{k} = 0,
\ees
on which we have no control, and we cannot prove that $\hat R^j_{k}$ is zero for small $k$. A simple way to see that is to note that if $\nu\neq2$, then we have
\bes
w_{1,0} = -\partial_x^2\ell_{0,2},
\ees
which is non-zero in the general case (unless the change of variable $U$ is a unitary matrix). Therefore, we cannot obtain a better estimate than
\bes
w = O(h).
\ees
The conclusion of this, which is also true for higher dimensions, is: to be able to solve the transport equations, we need at least one quadratic direction in $V$.

{\it Proof of $iv)$.} We consider the case where Assumption \ref{ass.case2} is satisfied. Let us recall that solving \eqref{eq:eiktilde} is equivalent to solving
\bes
\lll_0\tilde\ell_{0,j} + R^0_j = 0,
\ees
for all $j\geq2$, where $R^0_j\in\ppp^j_{hom}$ is a smooth polynomial of the $\tilde\ell_{0,k},\nabla\tilde\ell_{0,k}$ for $k<j$. And in this case
\bes
\lll_0 = 4\diag(|t_i|\delta_{\nu_i,2})_ix\cdot\nabla + 4|t_1|.
\ees
One can show that $\lll_0$ injective in this context. Indeed, $\diag(|t_i|\delta_{\nu_i,2})_i$ is a non negative matrix and $4|t_1|$ is a multiplication operator by a positive constant, hence with an argument similar to \cite[Lemma A.1]{BoLePMi22}, $\lll_0$ has positive spectrum, then it is invertible over $\ppp^j_{hom}$ for every $j\geq2$. Therefore we can formally solve \eqref{eq:eikapprox}. Using a Borel procedure, we can transform the formal construction of $\tilde\ell_0$ and thus $\tilde w_0$ into a proper one in $\ccc^\infty(\www)$.

Now for the transport equations, we recall \eqref{eq:transptilde} in this case
\bes
\lll\tilde\ell_j + \tilde R^j = 0
\ees
with
\bes
\lll = 2\Omega(\nabla\tilde V + \tilde\ell_0^{\nu_1-1}\nabla\tilde\ell_0)\cdot\nabla + (\nu_1-1)\tilde\ell_0^{\nu_1-2}\Omega\nabla\tilde\ell_0\cdot\nabla\tilde\ell_0,
\ees
and $\tilde R^j$ smooth, depending only on the $\tilde\ell_k$ for $k<j$. Notice that $\lll$ has a decomposition $\lll = \lll_0 + \lll_>$ where $\lll_>(p) = O(x^{k+1})$ for all $p\in\ppp^k_{hom}$ and all $k\geq1$. Then, just like the eikonal equation, we look for $\tilde\ell_j$ having a formal development in powers of $x$, $\tilde\ell_j \sim \sum_{k\geq0}\tilde\ell_{j,k}$, $\tilde\ell_{j,k}\in\ppp^k_{hom}$ which leads to $\tilde w_j \sim \sum_{k\geq0}\tilde w_{j,k}$, $\tilde w_{j,k}\in\ppp^k_{hom}$ and solving the transport equations \eqref{eq:transptilde} is formally equivalent to solving
\bes
\lll_0\tilde\ell_{j,k} + \hat R^j_k = 0
\ees
for all $k\in\N$, where $\hat R^j_k$ a smooth function of the $\partial^\alpha\tilde\ell_i$ and $\partial^\beta\tilde\ell_{j,r}$ for $i<j$, $|\alpha|\leq2$, $r<k$ and $|\beta|\leq1$. By the work done on \eqref{eq:eiktilde}, we know we can find a family of smooth functions $(\tilde\ell_{j,k})_{k\in\N}$ that solves \eqref{eq:transpapprox}. Using another Borel procedure we transform this formal construction of $\tilde\ell_j$ into a $\ccc^\infty(\www)$ one.

The final step is to obtain $\tilde w = O(x^\infty)$. We use a third Borel procedure, with respect to $h$ this time and not $x$, to give a proper definition of $\sum_{j\geq0}\ell_jh^j$. This gives the existence of an $\ell\in\ccc^\infty(\www)$ such that $\ell\sim\sum_{j\geq0}\ell_jh^j$ in the sense of \eqref{eq:classicexpansion}. This conclude the proof of $iv)$.

{\it On the invertibility of $\lll_0$.} In a general setting, let $X = (x_i)_{i\in\lb1,p\rb}$, $Y = (y_i)_{i\in\lb1,q\rb}$ two sets of variables and $Q\in\ppp^j_{hom}(X,Y)$. By taking aside the monomials of $Q$ depending only on $X$ or only on $Y$, we have a decomposition
\be\label{eq:decomp}
Q = \sum_iP_iQ_i + Q_X + Q_Y
\ee
with $Q_X\in\ppp^j_{hom}(X)$, $Q_Y\in\ppp^j_{hom}(Y)$ and the $(P_i)_i\subset\ppp_{hom}(X),(Q_i)_i\subset\ppp_{hom}(Y)$ are some families of homogeneous polynomials such that $\deg P_i + \deg Q_i = j$ and $\deg P_i,\deg Q_i < j$. Note that $Q_X$ and $Q_Y$ are unique in such a decomposition.

Now consider $\lll$ a linear operator such that there exists $k\in\N$, for all $j\in\N^*$, $\lll:\ppp^j_{hom}(Y)\to\ppp^{j+k}_{hom}(Y)$. Thus for all $P\in\ppp^j_{hom}(X)$, $[\lll,P] = 0$, hence we easily see the action of $\lll$ on $Q$
\bes
\lll Q = \sum_iP_i\lll Q_i + \lll Q_Y.
\ees
Usually, we are led to consider injective operators, and thus for all $j$
\be\label{eq:invertible}
\lll \mbox{ is invertible }\iff\dim\ppp^j_{hom}(Y) = \dim\ppp^{j+k}_{hom}(Y).
\ee
Provided that the left-hand side of \eqref{eq:invertible} is true, the resolution of the equation $\lll Q' = Q$ in $\ppp^{j}_{hom}(X,Y)$ is equivalent to showing that $Q_X = 0$ in its decomposition \eqref{eq:decomp}.

But we know that
\bes
\forall j>0\ \dim\ppp^j_{hom}(Y) = \dim\ppp^{j+k}_{hom}(Y) \iff k=0 \mbox{ or } q=\sharp Y = 1.
\ees
We investigate the case $k=0$ in $iv)$ and $q=1$ in $v)$. Notice that $iv)$ treats the case $k=0$ with the additional restriction that $p=\sharp X = 0$.

{\it Proof of $v)$.} Let us now consider a potential satisfying Assumption \ref{ass.case3}, then we can assume we have the directions on the following order
\bes
\forall i\geq2\ t_1<0<t_i
\ees
(which we already assumed) and
\bes
\exists J\in\lb2,d\rb,\ \nu_i = 2 \iff i\in\lb2,J\rb,
\ees
let us recall that $\hat\nu = \min_{i\notin\lb2,J\rb}\nu_i$ is the second lowest degeneracy order.

Moreover, we know that the eikonal equation split up into
\bes
\lll_0\tilde\ell_{0,j} + R^0_j = 0
\ees
with $\lll_0 = 4\diag(t_i\delta_{\nu_i,2})_ix\cdot\nabla$, and $R^0_j$ is a smooth polynomial which depends on the $\tilde\ell_{0,k}$ for $k<j$. To have a better view on what is inside $R^0_j$, from \eqref{eq:eiktilde}, we recall \eqref{eq:decomposition}
\bes
\tilde w_{0,j} = 2\sum_{p+q+r = j+2}\Omega_p\nabla\tilde V_q\cdot\nabla\tilde\ell_{0,r} + \sum_{\sum i_k +p+q+r = j+2}\tilde\ell_{0,i_1}\cdots\tilde\ell_{0,i_{\nu_1-1}}\Omega_p\nabla\tilde\ell_{0,q}\cdot\nabla\tilde\ell_{0,r}.
\ees
Let us describe more precisely what $\tilde w_{0,j}$ is made of for small $j$.

Because $\nabla\tilde\ell_{0,1}\bot\nabla\tilde V_2$ since $\nabla\tilde\ell_{0,1}$ only has one non-zero term on the first coordinate thanks to $i)$, for $j\leq1$, $\tilde w_{0,j} = 0$. Then, noticing that the principal term in the second sum in \eqref{eq:decomposition} is $\tilde\ell_{0,1}^{\nu_1-1}\Omega_0\nabla\tilde\ell_{0,1}\cdot\nabla\tilde\ell_{0,1}$ which is of order $\nu_1-1$, for $j < \hat\nu - 1$, this sum is null.

Moreover, for all $p\in\N$,
\be\label{eq:blabla}
2\Omega_p(\nu_1t_1x_1^{\nu_1-1},0,\ldots,0)\cdot\nabla\tilde\ell_{0,1} + \tilde\ell_{0,1}^{\nu_1-1}\Omega_p\nabla\tilde\ell_{0,1}\cdot\nabla\tilde\ell_{0,1} = 0.
\ee

Thus for $j = \hat\nu-1$, apart from the terms involving $\nabla\tilde V_2$ we only have $\Omega_0\nabla\tilde V_{\hat\nu}\cdot\nabla\tilde\ell_{0,1}$ and possibly $\tilde\ell_{0,1}^{\nu_1-1}\Omega_0\nabla\tilde\ell_{0,1}\cdot\nabla\tilde\ell_{0,1}$ if $\nu_1 = \hat\nu$. In this case, thanks to \eqref{eq:blabla}, there only remains $\big(\nabla\tilde V_{\nu_1} - (\nu_1t_1x_1^{\nu_1-1},0,\ldots,0)\big)\cdot\nabla\tilde\ell_{0,1}$ which again is $0$ because those two vectors are orthogonal. And if $\nu_1\neq\hat\nu$ then $\nabla\tilde V_{\hat\nu}\bot\nabla\tilde\ell_{0,1}$, all this leads to
\bes
\forall j\in\lb2,\hat\nu-1\rb,\ \tilde w_{0,j} = 2\nabla\tilde V_2\cdot\sum_{k=1}^{j}\Omega_{j-k}\nabla\tilde\ell_{0,k}
\ees
using that the matrices $\Omega_k$ are symmetric because $\Omega$ is. Having that $\lll_0\tilde\ell_{0,j} = 2\nabla\tilde V_2\cdot\Omega_0\nabla\tilde\ell_{0,j}$ which is the last term of the previous sum, this therefore means that
\be\label{eq:R0jdecomp}
\forall j\in\lb2,\hat\nu-1\rb,\ R^0_j = 2\nabla\tilde V_2\cdot\sum_{k=1}^{j-1}\Omega_{j-k}\nabla\tilde\ell_{0,k}.
\ee
Now we want to make sure that $R^0_j\in\Im\lll_0$ so that we can obtain $\tilde w_{0,j} = 0$ for all $j\leq\hat\nu-1$.

Using the preliminaries about $\ppp^j_{hom}(X,Y)$ stated above applied to
\bes
\lll_0 = 4\diag(t_i\delta_{\nu_i,2})_ix\cdot\nabla
\ees
which is indeed injective on $\ppp^j_{hom}(Y)$, $Y = (x_i)_{i\in\lb2,J\rb}$ (see \cite[Lemma A.1]{BoLePMi22}) and having that $\lll_0$ is an endomorphism, there only remains to show that the $Q_X$ associated with $R^0_j$ in \eqref{eq:decomp} is zero. From \eqref{eq:R0jdecomp}, notice that we can write $R^0_j = Y\cdot (P_2,\ldots,P_J)$ for some polynomials $P_i\in\ppp^{j-1}_{hom}(X,Y)$, which gives the result.

Moreover, we see that we can choose by induction a $\tilde\ell_{0,j}$ solving $\lll_0\tilde\ell_{0,j} + R^0_j = 0$ such that $Q_X = 0$ in its decomposition \eqref{eq:decomp}, therefore
\be\label{eq:boundell0j}
\forall j\in\lb2,\hat\nu-1\rb,\ \ \tilde\ell_{0,j} = O(|x|^{j-1}\sum_{i=2}^J|x_i|).
\ee

The last thing to do to prove $iv)$, is to obtain estimates on $\tilde w_{0,j}$ for $j\geq\hat\nu$. Looking at \eqref{eq:decomposition}, we have the following terms appearing
\begin{itemize}
    \item $\tilde\ell_{0,i_1}\cdots\tilde\ell_{0,i_{\nu_1-1}}\Omega_p\nabla\tilde\ell_{0,q}\cdot\nabla\tilde\ell_{0,r}$ with at least one $i_k\geq2$, let us denote $n$ the number of such $i_k$. Therefore, using \eqref{eq:boundell0j}
    \bes
    \tilde\ell_{0,i_1}\cdots\tilde\ell_{0,i_{\nu_1-1}}\Omega_p\nabla\tilde\ell_{0,q}\cdot \nabla\tilde\ell_{0,r} = O\bigg(\Big(|x|\sum_{i=2}^J|x_i|\Big)^n|x_1|^{\nu_1-n}|x|^{j-(\nu_1-n)-2n}\bigg).
    \ees
    Using that necessarily, for such a term to appear, we need to have $j\geq\nu_1$, then we obtain a term that is $O\bigg(\Big(\sum_{i=2}^J|x_i|\Big)^n|x_1|^{\nu_1-n}\bigg)$. Having that $n\in\lb1,\nu_1\rb$, we can estimate this term as follows for $|x|\leq1$
    \bes
    \begin{split}
        \Big(\sum_{i=2}^J|x_i|\Big)^n|x_1|^{\nu_1-n} &\lesssim \Big(\sum_{i=2}^J|x_i|\Big)^{\nu_1} + \sum_{i=2}^J|x_i||x_1|^{\nu_1-1}\\
        &\lesssim \sum_{i=2}^J|x_i|^{\nu_1} + |x||x_1|^{\nu_1-1}\\
        &\lesssim |x|^{\nu_1-\underline\nu+1}\sum_{i=2}^J|x_i|^{\nu_i-1} + |x|\sum_{i\notin\lb2,J\rb}|x_i|^{\nu_i-1}\\
        &\lesssim |x|^{\hat\nu-\underline\nu+1}\sum_{i=2}^J|x_i|^{\nu_i-1} + |x|\sum_{i\notin\lb2,J\rb}|x_i|^{\nu_i-1}\\
    \end{split}
    \ees
    where we used that for $i\in\lb2,J\rb$, we have $\nu_i=\underline\nu$ and also that $\hat\nu\leq\nu_1$;
    \item $\tilde\ell_{0,1}^{\nu_1-1}\Omega_p\nabla\tilde\ell_{0,q}\cdot\nabla\tilde\ell_{0,r}$ with $(q,r)\neq(1,1)$ which leads to a $O(|x_1|^{\nu_1-1}|x|)$;
    \item $2\Omega_{j-\nu_1+1}\nabla\tilde V_{\nu_1}\cdot\tilde\ell_{0,1} + \tilde\ell_{0,1}^{\nu_1-1}\Omega_{j-\nu_1+1}\nabla\tilde\ell_{0,1}\cdot\nabla\tilde\ell_{0,1}$ which is equal to
    \bes
    2\Omega_{j-\nu_1+1}\big(\nabla\tilde V_{\nu_1} - (\nu_1t_1x_1^{\nu_1-1},0,\ldots,0)\big)\cdot\nabla\tilde\ell_{0,1}
    \ees
    thanks to \eqref{eq:blabla}. But if $j-\nu_1+1 = 0$, then the whole term is $0$, thus it is $O(|x|\sum_{i\notin\lb2,J\rb}|x_i|^{\nu_i-1})$;
    \item $2\Omega_p\nabla\tilde V_{q}\cdot\nabla\tilde\ell_{0,r}$ with $q\geq\hat\nu$ and $(q,r)\neq (\nu_1,1)$, hence if $(p,r) = (0,1)$ then it is $0$ and else it is $O(|x|\sum_{\nu_i\geq q}|x_i|^{\nu_i-1}) = O(|x|\sum_{i\notin\lb2,J\rb}|x_i|^{\nu_i-1})$;
    \item $2\nabla\tilde V_2\cdot\sum_{k=1}^{\hat\nu-1}\Omega_{j-k}\nabla\tilde\ell_{0,k}$, which is $O(|x|^{\hat\nu-\underline\nu+1}\sum_{i=2}^J|x_i|^{\nu_i-1})$. 
\end{itemize}
Even though $\nu_j-1 = 1$, we keep it that way so that this discussion can easily be generalized for the next case. All this together we obtain that with $\tilde\ell = \tilde\ell_0 = \sum_{j=1}^{\hat\nu-1}\tilde\ell_{0,j}$,
\bes
\tilde w = O\bigg(|x|\sum_{i\notin\lb2,J\rb}|x_i|^{\nu_i-1} + |x|^{\hat\nu-\underline\nu+1}\sum_{i\in\lb2,J\rb}|x_i|^{\nu_i-1} + h\bigg).
\ees

{\it Proof of $vi)$.} Under Assumption \ref{ass.case5}, we have $\exists!i,\ \nu_i = \underline\nu$ and $i\neq1$. Without loss of generality, let us assume that the direction corresponding to the minimal exponent is the last one, namely
\bes
\forall i < d,\ \nu_i > \nu_d.
\ees
Moreover, in that case, from \eqref{eq:lll0-2}, we know that $\lll_0$ has the form
\bes
\lll_0 = 2\nu_dt_dx_d^{\nu_d-1}\partial_d,
\ees
which is indeed invertible from $\ppp^j_{hom}(x_d)$ to $\ppp^{j+\nu_d-2}_{hom}(x_d)$ for all $j\geq1$.

Following the proof of $v)$, this allows us to solve $\tilde w_{0,k} = 0$ for all $k < \hat\nu$, recalling $\hat\nu = \inf_{i<d}\nu_i$ is this case. Following the same discussion on the terms appearing in \eqref{eq:decomposition}, this leads to an estimate similar to $iv)$
\bes
\tilde w = O\big(|x|\sum_{j<d}|x_j|^{\nu_j-1} + |x|^{\hat\nu-\underline\nu+1}|x_d|^{\nu_d-1} + h\big).
\ees

{\it Proof of $vii)$.} Let us notice that in cases $i)$, $ii)$, and $iii)$, $\tilde\ell_0(x) = \tilde\ell_{0,1}(x) = (2\nu_1|t_1|)^{\frac1{\nu_1}}x_1$ which gives the result with no negligible terms. For case $iv)$, it is sufficient to prove the following preliminary result: for all $n\geq2$,
\be\label{eq:borneR0}\tag{$H_n$}
\tilde\ell_{0,n} = O(|x|\sum_{j=1}^d|x_j|^{\nu_j-1}).
\ee
Indeed, since $\nu_1=2$ in this setting, this implies that
\bes
\begin{split}
    \tilde V + \frac{\tilde\ell_0^2}{2} &= V(\s) + \sum_{i=1}^d t_ix_i^{\nu_i} + \frac{\tilde\ell_{0,1}^2}{2} + \tilde\ell_{0,1}O(|x|\sum_{j=1}^d|x_j|^{\nu_j-1}) + O(|x|\sum_{j=1}^d|x_j|^{\nu_j-1})^2\\
    &= V(\s) + \sum_{i=1}^d |t_i|x_i^{\nu_i} + O\big(|x_1||x|\sum_{j=1}^d|x_j|^{\nu_j-1} + |x|^2\sum_{j,k = 1}^d|x_j|^{\nu_j-1}|x_k|^{\nu_k-1}\big)\\
    &= V(\s) + \sum_{i=1}^d |t_i|x_i^{\nu_i} + O\big(|x|x_1^2 + |x|\sum_{j=1}^dx_j^{2\nu_j-2} + |x|^2\sum_{j=1}^dx_j^{2\nu_j-2}\big)\\
    &= V(\s) + \sum_{i=1}^d |t_i|x_i^{\nu_i} + O(|x|\sum_{j=1}^d|x_j|^{\nu_j}),
\end{split}
\ees
using the Young inequality and $\nu_j\geq2\Rightarrow2\nu_j-2\geq\nu_j$.

We now proceed to prove \eqref{eq:borneR0}. We recall the form of the eikonal equation \eqref{eq:eikapprox}
\bes
2\Omega\nabla\tilde V\cdot\nabla\tilde\ell_0 + \tilde\ell_0\Omega\nabla\tilde\ell_0\cdot\nabla\tilde\ell_0 = O(x^\infty).
\ees
From the Taylor expansion of $\tilde w_0 = 2\Omega\nabla\tilde V\cdot\nabla\tilde\ell_0 + \tilde\ell_0\Omega\nabla\tilde\ell_0\cdot\nabla\tilde\ell_0$, we recall the decomposition \eqref{eq:decomposition} of $\tilde w_{0,k}$ for all $k\in\N$
\be\label{eq:expanw0k}
\tilde w_{0,k} = 2\sum_{p+q+r = k+2}\Omega_p\nabla\tilde V_q\cdot\nabla\tilde\ell_{0,r} + \sum_{p+q+r+s = k+2}\tilde\ell_{0,p}\Omega_q\nabla\tilde\ell_{0,r}\cdot\nabla\tilde\ell_{0,s}.
\ee
We also recall that $\Omega_0 = \Id$ and $\tilde V_i\neq0\iff\exists j,\ \nu_j=i$, and since $\s\in\vvv^{(1)}$, all the $\nu_i$ are even. A direct consequence is that $\tilde V_{2i+1} = 0$ for all $i$, and we also have $\tilde V_0 = \tilde V_1 = 0$.

Now let us prove \eqref{eq:borneR0} by induction. First, recalling that $\lll_0 = 4\diag(|t_i|\delta_{\nu_i,2})_ix\cdot\nabla + 4|t_1|$, we see that $\lll_0$ is a diagonal operator on $\ppp^k_{hom}$, namely for all $\alpha\in\N^d$, $|\alpha|=k$, $\lll_0x^\alpha = c_\alpha x^\alpha$ for some $c_\alpha\in\R$. Hence using that for all $k\in\N$, $\tilde w_{0,k} = \lll_0\tilde\ell_{0,k} + \tilde R^0_k = 0$, we see that if for $n\geq2$, $\tilde R^0_n = O(|x|\sum_{j=1}^d|x_j|^{\nu_j-1})$ then \eqref{eq:borneR0} is true. But making $\lll_0\tilde\ell_{0,k}$ and $\tilde R^0_k$ appear in \eqref{eq:expanw0k}, we have
\bes
\lll_0\tilde\ell_{0,k} = 2\Omega_0\nabla\tilde V_2\cdot\nabla\tilde\ell_{0,k} + \tilde\ell_{0,k}\Omega_0\nabla\tilde\ell_{0,1}\cdot\nabla\tilde\ell_{0,1} + 2\tilde\ell_{0,1}\Omega_0\nabla\tilde\ell_{0,1}\cdot\nabla\tilde\ell_{0,k}
\ees
and hence
\be\label{eq:R0ksums}
\tilde R^0_k = 2\sum_{E^1_k}\Omega_p\nabla\tilde V_q\cdot\nabla\tilde\ell_{0,r} + \sum_{E^2_k}\tilde\ell_{0,p}\Omega_q\nabla\tilde\ell_{0,r}\cdot\nabla\tilde\ell_{0,s}
\ee
where
\bes
E^1_k = \{p,q,r\in\N, p+q+r = k+2\}\setminus\{(0,2,k)\}
\ees
and
\bes
E^2_k = \{p,q,r,s\in\N, p+q+r+s = k+2\}\setminus\{(k,0,1,1),(1,0,k,1),(1,0,1,k)\}.
\ees

Let us then prove \eqref{eq:borneR0} by induction. For $n=2$, we have
\bes
R^0_2 = \underbrace{2\nabla\tilde V_3\cdot\nabla\tilde\ell_{0,1}}_{=0} + 2\Omega_1\nabla\tilde V_2\cdot\nabla\tilde\ell_{0,1} + \tilde\ell_{0,1}\Omega_1\nabla\tilde\ell_{0,1}\cdot\nabla\tilde\ell_{0,1},
\ees
recalling $\tilde\ell_{0,1} = (2\nu_1|t_1|)^{\frac1{\nu_1}}x_1$, we therefore see that
\bes
R^0_2 = O(|x|\sum_{j=1}^d|x_j|^{\nu_j-1}) + O(|x_1||x|) = O(|x|\sum_{j=1}^d|x_j|^{\nu_j-1}),
\ees
thus so is $\tilde\ell_{0,2}$.

Now assume \eqref{eq:borneR0} is true for all $n<k$, for some $k\geq3$. If we study each term appearing in the sums in \eqref{eq:R0ksums}, we identify four distinct types of terms to analyze.
\begin{itemize}
    \item $\Omega_p\nabla\tilde V_q\cdot\nabla\tilde\ell_{0,r}$, with $p=0$ and $r=1$, hence $q>2$ and then $\nabla\tilde V_q \bot \nabla\tilde\ell_{0,1}$, thus $\Omega_p\nabla\tilde V_q\cdot\nabla\tilde\ell_{0,r} = 0$ because $\Omega_0 = \Id$;
    \item $\Omega_p\nabla\tilde V_q\cdot\nabla\tilde\ell_{0,r}$, with $p\geq1$ or $2\leq r<k$, hence $\Omega_p\nabla\tilde V_q\cdot\nabla\tilde\ell_{0,r} = O(|x|\sum_{j=1}^d|x_j|^{\nu_j-1})$
    \item[] because $\tilde\ell_{0,r}\in\ppp^r_{hom}\Rightarrow\nabla\tilde\ell_{0,r} = O(x^{r-1}) = O(x)$ for $r\geq2$;
    \item $\tilde\ell_{0,p}\Omega_q\nabla\tilde\ell_{0,r}\cdot\nabla\tilde\ell_{0,s}$, with $2\leq p<k$, since \eqref{eq:borneR0} holds for $n=p$, we deduce that
    \item[]$\tilde\ell_{0,p}\Omega_q\nabla\tilde\ell_{0,r}\cdot\nabla\tilde\ell_{0,s} = O(|x|\sum_{j=1}^d|x_j|^{\nu_j-1})$;
    \item $\tilde\ell_{0,p}\Omega_q\nabla\tilde\ell_{0,r}\cdot\nabla\tilde\ell_{0,s}$, with $p=1$, hence $q\geq1$ or $2\leq r<k$ or $2\leq s<k$, in any case this leads to $\Omega_q\nabla\tilde\ell_{0,r}\cdot\nabla\tilde\ell_{0,s} = O(x)$ and thus
    \bes
    \tilde\ell_{0,p}\Omega_q\nabla\tilde\ell_{0,r}\cdot\nabla\tilde\ell_{0,s} = O(x_1|x|) = O(|x|\sum_{j=1}^d|x_j|^{\nu_j-1}).
    \ees
\end{itemize}

For case $v)$, thanks to \eqref{eq:boundell0j}, for $k\in\lb2,\hat\nu-1\rb$, $\tilde\ell_{0,k} = O\big(|x|\sum_{j\in\lb2,J\rb}|x_j|\big)$ and thus
\bes
\tilde\ell_0^{\nu_1} = \tilde\ell_{0,1}^{\nu_1} + \tilde\ell_{0,1}^{\nu_1-1}O\big(|x|\sum_{j\in\lb2,J\rb}|x_j|\big) + O\big(|x|^2\sum_{j\in\lb2,J\rb}|x_j|^2\big).
\ees
Recalling that for $j\in\lb2,J\rb$, $\nu_j = 2$, the last rest is of the form $O\big(|x|\sum_{j=1}^d|x_j|^{\nu_j})$ so we only need to study the middle term. By homogeneity and for $|x|\leq1$,
\bes
\forall j\in\lb2,J\rb,\ \ |x_j| |x_1|^{\nu_1-1} \lesssim |x_j|^{\nu_1} + |x_1|^{\nu_1} \leq |x_j|^{\nu_j} + |x_1|^{\nu_1},
\ees
and therefore $\tilde\ell_0^{\nu_1} = \tilde\ell_{0,1}^{\nu_1} + O\big(|x|\sum_{j=1}^d|x_j|^{\nu_j})$. Knowing that $\tilde V + \tilde\ell_{0,1}^{\nu_1} = V(\s) + \sum_{j=1}^d|t_j|x_j^{\nu_j}$, this gives the result for this case, it remains to study $\ell_0$ constructed in $vi)$.

In this case, using that $\tilde\ell_{0,k} = O(|x|^{k-1}|x_d|) = O(|x||x_d|)$ for $k\geq2$, we have
\bes
\tilde\ell_0^{\nu_1} = \tilde\ell_{0,1}^{\nu_1} + O(\sum_{k=1}^{\nu_1}|x|^{k}|x_d|^{k}|x_1|^{\nu_1-k}).
\ees
Recalling that under Assumption \ref{ass.case5}, $\nu_1 = \overline\nu$, we obtain for $k\geq1$, $|x|\leq1$
\bes
|x|^k|x_1|^{\nu_1-k}|x_d|^k \lesssim |x|^{\nu_1}|x_d| \lesssim |x_d|\sum_{j=1}^d|x_j|^{\nu_j},
\ees
hence $\tilde\ell_0^{\nu_1} = \tilde\ell_{0,1}^{\nu_1} + O(|x|\sum_{j=1}^d|x_j|^{\nu_j})$.

\ep

\begin{remark}\label{rem:signell}
    Note that $-\ell$ solves the equations the same way $\ell$ does. For now, the sign does not matter, but we will fix it in the next section when properly constructing the quasimodes.
\end{remark}

\begin{remark}\label{rem:counterexample}
    If $\nu_1\neq2$, there are cases when $R^0_{j}\notin\Im(\lll_0)$. Therefore we still cannot go further in the development and solving of $\tilde w = 0$.
\end{remark}
Looking at the proof of $v)$ in Lemma \ref{lem:ell}, it is evident that a term of the form $\Omega_1\nabla V_{\hat\nu}\cdot\nabla\tilde\ell_{0,1}$ can appear in $R^0_{\hat\nu}$ and contains, for instance, $x_i^{\hat\nu}$, with $\nu_i=\hat\nu$ which clearly cannot be generated by the operator $\lll_0$. Therefore, consider a potential where $\nu_1 = \underline \nu > 2$. For example $V$ such that $\tilde V = -\frac{x^4}{4} + \frac{y^4}{4}$ around $\s=0$, with the coordinate transformation $U:(x,y)\mapsto(x,y(1+x))$. Solving $\tilde w_0 = 0$ to the principal order leads us to $\tilde\ell_{0,1} = 2^\frac14x$, and in this context, we observe that $\tilde w_{0,0} = \tilde w_{0,1} = \tilde w_{0,2} = 0$. Then, denoting $\Omega = \Id + \Omega_1 + O(\vvert{(x,y)}^2)$, with $\Omega_1$ linear in $(x,y)$, the next order is
\bes
\begin{split}
    \tilde w_{0,4} &= 2\begin{pmatrix}-x^3\\y^3\end{pmatrix}\cdot\nabla\tilde\ell_{0,2} + 2\tilde\ell_{0,1}^3\nabla\tilde\ell_{0,1}\cdot\nabla\tilde\ell_{0,2} + 3\tilde\ell_{0,2}\tilde\ell_{0,1}^2|\nabla\tilde\ell_{0,1}|^2\\
    &\phantom{******}+2\Omega_1\nabla\tilde V\cdot\nabla\tilde\ell_{0,1} + \tilde\ell_{0,1}^3\Omega_1\nabla\tilde\ell_{0,1}\cdot\nabla\tilde\ell_{0,1}.
\end{split}
\ees
Since $\tilde w_{0,3} = 2\nabla\tilde V\cdot\nabla\tilde\ell_{0,1} + \tilde\ell_{0,1}^3|\nabla\tilde\ell_{0,1}|^2 = 0$, we can simplify the second line: Denoting $\Omega_1=\begin{pmatrix}\omega^1_{11}&\omega^1_{12}\\\omega^1_{21}&\omega^1_{22}\end{pmatrix}$, we have
\bes
\begin{split}
    \Omega_1\nabla\tilde V\cdot\nabla\tilde\ell_{0,1} &= \omega^1_{11}\nabla\tilde V\cdot\nabla\tilde\ell_{0,1} + 2^\frac14\omega^1_{12}y^3,\\
    \Omega_1\nabla\tilde\ell_{0,1}\cdot\nabla\tilde\ell_{0,1} &= \omega^1_{11}|\nabla\tilde\ell_{0,1}|^2
\end{split}
\ees
hence
\bes
\tilde w_{0,4} = 2\begin{pmatrix}x^3\\y^3\end{pmatrix}\cdot\nabla\tilde\ell_{0,2} + 6x^2\tilde\ell_{0,2} + 2^\frac54\omega^1_{12}y^3.
\ees

Recalling that $U(x,y)=(x,y(1+x))$,
\bes
(d_{(x,y)}U^Td_{(x,y)}U)^{-1} = \begin{pmatrix}1&-\frac y{1+x}\\-\frac y{1+x}&\frac {1+y^2}{(1+x)^2}\end{pmatrix},
\ees
hence
\bes
\Omega_1 = \begin{pmatrix}0&-y\\-y&-2x\end{pmatrix}
\ees
(This confirms that $d_0U$ is unitary), therefore $\omega^1_{12}=-y$. To solve $\tilde  w_{0,4} = 0$, we must consider $\tilde\ell_{0,2}$ homogeneous of order $2$. First, we can easily see that we must eliminate the terms in $x^4$ and $x^3$, which forces $\tilde\ell_{0,2}$ to lack any $x$-dependence. Thus, taking $\tilde\ell_{0,2} = ay^2$ leads to
\bes
(4a-2^\frac54)y^4+6ax^2y^2 = 0
\ees
which admits no solution for $a\in\R$, leading to a contradiction.

\section{Geometric constructions}\label{sec:geo}

To construct proper quasimodes, we rely on the notions introduced in Definition \ref{def:labeling} as well as the subsequent labeling. We also assume that the generic condition \eqref{ass.gener} holds, which we recall below:
\bes\tag{Gener}
\begin{array}{l}
    (\ast)\mbox{ for any }\m\in\uuu^{(0)},\m \mbox{ is the unique global minimum of } V_{|E(\m)},\\
    (\ast)\mbox{ for all }\m\neq\m'\in\uuu^{(0)},\j(\m)\cap\j(\m')=\emptyset.
\end{array}
\ees
Given $\m\in\uuu^{(0)}\setminus\{\um\}$, we have $\bsigma(\m) = \sigma_i$ for some $i\geq2$. Since $\sigma_{i-1} > \sigma_i$, there exists a unique connected component of $\{V<\sigma_{i-1}\}$ containing $\m$, which we denote by $E_-(\m)$. Moreover, we define the following map:
\be\label{eq:defik}
\begin{array}{c|ccc}
     \k:&\uuu&\to&(0,\frac d4]\\
     &x^*&\mapsto&\sum_{i=1}^d\frac1{2\nu_i^{x^*}}.
\end{array}
\ee
This function characterizes the degeneracy at the point $x^*$, and it will appear when rescaling $V$ around that point. We also define a partition of the set $\j(\m)$,
\be\label{eq:defjautre}
\j(\m) = \j_\infty(\m)\sqcup\j_{1}(\m)\sqcup\j_{2}(\m),
\ee
where, for $\s\in\j(\m)$, we assign:
\begin{itemize}
    \item $\s\in\j_\infty(\m)$ if $\s$ satisfies Assumption \ref{ass.case1} or \ref{ass.case2};
    \item $\s\in\j_{1}(\m)$ if $\s\notin\j_{\infty}(\m)$ but $\s$ satisfies Assumption \ref{ass.case3} or \ref{ass.case5}; and
    \item $\s\in\j_{2}(\m)$ contains all remaining saddle points of $\j(\m)$.
\end{itemize}

Notice that in the Morse case, for all $\m\in\uuu^{(0)}$, we have $\j(\m) = \j_\infty(\m)$ and $\k\equiv\frac d4$.

We then follow the construction as described in \cite[Section 4]{BoLePMi22}.

\begin{figure}[h]
\centering
\tikzset{every picture/.style={line width=0.75pt}}
\begin{tikzpicture}[x=0.75pt,y=0.75pt,yscale=-1,xscale=1]
\draw    (47,110) .. controls (292,51) and (275,353) .. (518,298) ;
\draw    (73,320) .. controls (262,311) and (329,74) .. (504,72) ;
\draw [dashed]   (46.7,76) .. controls (191.37,54.1) and (260.37,168.1) .. (292.37,161.1) .. controls (324.37,154.1) and (399.37,45.1) .. (493.37,47.1) ;
\draw [dashed]   (82.37,338.1) .. controls (227.03,316.2) and (244.37,240.1) .. (283.37,239.1) .. controls (322.37,238.1) and (281.37,322.1) .. (514.37,320.1) ;
\draw [dotted]   (309.37,43.1) -- (261.37,343.1) ;
\draw [dashed]   (327.8,137.7) -- (308.37,259.1) ;
\draw [dashed]   (257.8,145.7) -- (238.37,267.1) ;
\draw    (50.1,145.7) -- (72.42,83.58) ; \draw [shift={(73.1,81.7)}, rotate = 109.77] [color={rgb, 255:red, 0; green, 0; blue, 0 }  ][line width=0.75]    (10.93,-3.29) .. controls (6.95,-1.4) and (3.31,-0.3) .. (0,0) .. controls (3.31,0.3) and (6.95,1.4) .. (10.93,3.29)   ;
\draw    (506.1,125.7) -- (463.26,65.33) ; \draw [shift={(462.1,63.7)}, rotate = 54.64] [color={rgb, 255:red, 0; green, 0; blue, 0 }  ][line width=0.75]    (10.93,-3.29) .. controls (6.95,-1.4) and (3.31,-0.3) .. (0,0) .. controls (3.31,0.3) and (6.95,1.4) .. (10.93,3.29)   ;
\draw    (224,311) -- (259.13,237.8) ; \draw [shift={(260,236)}, rotate = 115.64] [color={rgb, 255:red, 0; green, 0; blue, 0 }  ][line width=0.75]    (10.93,-3.29) .. controls (6.95,-1.4) and (3.31,-0.3) .. (0,0) .. controls (3.31,0.3) and (6.95,1.4) .. (10.93,3.29)   ;
\draw (40,155) node{$E_{\m,\tau,\delta}^+$};
\draw (520,140) node{$E_{\m,\tau,\delta}^-$};
\draw (110,200) node{$\times$} node[below right]{$\m$};
\draw (295,202) node{$\s$};
\draw (350,180) node{$\times$} node[below right]{$\m'$};
\draw (220,330) node{$\ccc_{\s,\tau,\delta}$};
\draw (320,35) node{$\s+H$};
\draw (140,155) node{$E(\m)$};
\draw (100,302) node{$\partial E(\m)$};
\end{tikzpicture}
\caption{Representation of the potential $V$ near a point $\s\in\j(\m)$. Here $H$ denotes the hypersurface defined by $\s + H = \{U^{-1}(x)_1 = \s_1\}$}.
\label{fig:levelsets}
\end{figure}

Recall that $\um$ is the unique global minimum of $V$ (uniqueness follows from \eqref{ass.gener}). Now consider an arbitrary $\m\in\uuu^{(0)}\setminus\{\um\}$. For each $\s\in\j(\m)$ and any $\tau,\delta>0$, we define the sets $\bbb_{\s,\tau,\delta}$, $\ccc_{\s,\tau,\delta}$ and $E_{\m,\tau,\delta}$ by
\bes
\bbb_{\s,\tau,\delta} = \{V\leq V(\s)+\delta\}\cap\{x\in\R^d,(2\nu_1^\s|t_1^\s|)^{\frac1{\nu_1^\s}}|U^{-1}(x)_1 - \s_1|\leq\tau\},
\ees
\be\label{eq:ccc}
\ccc_{\s,\tau,\delta} \mbox{ the connected component of } \bbb_{\s,\tau,\delta} \mbox{ containing } \s
\ee
and
\bes
E_{\m,\tau,\delta} = \big(E_-(\m)\cap\{V<V(\j(\m))+\delta\}\big)\setminus\bigcup_{\s\in\j(\m)}\ccc_{\s,\tau,\delta}.
\ees
We have the following result
\begin{lemma}
For any $\m\in\uuu^{(0)}\setminus\{\um\}$ and $\s\in\j(\m)$, there exists a neighborhood $\www$ of $\s$ such that
\bes
\forall x\in\www\setminus\{\s\},\ \ \ (U^\s)^{-1}(x)_1 - \s_1 = 0 \Rightarrow V(x) > V(\s).
\ees
\end{lemma}

\bp This is obvious seeing Assumption \ref{ass.morse} and recalling $\s\in\uuu^{(1)}$.

\ep

For $\tau,\delta>0$ small enough, this leads to a partition $E_{\m,\tau,\delta} = E^+_{\m,\tau,\delta}\sqcup E^-_{\m,\tau,\delta}$, where $E^+_{\m,\tau,\delta}$ is defined as the connected component of $E_{\m,\tau,\delta}$ containing $\m$. Observe that any path connecting $E^+_{\m,\tau,\delta}$ to $E^-_{\m,\tau,\delta}$ within the set $\{V<V(\j(\m))+\delta\}$ must intersect at least one $\ccc_{\s,\tau,\delta}$ due to Definition \ref{def:labeling} of the separating saddle points.

For $h>0$ and $\tau,\delta$ small enough, we define the function $\theta_\m$ on the sublevel set $E_-(\m)\cap\{V<V(\j(\m)) + 3\delta\}$ as follows. On the disjoint open sets $E^+_{\m,3\tau,3\delta}$ and $E^-_{\m,3\tau,3\delta}$, we define
\bes
\theta_\m(x)=\left\{\begin{aligned}
    1\mbox{ for } x\in E^+_{\m,3\tau,3\delta},\\
    -1\mbox{ for } x\in E^-_{\m,3\tau,3\delta}.
\end{aligned}\right.
\ees
Furthermore, for every $\s\in\j(\m)$, the change of variable $U^\s$ provided by Assumption \ref{ass.morse} is defined on $\ccc_{\s,3\tau,3\delta}$ for sufficiently small $\tau,\delta>0$. Thus, for $x\in\ccc_{\s,3\tau,3\delta}$ we can set
\bes
\theta_\m(x) = C_{\s,h}^{-1}\int_0^{\ell_\s(x)}\zeta(r/\tau)e^{-\frac{r^{\nu_1^\s}}{\nu_1^\s h}}dr,
\ees
where the function $\ell_\s$ is the one constructed in the previous section and its sign (see Remark \ref{rem:signell}) is chosen so that there exists a neighborhood $\www$ of $\s$ with
\bes
E(\m)\cap\www\subset\{(U^\s)^{-1}(x)_1 - \s_1>0\}.
\ees
Recall that $\zeta\in\ccc_c^\infty(\R,[0,1])$ is even, satisfies $\zeta = 1$ on $[-1,1]$, and $\zeta=0$ outside $[-2,2]$. We define the normalizing constant
\bes
C_{\s,h}=\frac 12\int_{-\infty}^{+\infty}\zeta(r/\tau)e^{-\frac{r^{\nu_1^\s }}{\nu_1^\s h}}dr.
\ees
Therefore, since, for every $\tau>0$ and then $\delta>0$ small enough, the sets $E^+_{\m,3\tau,3\delta}$, $E^-_{\m,3\tau,3\delta}$ and $\ccc_{\s,3\tau,3\delta}$ for $\s\in\j(\m)$ are mutually disjoint, $\theta_\m$ is well defined on their union, which forms $E_-(\m)\cap\{V<V(\j(\m))+3\delta\}$. Moreover, in a small neighborhood of the boundary shared by $\ccc_{\s,3\tau,3\delta}$ and $E_{\m,3\tau,3\delta}$, $(2\nu_1^\s|t_1^\s|)^{\frac1{\nu_1^\s}}|(U^\s)^{-1}(x)_1 - \s_1|\geq\frac52\tau$ or in other words, $|\ell_{\s,0,1}|\geq\frac52\tau$. Using that $\ell_\s = \ell_{\s,0,1} + O(|x-\s|^2 + h)$ and having that $|x-\s| = O(\delta)$ in this neighborhood, we thus obtain that, for every sufficiently small $\tau > 0$, and subsequently for $\delta,h > 0$ small enough, it holds that $|\ell_\s| \geq 2\tau$ in a neighborhood of the boundary shared by $\ccc_{\s,3\tau,3\delta}$ and $E_{\m,3\tau,3\delta}$. This shows that $\theta_\m$ is $\ccc^\infty$ on $E_-(\m)\cap\{V<V(\j(\m))+3\delta\}$.

Note also that there exist $\gamma,\beta>0$ such that
\bes
\begin{split}
    C_{\s,h} &= \frac 12\int_{-\infty}^{+\infty}\zeta(r/\tau)e^{-\frac{r^{\nu_1^\s }}{\nu_1^\s h}}dr = \int_0^{+\infty}\zeta(r/\tau)e^{-\frac{r^{\nu_1^\s }}{\nu_1^\s h}}dr\\
    &= \int_0^{+\infty}e^{-\frac{r^{\nu_1^\s }}{\nu_1^\s h}}dr + \int_0^{+\infty}(\zeta(r/\tau)-1)e^{-\frac{r^{\nu_1^\s }}{\nu_1^\s h}}dr\\
    &= \frac1{\nu_1^\s }\int_0^{+\infty}e^{-\frac{r}{\nu_1^\s h}}r^{\frac1{\nu_1^\s}-1}dr + O\big(\int_\gamma^{+\infty}e^{-\frac{r^{\nu_1^\s }}{\nu_1^\s h}}dr\big)\\
    &= \frac{(\nu_1^\s h)^{\frac1{\nu_1^\s}}}{\nu_1^\s}\int_0^{+\infty}e^{-r}r^{\frac1{\nu_1^\s}-1}dr + O(e^{-\beta/h})\\
    &= \frac{(\nu_1^\s h)^{\frac1{\nu_1^\s}}}{\nu_1^\s}\Gamma\big(\frac1{\nu_1^\s}\big) +  O(e^{-\beta/h}).
\end{split}
\ees
Hence
\be\label{eq:ch}
\exists\beta>0,\ C_{\s,h}^{-1}=\frac{\nu_1^\s}{(\nu_1^\s h)^{\frac1{\nu_1^\s}}}\Gamma\big(\frac1{\nu_1^\s}\big)^{-1}(1+O(e^{-\beta/h})),
\ee
when $\nu_1^\s = 2$, this simplifies to $C_{\s,h}^{-1} = \sqrt{\frac{2}{\pi h}}(1+O(e^{-\beta/h}))$.

We now extend $\theta_\m$ to a cutoff function defined on $\R^d$. Considering a smooth function $\chi_\m$ such that
\bes
\chi_\m(x) = \left\{\begin{aligned}
    &1\mbox{ for } x\in E_-(\m)\cap\{V \leq V(\j(\m)) + 2\delta\},\\
    &0\mbox{ for } x\in\R^d\setminus\big(E_-(\m)\cap\{V < V(\j(\m)) + 3\delta\}\big),
\end{aligned}\right.
\ees
then, $\chi_\m\theta_\m$ belongs to $\ccc_c^\infty(\R^d,[-1,1])$, and
\bes
\supp(\chi_\m\theta_\m)\subset E_-(\m)\cap\{V < V(\j(\m)) + 3\delta\}.
\ees

\begin{defin}\label{def:quasimodes}
For $\tau>0$ and then $\delta,h>0$ small enough, we define the quasimodes 
\bes
\left\{\begin{array}{l}
    \psi_{\um}(x)=2e^{-\frac{V(x)-V(\um)}{h}}\\
    \psi_{\m}(x)=\chi_\m(x)(\theta_\m(x)+1)e^{-\frac{V(x)-V(\m)}{h}}\quad\mbox{ for }\m\in\uuu^{(0)}\setminus\{\um\}.
\end{array}\right.
\ees
we also define the normalized quasimodes for $\m\in\uuu^{(0)}$ by
\bes
\phii_{\m}=\frac{\psi_{\m}}{\vvert{\psi_{\m}}}.
\ees
\end{defin}

Based on this definition, we state the following lemma, which provides an initial relationship between the quasimodes:
\begin{lemma}\label{lem:supppsi}Let $\m,\m'\in\uuu^{(0)}$ with $\m\neq\m'$. Then

    If $\bsigma(\m) = \bsigma(\m')$ and $\j(\m)\cap\j(\m')=\emptyset$, then $\supp(\psi_\m)\cap\supp(\psi_{\m'}) = \emptyset$.

    If $\bsigma(\m) > \bsigma(\m')$, then one of the following holds
    \begin{itemize}
        \item[$\star$] $\supp(\psi_\m)\cap\supp(\psi_{\m'}) = \emptyset$,
        \item[$\star$] $\psi_\m = 2e^{-(V-V(\m))/h}$ on $\supp(\psi_{\m'})$.
    \end{itemize}
\end{lemma}

\bp Since $E(\m)$ is a connected component of $\{V<\bsigma(\m)\}$, the boundary of $\overline{E(\m)}$ consists of non-critical points of $V$ (points where $\nabla V\neq0$) as well as separating saddle points $\s\in\j(\m)$, by the definition of separating saddle points. Therefore, from the definition of $\psi_\m$, we observe that for all $\eps>0$,
\bes
\supp\psi_\m\subset \overline{E^+_{\m,3\tau,3\delta}\cup\bigcup_{\s\in\j(\m)}\ccc_{\s,3\tau,3\delta}}\subset\overline{E(\m)}+B(0,\eps)
\ees
for $\tau,\delta$ small enough.

Now if $\bsigma(\m) = \bsigma(\m')$ then necessarily $E(\m)\cap E(\m') = \emptyset$. If this were not the case, then since they are critical components of $\{V\leq\bsigma(\m)\}$, they would be the same which is in a contradiction with the construction of $E$. When in addition $\j(\m)\cap\j(\m')=\emptyset$, then
\bes
\overline{E(\m)}\cap\overline{E(\m')} = \partial E(\m)\cap\partial E(\m') = \j(\m)\cap\j(\m') = \emptyset,
\ees
hence with $\eps$ sufficiently small we have $\supp(\psi_\m)\cap\supp(\psi_{\m'}) = \emptyset$.

If $\bsigma(\m) > \bsigma(\m')$, then either $\m'\notin E(\m)$, in which case we conclude with the same argument that $\supp(\psi_\m)\cap\supp(\psi_{\m'}) = \emptyset$. If instead $\m'\in E(\m)$, then $\overline{E(\m')}\subset E_-(\m')\subset E(\m)$. Moreover, $\chi_\m\theta_\m\equiv1$ in a neighborhood of $\overline{E(\m)}\setminus\bigcup_{\s\in\j(\m)}\ccc_{\s,3\tau,3\delta}$, and since $\j(\m)\cap\supp(\psi_{\m'})=\emptyset$, the result follows.

\ep

Recall $\k$ is introduced in \eqref{eq:defik}, we define
\be\label{eq:defmu}
\mu(\m) = \inf_{\s\in\j(\m)}(2+2\k(\s)-2\k(\m)-\frac2{\nu_1^\s}),
\ee
which is an exponent that will appear in subsequent estimates. We also define $\tilde\j(\m)\subset\j(\m)$ as the set of saddle points at which the infimum is attained, that is
\be\label{eq:deftildej}
\s_0\in\tilde\j(\m) \iff 2\k(\s_0)-\frac2{\nu_1^{\s_0}} = \min_{\s\in\j(\m)}\big(2\k(\s)-\frac2{\nu_1^\s}\big).
\ee
Even though $\tilde\j(\m)\subset\j(\m)$ it has no reason to be comparable to any set of the partition \eqref{eq:defjautre}. For all $\m\in\uuu^{(0)}$ and all $\s\in\j(\m)$, we define the following quantity
\be\label{eq:defgamma}
\gamma(\s) = \left\{\begin{array}{ll}
    +\infty &\mbox{ if }\s\in\j_\infty(\m),\\
    \max\Big(-\frac2{\underline{\nu}_*^\s},\min\big(-\frac{2}{\hat\nu^\s},2\frac{\hat\nu^\s-\underline\nu^\s}{\overline\nu^\s}-\frac{2}{\underline\nu^\s}\big)\Big) + \frac2{\overline{\nu}^\s} + \frac2{\overline{\nu}}&\mbox{ if }\s\in\j_1(\m),\\
    \frac2{\overline{\nu}^\s} - \frac2{\underline{\nu}_*^\s} + \frac2{\overline{\nu}} &\mbox{ if }\s\in\j_2(\m),
\end{array}\right.
\ee
where the sets $\j_\infty,\j_1$, and $\j_2$ are defined after \eqref{eq:defjautre}, and the notation for the $\nu$ parameters follows \eqref{eq:notationnu}. Finally, we define
\be\label{eq:defalpha}
\alpha(\m) = \min_{\s\in\j(\m)}\gamma(\s).
\ee

\begin{proposition}\label{prop:interactionmatrix}
Under Assumption \eqref{ass.gener}, there exist constants $C,\beta>0$, such that for all $\tau>0$ and then $\delta,h>0$ small enough, the following holds for all $\m,\m'\in\uuu^{(0)}$,
\begin{itemize}
    \item[$i)$] $\<\phii_\m,\phii_{\m'}\> = \delta_{\m,\m'}+O(e^{-C/h}),$
    \item[$ii)$] $\<\Delta_V\phii_\m,\phii_\m\> = \sum_{\s\in\tilde\j(\m)}z(\m,\s) h^{\mu(\m)}e^{-2S(\m)/h}(1+O(h^\beta))$, where
    \be\label{eq:zms}
    z(\m,\s) = \bigg(\frac{\nu_1^\s(2|t_1^\s|)^{\frac{1}{\nu_1^\s}}}{2\Gamma(\frac{1}{\nu_1^\s})}\bigg)^2\prod_{i=1}^d\frac{\Gamma(\frac{1}{\nu_i^\s})}{\Gamma(\frac{1}{\nu_i^\m})}\frac{\nu_i^\m(2t_i^\m)^{\frac{1}{\nu_i^\m}}}{\nu_i^\s(2|t_i^\s|)^{\frac{1}{\nu_i^\s}}} > 0
    \ee
    and $\mu$ defined in \eqref{eq:defmu}.
    \item[$iii)$] $\vvert{\Delta_V\phii_\m}^2 = O(h^{\alpha(\m)+2-\frac{2}{\overline\nu}})\<\Delta_V\phii_\m,\phii_\m\>$, where $\alpha$ is defined in \eqref{eq:defalpha}.
\end{itemize}
Here, $S(\m)=V(\j(\m))-V(\m)$ for $\m\neq\um$, and $S(\um)=+\infty$, as defined in \eqref{eq:defS}.
\end{proposition}

\begin{remark}
    Observe that if $V$ has the same degeneracy orders at $\m$ and at all saddle poitns $\s\in\j(\m)$ i.e. $\forall \s\in\j(\m),\ \forall 1\leq i\leq d$, $\nu_i^\s = \nu_i^\m$, then point $ii)$ simplifies significantly:
    \bes
    \<\Delta_V\phii_\m,\phii_\m\> = \sum_{\s\in\j(\m)}z(\m,\s) h^{2-\frac2{\nu_1^\s}}e^{-2S(\m)/h}(1+O(h^\beta))
    \ees
    with
    \bes
    z(\m,\s) = \bigg(\frac{\nu_1^\s(2|t_1^\s|)^{\frac{1}{\nu_1^\s}}}{2\Gamma(\frac{1}{\nu_1^\s})}\bigg)^2\bigg(\prod_{i=1}^d\frac{t_i^\m}{|t_i^\s|}\bigg)^{\frac{1}{\nu_i^\m}}.
    \ees
    And thus if $V$ is Morse we recover \cite[Proposition 5.1]{BoLePMi22}.
\end{remark}

\bp We follow the argument of \cite[Proposition 5.1]{BoLePMi22} and omit details of the steps that remain unchanged.

In this proof, we apply several versions of the degenerate Laplace method (DLM), as presented in Proposition \ref{LaplaceMethod}, justified by Assumptions \ref{ass.morse}, \ref{ass.morsebis}, and Remark \ref{rem:relax}.

Since $V$ attains its unique global minimum at $\m$ on $\supp\psi_\m$, applying a DLM to $2V$ yields
\bes
\begin{split}
    \vvert{\psi_\m}^2 &= \int_{\R^d}\chi_\m^2(\theta_\m+1)^2e^{-2\frac{V-V(\m)}{h}}\\
    &= \chi_\m^2(\m)(\theta_\m(\m)+1)^2h^{\sum_{i=1}^d\frac1{\nu_i^\m}}\prod_{i=1}^d\frac{2\Gamma(\frac{1}{\nu_i^\m})}{\nu_i^\m(2t_i^\m)^{\frac{1}{\nu_i^\m}}}(1+O(h^{\frac{2}{\overline{\nu}^\m}})).
\end{split}
\ees
This leads to
\be\label{eq:eclpsi}
\vvert{\psi_\m}=2h^{\k(\m)}\sqrt{\prod_{i=1}^d\frac{2\Gamma(\frac{1}{\nu_i^\m})}{\nu_i^\m(2t_i^\m)^{\frac{1}{\nu_i^\m}}}}(1+O(h^{\frac{2}{\overline{\nu}^\m}})).
\ee
In the Morse case, $\forall i\in\lb1,d\rb,\ \nu_i^\m=2$ and thus we retrieve
\bes
\prod_{i=1}^d\frac{2\Gamma(\frac{1}{\nu_i^\m})}{\nu_i^\m(2t_i^\m)^{\frac{1}{\nu_i^\m}}} = \prod_{i=1}^d\frac{\sqrt{\pi}}{\sqrt{2t_i^\m}} = \pi^{\frac d2}|\det\Hess_\m V|^{-\frac12},
\ees
then, the proof of $i)$ is exactly the same as in \cite{BoLePMi22} using Lemma \ref{lem:supppsi}.

For $ii)$, using that $\d_V = e^{-V/h}\circ h\nabla\circ e^{V/h}$,
$\d_V\psi_\m = h\nabla(\chi_\m(\theta_\m+1))e^{-\frac{V-V(\m)}{h}}$, recalling the definitions of $\chi_\m$ and $\theta_\m$,
\bes
\supp\nabla\chi_\m\subset \{V\geq V(\j(\m)) + 2\delta\}\ \ \mbox{ and }\ \ \supp\nabla\theta_\m\subset\bigcup_{\s\in\j(\m)}\ccc_{\s,3\tau,3\delta},
\ees
with $\nabla\theta_\m = C_{\s,h}^{-1}\zeta(\ell_\s/\tau)e^{-\frac{\ell_\s^{\nu_1^\s}}{\nu_1^\s h}}\nabla\ell_\s$ on each $\ccc_{\s,3\tau,3\delta}$. Therefore
\bes
\begin{split}
    \<\Delta_V\psi_\m,\psi_\m\> &= \vvert{\d_V\psi_\m}^2\\
    &= h^2\int_{\bigcup_{\s\in\j(\m)}\ccc_{\s,3\tau,3\delta}}\chi_\m^2|\nabla\theta_\m|^2e^{-2\frac{V-V(\m)}{h}}\\&\phantom{*********} + h^2\int_{\{V\geq V(\j(\m)) + 2\delta\}}\nabla\chi_\m\cdot\nabla(\chi_\m(\theta_\m+1)^2)e^{-2\frac{V-V(\m)}{h}}\\
    &= h^2\sum_{\s\in\j(\m)}C_{\s,h}^{-2}\int_{\ccc_{\s,3\tau,3\delta}}\chi_{\m}^2\zeta(\ell_\s/\tau)^2|\nabla\ell_\s|^2e^{-2\big(V+\frac{\ell^{\nu_1^\s}}{\nu_1^\s}-V(\m)\big)/h}\\ &\phantom{*********} + O(e^{-2(S(\m)+2\delta)/h}).
\end{split}
\ees
By the construction of $\ell$ and Lemma \ref{lem:ell} $vii)$, $\s$ is the unique minima of $V+\frac{\ell^{\nu_1^\s}_{\s,0}}{\nu_1^\s}-V(\m)$ on $\ccc_{\s,3\tau,3\delta}$ and
\bes
\big(V+\frac{\ell^{\nu_1^\s}_{\s,0}}{\nu_1^\s}-V(\m)\big)(\s) = V(\s) - V(\m) = V(\j(\m)) - V(\m) = S(\m).
\ees
Moreover, thanks to Remark \ref{rem:relax} and Lemma \ref{lem:ell} $vii)$, $V+\frac{\ell^{\nu_1^\s}_{\s,0}}{\nu_1^\s}$ satisfies Assumption \ref{ass.morse}, thus we can apply a DLM and using \eqref{eq:ch}, we deduce
\be\label{eq:Ppsipsi}
\begin{split}
    \<\Delta_V\psi_\m,\psi_\m\> &= h^2\sum_{\s\in\j(\m)}C_{\s,h}^{-2}\chi_\m^2(\s)\zeta(\ell_{\s,0}(\s)/\tau)^2|\nabla\ell_{\s,0}(\s)|^2h^{2\k(\s)}\prod_{i=1}^d\frac{2\Gamma(\frac{1}{\nu_i^\s})}{\nu_i^\s(2|t_i^\s|)^{\frac{1}{\nu_i^\s}}}\\
    &\phantom{******}\times e^{-2S(\m)/h}(1+O(h^{\frac{2}{\overline{\nu}^\s}}))\\
    &= \sum_{\s\in\j(\m)}(z_1(\s)h^{2\k(\s)-\frac2{\nu_1^\s}})h^2e^{-2S(\m)/h}(1+O(h^{\frac{2}{\overline{\nu}^{\j(\m)}}}))\\
    &= \sum_{\s\in\tilde\j(\m)}(z_1(\s)h^{2\k(\s)-\frac2{\nu_1^\s}})h^2e^{-2S(\m)/h}(1+O(h^{\beta_1})),
\end{split}
\ee
where we recall $\overline \nu^{\s} = \sup_{1\leq i\leq d}\nu_i^\s$, $\overline{\nu}^{\j(\m)} = \sup_{\s\in\j(\m)}\overline\nu^{\s}$ and $\overline{\nu} = \sup_{x^*,i}\nu_i^{x^*}$ from \eqref{eq:notationnu}, with
\bes
z_1(\s) = \frac{(\nu_1^\s)^2}{(\nu_1^\s)^{\frac2{\nu_1^\s}}}\Gamma\big(\frac1{\nu_1^\s}\big)^{-2}(2\nu_1^\s |t_1|^\s)^\frac2{\nu_1^\s}\prod_{i=1}^d\frac{2\Gamma(\frac{1}{\nu_i^\s})}{\nu_i^\s(2|t_i^\s|)^{\frac{1}{\nu_i^\s}}} = \bigg(\frac{\nu_1^\s(2|t_1^\s|)^{\frac{1}{\nu_1^\s}}}{\Gamma(\frac{1}{\nu_1^\s})}\bigg)^2\prod_{i=1}^d\frac{2\Gamma(\frac{1}{\nu_i^\s})}{\nu_i^\s(2|t_i^\s|)^{\frac{1}{\nu_i^\s}}}
\ees
and denoting $\beta_1 = \min\big(\frac{2}{\overline{\nu}^{\j(\m)}},\beta'\big)$, with
\bes
\beta' = \min_{\s\in\j(\m)\setminus\tilde\j(\m)}\big(2\k(\s)-\frac2{\nu_1^\s}\big) - \min_{\s\in\j(\m)}\big(2\k(\s)-\frac2{\nu_1^\s}\big).
\ees
By the definition of $\tilde\j(\m)$, $\beta'>0$. We use the convention that $\min\emptyset = +\infty$ so that $\beta' = +\infty$ when $\tilde\j(\m) = \j(\m)$, in which case, we just have $\beta_1 = \frac{2}{\overline{\nu}^{\j(\m)}}$.

When $V$ is Morse, we recover
\bes
z_1(\s) = 4\bigg(\frac{(2|t_1^\s|)^{\frac{1}{2}}}{\Gamma(\frac{1}{2})}\bigg)^2\prod_{i=1}^d\frac{\Gamma(\frac{1}{2})}{(2|t_i^\s|)^{\frac{1}{2}}} = \frac{2|4t_1^\s|}{\pi}\pi^\frac d2 |\det\Hess_\s V|^{-\frac12}
\ees
and we recognize exactly the prefactor in \cite[(5.4)]{BoLePMi22}.

Combining \eqref{eq:Ppsipsi} and \eqref{eq:eclpsi}, we obtain
\bes
\begin{split}
    \<\Delta_V\phii_\m,\phii_\m\> &= \frac14h^{-2\k(\m)}\prod_{i=1}^d\frac{\nu_i^\m(2t_i^\m)^{\frac{1}{\nu_i^\m}}}{2\Gamma(\frac{1}{\nu_i^\m})}\sum_{\s\in\tilde\j(\m)}z_1(\s)h^{2\k(\s)-\frac2{\nu_1^\s}}\\&\phantom{********}\times h^2e^{-2S(\m)/h}(1+O(h^\beta))(1+O(h^{\frac{2}{\overline{\nu}^\m}}))\\
    &= \sum_{\s\in\tilde\j(\m)}z(\m,\s) h^{2+2\k(\s)-2\k(\m)-\frac2{\nu_1^\s}}e^{-2S(\m)/h}(1+O(h^{\min(\beta_1,\frac{2}{\overline{\nu}^\m})}))
\end{split}
\ees
with 
\bes
z(\m,\s) = \bigg(\frac{\nu_1^\s(2|t_1^\s|)^{\frac{1}{\nu_1^\s}}}{2\Gamma(\frac{1}{\nu_1^\s})}\bigg)^2\prod_{i=1}^d\frac{\Gamma(\frac{1}{\nu_i^\s})}{\Gamma(\frac{1}{\nu_i^\m})}\frac{\nu_i^\m(2t_i^\m)^{\frac{1}{\nu_i^\m}}}{\nu_i^\s(2|t_i^\s|)^{\frac{1}{\nu_i^\s}}},
\ees
this proves $ii)$ with
\be\label{eq:beta}
\beta = \beta(\m) = \min(\beta_1,\frac{2}{\overline{\nu}^\m}).
\ee

We now prove point $iii)$. We begin as in the proof of \cite[Proposition 5.1, $iii)$]{BoLePMi22},
\be\label{eq:Dvpsim}
\vvert{\Delta_V\psi_\m}^2 = \vvert{\Delta_V(\theta_\m e^{-(V-V(\m))/h})}^2_{L^2(\supp\chi_\m)} + O(e^{-2(S(\m)+2\delta)/h}).
\ee
On $\supp\chi_\m$, the function $\Delta_V(\theta_\m e^{-(V-V(\m))/h})$ is supported in $\bigcup_{\s\in\j(\m)}\ccc_{\s,3\tau,3\delta}$. Thus, using \eqref{eq:wplusr}, we need to consider three cases:

\textbullet\ Let $\s\in\j_\infty(\m)$. If it satisfies Assumption \ref{ass.case1}, then by \eqref{eq:PchiCase1} and applying a DLM, this yields
\bes
\vvert{\Delta_V(\theta_\m e^{-(V-V(\m))/h})}^2_{L^2(\ccc_{\s,3\tau,3\delta})} = O(h^\infty)e^{-2S(\m)/h} = O(h^\infty)\<\Delta_V\psi_\m,\psi_\m\>_{L^2(\ccc_{\s,3\tau,3\delta})}.
\ees
If $\s$ satisfies Assumption \ref{ass.case2}, then by \eqref{eq:PchiCase2} and a DLM we also have
\bes
\begin{split}
    \vvert{\Delta_V(\theta_\m e^{-(V-V(\m))/h})}^2_{L^2(\ccc_{\s,3\tau,3\delta})} &= \int_{\ccc_{\s,3\tau,3\delta}} O(x^\infty)e^{-2\big(V-V(\m)+\frac{\ell_\s^{\nu_1^\s}}{\nu_1^\s}\big)/h}\\
    &= O(h^\infty)e^{-2S(\m)/h}\\
    &= O(h^\infty)\<\Delta_V\psi_\m,\psi_\m\>_{L^2(\ccc_{\s,3\tau,3\delta})},
\end{split}
\ees
The only difference from the previous case lies in the DLM: here, we must account for the $O(x^\infty)$ term, which requires using \eqref{eq:Laplace3} instead of \eqref{eq:Laplace2}.

\textbullet\ Let $\s\in\j_1(\m)$ which satisfies Assumption \ref{ass.case3}. Recalling \eqref{eq:PchiCase3} and \eqref{eq:ch}, we have on $\ccc_{\s,3\tau,3\delta}$
\bes
\begin{split}
    \Delta_V(\theta_\m e^{-(V-V(\m))/h}) &= h^{1-\frac{1}{\nu_1^\s}}\bigg(|U^{-1}(x)|\sum_{j\notin\lb2,J\rb}|U^{-1}(x)_j|^{\nu_j^\s-1}\\&\phantom{*******} + |U^{-1}(x)|^{\hat\nu^\s-1}\sum_{j\in\lb2,J\rb}|U^{-1}(x)_j| + h\bigg)e^{-\big(V-V(\m)+\frac{\ell_\s^{\nu_1^\s}}{\nu_1^\s}\big)/h}
\end{split}
\ees
hence we get
\bes
\begin{split}
    \vvert{\Delta_V(\theta_\m e^{-(V-V(\m))/h})}^2_{L^2(\ccc_{\s,3\tau,3\delta})} &= O\Big(h^{\frac{2}{\overline\nu^\s}}h^{2-\frac{2}{\hat\nu^\s}} + h^{2\frac{\hat\nu^\s-1}{\overline\nu^\s}}h + h^2\Big)h^{2 + 2\k(\s) - \frac{2}{\nu_1^\s}}e^{-2S(\m)/h}\\
    &= O(h^{\min(2+\frac{2}{\overline\nu^\s}-\frac{2}{\hat\nu^\s},1+2\frac{\hat\nu^\s-1}{\overline\nu^\s})})\<\Delta_V\psi_\m,\psi_\m\>_{L^2(\ccc_{\s,3\tau,3\delta})}.
\end{split}
\ees
Here, each term $|U^{-1}(x)_j|$ contributes a factor of $h^{\frac{2}{\sup_{\lb2,J\rb}\nu_j^\s}} = h$, since $\nu_j^\s = 2$ for all $j\in\lb2,J\rb$. So we can write it the following way
\bes
\vvert{\Delta_V(\theta_\m e^{-(V-V(\m))/h})}^2_{L^2(\ccc_{\s,3\tau,3\delta})} = h^{2+\frac{2}{\overline\nu^\s}}O(h^{\min(-\frac{2}{\hat\nu^\s},2\frac{\hat\nu^\s-\underline\nu^\s}{\overline\nu^\s}-\frac{2}{\underline\nu^\s})})\<\Delta_V\psi_\m,\psi_\m\>_{L^2(\ccc_{\s,3\tau,3\delta})}.
\ees

Similarly, if $\s$ satisfies Assumption \ref{ass.case5}, recalling \eqref{eq:PchiCase5} and \eqref{eq:ch}, we have on $\ccc_{\s,3\tau,3\delta}$
\bes
\begin{split}
    \Delta_V(\theta_\m e^{-(V-V(\m))/h}) &= h^{1-\frac{1}{\nu_1^\s}}\bigg(|U^{-1}(x)|\sum_{j<d}|U^{-1}(x)_j|^{\nu_j^\s-1}\\&\phantom{*******} + |U^{-1}(x)|^{\hat\nu^\s-\underline\nu^\s+1}|U^{-1}(x)_d|^{\nu_d^\s-1} + h\bigg)e^{-\big(V+\frac{\ell^{\nu_1}}{\nu_1}\big)/h}
\end{split}
\ees
which leads to
\bes
\begin{split}
    \vvert{\Delta_V(\theta_\m e^{-(V-V(\m))/h})}^2_{L^2(\ccc_{\s,3\tau,3\delta})} &= O\Big(h^{\frac{2}{\overline\nu^\s}}h^{2-\frac{2}{\hat\nu^\s}} + h^{2\frac{\hat\nu^\s-\underline\nu^\s+1}{\overline\nu^\s}}h^{2-\frac{2}{\underline\nu^\s}}\Big)h^{2 + 2\k(\s) - \frac{2}{\nu_1^\s}}e^{-2S(\m)/h}\\
    &= h^{2+\frac{2}{\overline\nu^\s}}O(h^{\min(-\frac{2}{\hat\nu^\s},2\frac{\hat\nu^\s-\underline\nu^\s}{\overline\nu^\s}-\frac{2}{\underline\nu^\s})})\<\Delta_V\psi_\m,\psi_\m\>_{L^2(\ccc_{\s,3\tau,3\delta})}.
\end{split}
\ees

\textbullet\ Finally, for $\s\in\j_2(\m)$, recalling \eqref{eq:PchiCase0} along with \eqref{eq:ch}, we have on $\ccc_{\s,3\tau,3\delta}$
\bes
\begin{split}
    \Delta_V(\theta_\m e^{-(V-V(\m))/h}) &= h^{1-\frac{1}{\nu_1^\s}}O\Big(\sum_{j\geq2}|U^{-1}(x)||U^{-1}(x)_j|^{\nu_j-1} + h\Big)e^{-\big(V-V(\m)+\frac{\ell_\s^{\nu_1^\s}}{\nu_1^\s}\big)/h}
\end{split}
\ees
hence we get
\bes
\begin{split}
    \vvert{\Delta_V(\theta_\m e^{-(V-V(\m))/h})}^2_{L^2(\ccc_{\s,3\tau,3\delta})} &= h^{2-\frac{2}{\nu_1^\s}}O\Big(\sum_{j\geq2}h^{\frac2{\overline{\nu}^\s}}h^{2-\frac{2}{\nu_{j}^\s}} + h^2\Big)h^{2\k(\s)}e^{-2S(\m)/h}\\
    &= O\Big(h^{2 + \frac2{\overline{\nu}^\s}}\sum_{j\geq2}h^{-\frac{2}{\nu_{j}^\s}}\Big)h^{2 + 2\k(\s) - \frac{2}{\nu_1^\s}}e^{-2S(\m)/h}\\
    &= O(h^{2+\frac2{\overline{\nu}^\s}-\frac{2}{\underline{\nu}_*^\s}})\<\Delta_V\psi_\m,\psi_\m\>_{L^2(\ccc_{\s,3\tau,3\delta})}.
\end{split}
\ees
Note that for saddle points in $\j_1$, we could alternatively apply the method used for $\j_2$. However, we select the approach that yields the highest power of $h$, which justifies the $\max$ in \eqref{eq:defgamma}. Combining these results with \eqref{eq:Dvpsim} and \eqref{eq:defalpha}, we have thus proved point $iii)$.

Furthermore, if Assumption \ref{ass.case4} holds, we have using \eqref{eq:PchiCase4},
\bes
\begin{split}
    \vvert{\Delta_V(\theta_\m e^{-(V-V(\m))/h})}^2_{L^2(\ccc_{\s,3\tau,3\delta})} &= h^{2 + 2\k(\s) - \frac{2}{\nu_1^\s}}O(h^4)e^{-2S(\m)/h}\\
    &= O(h^4)\<\Delta_V\psi_\m,\psi_\m\>_{L^2(\ccc_{\s,3\tau,3\delta})},
\end{split}
\ees
resulting in
\be\label{eq:estimecase4}
\vvert{\Delta_V\phii_\m}^2 = O(h^4)\<\Delta_V\phii_\m,\phii_\m\>.
\ee

\ep

\subsection{Proof of Theorem \ref{thm:2} and graded matrices}\label{ssec:proofthm2}
The proof follows the same strategy as in \cite{BoLePMi22}. We refer to that work and to \cite{LePMi20} for the details, and here we briefly outline the main arguments. From now on, we relabel the minima by increasing saddle height:
\bes
\forall j\in\lb1,n_0-1\rb,\ \ S(\m_j)\leq S(\m_{j+1}),
\ees
recalling $S$ is defined in \eqref{eq:defS}. Thus, $\m_{n_0} = \um$ because $S(\um) = +\infty$ and $\um$ is the only minimum with this property by construction. For brevity, we introduce the following notations
\bes
\forall j\in\lb1,n_0\rb,\ \ S_j = S(\m_j),\ \ \phii_j = \phii_{\m_j},\ \mbox{ and } \ \tilde{\lambda}_j = \<\Delta_V\phii_j,\phii_j\>.
\ees
It follows that
\be\label{eq:Pmm'}
\forall j,k\in\lb1,n_0\rb,\ \ \<\Delta_V\phii_j,\phii_k\> = \delta_{j,k}\tilde{\lambda}_j.
\ee
The identity is immediate when $j = k$. When $j \neq k$, Lemma \ref{lem:supppsi} implies either $\supp \phii_j \cap \supp \phii_k = \emptyset$ or $\phii_j = c_h e^{-(V - V(\m_j))/h}$ on $\supp \phii_k$, where $c_h$ is a normalizing constant. The same applies with $j$ and $k$ swapped, due to the self-adjointness of $\Delta_V$. Since $\Delta_V e^{-V/h} = 0$, this confirms the validity of \eqref{eq:Pmm'}. Hence, denoting
\bes
\Pi = \frac{1}{2i\pi}\int_{\partial D(0,\frac\eps2h^{2-\frac2{\overline{\nu}}})}(z-\Delta_V)^{-1}dz
\ees
the spectral projector on the small eigenvalues of $\Delta_V$, where $\eps$ is given by Theorem \ref{thm:1}, and $u_j = \Pi\phii_j$ for $j\in\lb1,n_0\rb$ (note that $u_{n_0} = \phii_{\um}$), we obtain the following proposition.

\begin{proposition}
    There exists $c>0$ such that for all $j,k\in\lb1,n_0\rb$ and sufficiently small $h>0$, one has
    \be\label{eq:quasiortho}
    \<u_j,u_k\> = \delta_{j,k} + O(e^{-c/h})
    \ee
    and
    \be\label{eq:interaction}
    \<\Delta_Vu_j,u_k\> = \delta_{j,k}\tilde{\lambda}_j + O\Big(h^{\frac12(\alpha(\m_j)+\alpha(\m_k))}\sqrt{\tilde{\lambda}_j\tilde{\lambda}_k}\Big),
    \ee
    where $\alpha$ is defined in \eqref{eq:defalpha}.
\end{proposition}

\bp Using that
\bes
\Pi-1 = \frac{1}{2i\pi}\int_{\partial D(0,\frac\eps2h^{2-\frac2{\overline{\nu}}})}(z-\Delta_V)^{-1}\Delta_V\frac{dz}{z},
\ees
for $u\in D(\Delta_V)$, $\vvert{(\Pi-1)u}\leq\sup_{\partial D(0,\frac\eps2h^{2-\frac2{\overline{\nu}}})}\Vert(z-\Delta_V)^{-1}\Vert\vvert{\Delta_Vu}$, and thus
\bes
\<u_j,u_k\> = \<\phii_j,\phii_k\> + O\Big(\frac{\Vert\Delta_V\phii_j\Vert + \Vert\Delta_V\phii_k\Vert}{d(\frac\eps2h^{2-\frac2{\overline{\nu}}},\sigma(\Delta_V))}\Big) + O\Big(\frac{\Vert\Delta_V\phii_j\Vert\Vert\Delta_V\phii_k\Vert}{d(\frac\eps2h^{2-\frac2{\overline{\nu}}},\sigma(\Delta_V))^2}\Big)
\ees
using that $\Delta_V$ is self-adjoint, which yields the first result, by Proposition \ref{prop:interactionmatrix}. Thanks to Theorem \ref{thm:1} and \eqref{eq:Pmm'}, we then deduce
\bes
\begin{split}
    \<\Delta_Vu_j,u_k\> &= \<\Delta_V\phii_j,\phii_k\> + \<\Delta_V(\Pi-1)\phii_j,\phii_k\> + \<\Delta_V\Pi\phii_j,(\Pi-1)\phii_k\>\\
    &= \delta_{j,k}\tilde{\lambda}_j + O\Big(\frac{\Vert\Delta_V\phii_j\Vert\Vert\Delta_V\phii_k\Vert}{d(\frac\eps2h^{2-\frac2{\overline{\nu}}},\sigma(\Delta_V))}\Big)\\
    &= \delta_{j,k}\tilde{\lambda}_j + O\Big(h^{\frac12(\alpha(\m_j)+2-\frac2{\overline{\nu}})}\sqrt{\tilde{\lambda}_{j}}\ h^{-2+\frac2{\overline{\nu}}}\ h^{\frac12(\alpha(\m_k)+2-\frac2{\overline{\nu}})}\sqrt{\tilde{\lambda}_{k}}\Big).
\end{split}
\ees

\ep

For every $j,k\in\lb1,n_0\rb$, we define
\bes
\alpha_{j,k} = \frac{\alpha(\m_j) + \alpha(\m_k)}{2}
\ees
hence \eqref{eq:interaction} can be rewritten
\bes
\<\Delta_Vu_j,u_k\> = \delta_{j,k}\tilde{\lambda}_j + O\Big(h^{\alpha_{j,k}}\sqrt{\tilde{\lambda}_j\tilde{\lambda}_k}\Big).
\ees

Thus we need to ensure that $\forall j,k\in\lb1,n_0\rb,\ \alpha_{j,k} > 0$ in order to be able to correctly compute the eigenvalues, hence we define
\be\label{eq:alpha0}
\alpha_0 = \inf_{j\in\lb1,n_0\rb}\alpha_{j,j}
\ee
and we recall that Assumption \ref{ass.alpha} consists in requiring $\alpha_0 > 0$.
\begin{remark}
    In the Morse case, we have $\j(\m) = \j_\infty(\m)$, so formally $\alpha_0 = +\infty$. This indicates that there is no obstruction to extracting the low-lying eigenvalues of $\Delta_V$. On the other hand, Assumption \ref{ass.case2} is satisfied, and thus we can solve \eqref{eq:eikapprox} and \eqref{eq:transpapprox} leading to a very sharp description of the eigenvalues \cite{BoLePMi22}, thus the formal $\alpha_0 = +\infty$ makes sense.
\end{remark}

Therefore \eqref{eq:interaction} implies that $\<\Delta_Vu_j,u_k\> = \delta_{j,k}\tilde{\lambda}_j + O\Big(h^{\alpha_0}\sqrt{\tilde{\lambda}_j\tilde{\lambda}_k}\Big)$ and thanks to Assumption \ref{ass.alpha}, this endows the interaction matrix with a graded structure, which we develop below. Then we use the Gram-Schmidt process to transform the basis $(u_{n_0-j+1})_{1\leq j\leq n_0}$ into an orthonormal basis $(e_{n_0-j+1})_{1\leq j\leq n_0}$ of $\Ran\Pi$. Moreover, thanks to \eqref{eq:quasiortho}, we have, for all $j\in\lb1,n_0\rb$,
\bes
e_j = u_j + O(e^{-c/h}),
\ees
see \cite[Lemma 4.11]{LePMi20} for details. And thus, using \eqref{eq:interaction} along with \cite[Proposition 4.12]{LePMi20},
\be\label{eq:Pee}
\forall j\in\lb1,n_0\rb,\ \ \<\Delta_Ve_j,e_k\> = \delta_{j,k}\tilde{\lambda}_j + O\Big(h^{\alpha_0}\sqrt{\tilde{\lambda}_j\tilde{\lambda}_k}\Big).
\ee
Using its graded structure, we can now compute the eigenvalues of the interaction matrix
\be\label{eq:matrix}
M := (\<\Delta_Ve_j,e_k\>)_{j,k} = {\Delta_V}_{|\Ran\Pi}.
\ee
We recall the results stated in \cite{LePMi20}:

We denote by $\Dr_0(E)$ the set of invertible and diagonalizable complex matrices on a Euclidean space $E$.

\begin{defin}\cite[Definition A.1]{LePMi20}
    Let $\Er = (E_j)_{1\leq j \leq p}$ be a sequence of vector spaces $E_j$ of finite dimension $r_j>0$, let $E = \bigoplus_{j=1}^pE_j$ and let $\tau = (\tau_i)_{2\leq i\leq p}\in(\R_+^*)^{p-1}$. Suppose that $(h,\tau)\mapsto \mmm_h(\tau)$ is a map from $(0,1]\times(\R_+^*)^{p-1}$ to the set of complex matrices acting on $E$.

    We say that $\mmm_h(\tau)$ is an $(\Er,\tau,h)$-graded matrix if there exists $\mmm'\in\Dr_0(E)$ independent of $(h,\tau)$ such that $\mmm_h(\tau) = \Omega(\tau)(\mmm' + O(h))\Omega(\tau)$, where $\Omega(\tau)$ and $\mmm'$ satisfy
    \begin{itemize}
        \item $\mmm' = \diag(M_j,1\leq j\leq p)$ with $M_j\in\Dr_0(E_j)$,
        \item $\Omega(\tau) = \diag(\eps_j(\tau)I_{r_j},1\leq j\leq p)$ with $\eps_1(\tau) = 1$ and $\eps_j(\tau) = \prod_{k=2}^j\tau_k$ for  $j\geq2$.
    \end{itemize}
\end{defin}

\begin{theorem}\cite[Theorem A.4]{LePMi20}\label{thm:graded}
    Suppose that $\mmm_h(\tau)$ is $(\Er,\tau,h)$-graded. Then, there exist $\tau_0, h_0>0$ such that, for all $0<\tau_j<\tau_0$ and $h\in(0,h_0]$, one has
    \bes
    \sigma(\mmm_h(\tau))\subset\bigsqcup_{j=1}^p\eps_j(\tau)^2(\sigma(M_j) + O(h)).
    \ees
    Moreover, for any eigenvalue $\lambda$ of $M_j$ with multiplicity $m_j(\lambda)$, there exists $K>0$ such that, denoting $D_j = \{z\in\C,\ |z-\eps_j(\tau)^2\lambda| < \eps_j(\tau)^2Kh\}$, one has
    \bes
    n(D_j;\mmm_h(\tau)) = m_j(\lambda),
    \ees
    where $n(D_j;\mmm_h(\tau))$ denotes the rank of the Riesz projector associated with $\mmm_h(\tau)$ along the contour $D_j$.
\end{theorem}

To apply this theorem to the interaction matrix \eqref{eq:matrix}, we must first show that it is graded. Let us first notice that
\bes
M = \begin{pmatrix}M'&0\\0&0\end{pmatrix}
\ees
where $M'\in\Mr_{n_0-1}(\R)$, since by construction $e_{n_0}$ is the ground state of $\Delta_V$. Now, for all $j\in\lb1,n_0-1\rb$, define
\be\label{eq:defvmu}
v(\m_j) = \sum_{\s\in\tilde\j(\m_j)}z(\m_j,\s),\ \mbox{ and }\ \mu_j = 2-2\k(\m_j)+2\k(\s)-\frac{2}{\nu_1^{\s}} = \mu(\m_j)
\ee
for some $\s\in\tilde\j(\m_j)$ (recalling $\tilde \j$ is introduced in \eqref{eq:deftildej}), where $z(\m_j,\s)$ is defined in \eqref{eq:zms}. Therefore,
\be\label{eq:lamdbav}
\tilde\lambda_j = v(\m_j)h^{\mu_j}e^{-2S_j/h}\big(1 + O(h^{\beta_j})\big),
\ee
denoting $\beta_j = \beta(\m_j)$ where $\beta$ is defined in \eqref{eq:beta}. Since the sequence $(S_j)_j$ is non-decreasing, there exists a partition $J_1\sqcup\ldots\sqcup J_p$ of $\lb1,n_0-1\rb$ such that, for all $k\in\lb1,p\rb$, there exists $\iota(k)\in\lb1,n_0-1\rb$ satisfying
\bes
\forall j\in J_k,\ \ S_j = S_{\iota(k)},\ \ \mbox{ and }\ \ \forall 1\leq k< k'\leq p,\ \ S_{\iota(k)} < S_{\iota(k')}
\ees
Therefore, using \eqref{eq:Pee} and \eqref{eq:lamdbav}, the matrix $\big(h^{\mu_1}e^{-2S_1/h}\big)^{-1}M'$ is $(\Er,\tau,h^{\beta''})$-graded with $\Er = (\R^{J_k})_{1\leq k\leq p}$, $\tau_k = h^{\frac{\mu_{\iota(k)}-\mu_{\iota(k-1)}}{2}} e^{-\frac{S_{\iota(k)}-S_{\iota(k-1)}}{h}}$ for $k\geq 2$ and
\bes
\beta'' = \min(\alpha_0,\min_{j\in\lb1,n_0-1\rb}\beta_j).
\ees
We can now apply Theorem \ref{thm:graded} and we obtain that
\bes
\sigma(M') \subset \bigsqcup_{k=1}^ph^{\mu_1}e^{-2S_1/h}\eps_k^2(\sigma(M_k) + O(h^{\beta''}))
\ees
with $M_k = \diag(v_j, j\in J_k)$, $h^{\mu_1}e^{-2S_1/h}\eps_k^2 = h^{\mu_{\iota(k)}}e^{-2S_{\iota(k)}/h}$ and corresponding multiplicities. This completes the proof of the stated result.

\section{Examples of admissible potentials}\label{sec:example}

A natural question arises: is it easy to satisfy Assumption \ref{ass.alpha}, and how can one verify it for a given potential?

In its utmost generality, one would have to check by hand if a given potential satisfies Assumption \ref{ass.alpha}. However, in this section, we present several criteria to ensure its validity. These criteria cover broad classes of potentials, and in most applications of this paper, the potential under consideration is likely to satisfy at least one of them.

Let us first consider the one-dimensional situation. That is when Assumption \ref{ass.case4} is satisfied. Thanks to \eqref{eq:estimecase4}, \eqref{eq:interaction} becomes
\bes
\<\Delta_vu_j,u_k\> = \delta_{j,k}\tilde\lambda_j + O\Big(h^\frac{2}{\overline{\nu}}\sqrt{\tilde\lambda_j\tilde\lambda_k}\Big).
\ees
We therefore see that when $d=1$, any potential is admissible.

Recall that for $\m \in \uuu^{(0)}$ and $\s \in \j(\m)$, the quantity $\gamma(\s)$ is defined in \eqref{eq:defgamma} as follows:
\bes
\gamma(\s) = \left\{\begin{array}{ll}
    +\infty &\mbox{ if }\s\in\j_\infty(\m),\\
    \max\Big(-\frac2{\underline{\nu}_*^\s},\min\big(-\frac{2}{\hat\nu^\s},2\frac{\hat\nu^\s-\underline\nu^\s}{\overline\nu^\s}-\frac{2}{\underline\nu^\s}\big)\Big) + \frac2{\overline{\nu}^\s} + \frac2{\overline{\nu}}&\mbox{ if }\s\in\j_1(\m),\\
    \frac2{\overline{\nu}^\s} - \frac2{\underline{\nu}_*^\s} + \frac2{\overline{\nu}} &\mbox{ if }\s\in\j_2(\m).
\end{array}\right.
\ees
Here, the notation for the $\nu$'s follows \eqref{eq:notationnu}, and the sets $\j_\infty$, $\j_1$, and $\j_2$ are defined in \eqref{eq:defjautre}. A straightforward computation shows that
\be\label{eq:conditionmin}
\min\big(-\frac{2}{\hat\nu^\s},2\frac{\hat\nu^\s-\underline\nu^\s}{\overline\nu^\s}-\frac{2}{\underline\nu^\s}\big) = -\frac{2}{\hat\nu^\s} \iff \hat\nu^\s\underline\nu^\s\geq\overline\nu^\s
\ee
because for $\s\in\j_1(\m)$, $\hat\nu^\s > \underline\nu^\s$.

Since $\alpha_0 = \min_{\m\in\uuu^{(0)}}\min_{\s\in\j(\m)}\gamma(\s)$, Assumption \ref{ass.alpha} is equivalent to
\bes
\forall \m\in\uuu^{(0)},\ \forall\s\in\j(\m),\ \gamma(\s) > 0.
\ees
Therefore, a potential $V$ satisfies Assumption \ref{ass.alpha} if, for all $\m \in \uuu^{(0)}$ and all $\s \in \j(\m)$, at least one of the following conditions holds. (We continue to assume that the first direction corresponds to the drop direction.) The specific criterion need not be the same for all saddle points.

\textbullet\ $\s$ satisfies Assumption \ref{ass.case1}; that is, the change of variable $U^\s$ is a unitary matrix near $\s$.

\textbullet\ $\s$ satisfies Assumption \ref{ass.case2}; that is, the drop direction is not degenerate, that is $\nu_1 = 2$.

In these two criteria, $\s\in\j_\infty(\m)$, hence the constructions made in the previous section are precise up to a factor $1 + O(h^\infty)$. For the remaining cases, this is no longer true; the factor instead takes the form $1+ O(h^\beta)$ with some $\beta\in\R$. We now present criteria that ensure $\beta > 0$.

\textbullet\ All degeneracies (except for the drop direction) are greater than half the maximum degeneracy among all critical points; that is, $\underline{\nu}_*^\s > \frac{\overline{\nu}}{2}$. Indeed, in this case, then, $\gamma(\s) > \frac2{\overline{\nu}^\s} - \frac2{\overline{\nu}} \geq 0$.

\textbullet\ $d=2$. To prove this, we consider two sub-cases:

$\star$ Either $\nu_1^\s \leq \nu_2^\s$, therefore we have $\overline{\nu}^\s = \underline{\nu}_*^\s$, and thus $\gamma(\s) \geq \frac2{\overline{\nu}}$.

$\star$ Or $\nu_2^\s < \nu_1^\s$, which is Assumption \ref{ass.case5}, this means that $\s\in\j_1(\m)$ and having that $\hat\nu^\s = \overline\nu^\s = \nu_1^\s$, we obtain from \eqref{eq:conditionmin} that
\bes
\gamma(\s) \geq -\frac{2}{\hat\nu^\s} + \frac2{\overline\nu^\s} + \frac2{\overline\nu} = \frac2{\overline\nu}.
\ees

\textbullet\ Similarly, we can generalize the previous point to the case where all the directions have the same degeneracy except for one, for which the degeneracy order is lower than the others. In other words $(\nu_i^\s)_i \subset \{\underline\nu^\s,\overline\nu^\s\}$, $\underline\nu^\s\neq\overline\nu^\s$ and $\exists!\nu_{i_0}^\s = \underline\nu^\s$.

$\star$ If $i_0 = 1$, then we have $\overline{\nu}^\s = \underline{\nu}_*^\s$, and thus $\gamma(\s) \geq \frac2{\overline{\nu}}$.

$\star$ Otherwise $\s$ satisfies Assumption \ref{ass.case5}, thus $\s\in\j_1(\m)$ and $\hat\nu^\s = \overline\nu^\s$, hence $\gamma(\s) \geq \frac2{\overline{\nu}}$.

\textbullet\ $V$ has only two different degeneracy degree at $\s$, one of which is $2$. From the second point, we know that if $\nu_1^\s=2$ then it is done, otherwise, Assumption \ref{ass.case3} is satisfied, hence $\s\in\j_1(\m)$. We then have $\forall i\in\lb1,d\rb,\ \nu_i^\s\in\{2,\overline\nu^\s\}$ and $\nu_1^\s = \overline\nu^\s$. Once again, since $\hat\nu^\s = \overline\nu^\s$, we obtain $\gamma(\s) \geq \frac{2}{\overline{\nu}}$.

\textbullet\ The saddle point is homogeneous after the change of variables, this means that $\forall i,j\in\lb1,d\rb,\ \nu_i^\s = \nu_j^\s$, therefore $\overline{\nu}^\s = \underline{\nu}_*^\s$, and thus $\gamma(\s) \geq \frac2{\overline{\nu}}$.

Although these criteria provide numerous ways to satisfy Assumption \ref{ass.alpha}, it is still relatively easy to construct a potential that fails to meet the requirement. For example, in dimension $d = 3$, consider a potential $V$ such that $V \circ U(x, y, z) = -x^4 + y^4 + z^8$, where $U$ is a change of variables which Taylor expansion lacks good properties. In this case, it is possible for $\gamma(0) \leq 0$ to occur.

\section{Triple well, same height}\label{sec:triple}

\begin{figure}[ht]
    \centering

\tikzset{every picture/.style={line width=0.75pt}} 

\begin{tikzpicture}[x=0.75pt,y=0.75pt,yscale=-1,xscale=1,scale=0.8]

\draw   (195,262.5) .. controls (157,312.5) and (132,333.5) .. (64,296.5) .. controls (-4,259.5) and (28,199.5) .. (45,151.5) .. controls (62,103.5) and (136,67.5) .. (202,113.5) .. controls (268,159.5) and (205,306.5) .. (313,316.5) .. controls (421,326.5) and (382,209.5) .. (441,149.5) .. controls (500,89.5) and (528,104.5) .. (573,148.5) .. controls (618,192.5) and (635,310.5) .. (502,322.5) .. controls (369,334.5) and (434,199.5) .. (400,136.5) .. controls (366,73.5) and (312,65.5) .. (268,124.5) .. controls (224,183.5) and (233,212.5) .. (195,262.5) -- cycle ;
\draw  [dash pattern={on 0.84pt off 2.51pt}] (77,95.5) .. controls (164,51.5) and (213,122.5) .. (242,121.5) .. controls (271,120.5) and (284,69.5) .. (330,65.5) .. controls (376,61.5) and (405,118.5) .. (429,122.5) .. controls (453,126.5) and (491,75.5) .. (544,103.5) .. controls (597,131.5) and (680.45,275.39) .. (560,320.5) .. controls (439.55,365.61) and (423,294.5) .. (404,290.5) .. controls (385,286.5) and (396,324.5) .. (323,331.5) .. controls (250,338.5) and (240,247.5) .. (223,247.5) .. controls (206,247.5) and (183,323.5) .. (147,327.5) .. controls (111,331.5) and (63,324.5) .. (34,295.5) .. controls (5,266.5) and (-10,139.5) .. (77,95.5) -- cycle ;
\draw    (215,37.5) -- (197,358.5) ;
\draw    (273,41.5) -- (255,362.5) ;
\draw    (422,62.5) -- (345,372.5) ;
\draw    (478,73.5) -- (401,383.5) ;

\draw[->] (40,70)--(75,90) ;
\draw (12,52) node[anchor=north west][inner sep=0.75pt]{$\supp\chi$};

\draw[->] (600,325)--(555,312) ;
\draw (580,330) node[anchor=north west][inner sep=0.75pt]{$\{V = \bsigma\}$};

\draw (150,150) node{$\times$} node[below right]{$\um$};
\draw (235.5,180.75) node{$\times$} node[below right]{$\s_1$};
\draw (320,240) node{$\times$} node[below right]{$\m_1$};
\draw (411,200.63) node{$\times$} node[below right]{$\s_2$};
\draw (540,200) node{$\times$} node[below right]{$\m_2$};

\draw (110,350) node{$\underline{E}$};
\draw (220,350) node{$\ccc_1$};
\draw (300,350) node{$E_1$};
\draw (370,350) node{$\ccc_2$};
\draw (440,350) node{$E_2$};

\end{tikzpicture}

    \caption{Representation of a triple well function near its critical points.}
    \label{fig:TripleWell}
\end{figure}

We could generalize the previous method and results following the work of \cite{Mi19} for potentials that do not satisfy \eqref{ass.gener}, but here we only consider an example to illustrate the method and highlight which results can or cannot be expected. We consider a potential $V$ that still satisfies Assumption \ref{ass.confin}, \ref{ass.morse}, and \ref{ass.morsebis}, such that $\uuu = \uuu^{(0)}\cup\uuu^{(1)}$, with $\uuu^{(0)} = \{\um,\m_1,\m_2\}$ and $\uuu^{(1)} = \{\s_1,\s_2\}$. Moreover, we assume that $V(\um) = V(\m_1) = V(\m_2) =: V_\m$, and $V(\s_1) = V(\s_2) =: \bsigma$. We immediately observe that \eqref{ass.gener} cannot be satisfied since one of its consequences is that $V$ must admit a unique global minimum. Such a case, assuming moreover that $V$ is Morse, has already been studied in \cite{Mi19}, where interactions between the wells were explicitly described, unlike in the generic case where \eqref{ass.gener} holds and no interaction occurs. Specifically, \cite[Theorem 1.1]{Mi19} states:
\begin{theorem}
    Let $V$ be such a potential that is also Morse. Then there exist $\eps,h_0>0$ such that for all $h\in(0,h_0]$, the eigenvalues of $\Delta_V$ in $[0,\eps h]$ consist of $0$ and others of the form
    \bes
    \lambda_k = hb_k(h)e^{-2S/h}.
    \ees
    The prefactors $b_k(h)$ satisfy $b_k(h) \sim \sum_{j\geq0}h^jb_{k,j}$, and the leading terms $b_{k,0}$ correspond to the non-zero eigenvalues of the discrete Laplacian of the graph whose vertices are the minima of $V$ and whose edges are the saddle points connecting two minima.
\end{theorem}

In this theorem, we observe that the principal term in the semiclassical expansion of the eigenvalues does not necessarily coincide with the principal term of the diagonal elements of the interaction matrix $\mmm = (\<\Delta_V\phii_{\m},\phii_{\m'}\>)_{\m,\m'\in\uuu^{(0)}}$. In this section, we show that for a non-Morse potential, this behavior may disappear, and under appropriate conditions, no interaction between the wells can be observed. More precisely, we have the following result:
\begin{assumption}\label{ass.triple}~

    \begin{itemize}
    \item[$i)$] $\k(\um) = \k(\m_1) = \k(\m_2) =: \k(\m)$, where $\k$ is defined in \eqref{eq:defik};
    \item[$ii)$] $\s_1$ and $\s_2$ are homogeneous, meaning that for each $\s\in\{\s_1,\s_2\}$, the values $\nu_i^\s$ are constant with respect to $i\in\lb1,d\rb$. We will denote them by $\nu^\s$ for short;
    \item[$iii)$] $\nu^{\s_2} < \nu^{\s_1} = \overline{\nu}$, $d>2$, and $\frac{\nu^{\s_2}}{\nu^{\s_1}} > 1-\frac 2d$.
\end{itemize}
\end{assumption}
\begin{theorem}\label{thm.triple}
    Let $V$ be a potential satisfying Assumption \ref{ass.triple}. Then the interaction matrix $\mmm = (\<\Delta_V\phii_{\m},\phii_{\m'}\>)_{\m,\m'\in\uuu^{(0)}}$, with the quasimodes set in definition \ref{def:triplephii}, has the form:
    \bes
    \mmm = h^{2-2\k(\m)+2\k(\s_1)-\frac2{\nu^{\s_1}}}e^{-2S/h}\begin{pmatrix}
        0&0\\
        0&M
    \end{pmatrix}
    \ees
    where $M$ has the following expansion
    \be\label{eq:thmtriple}
    M = \begin{pmatrix}\alpha&0\\0&0\end{pmatrix} + h^{(d-2)\big(\frac1{\nu^{\s_2}}-\frac1{\nu^{\s_1}}\big)}M_1 + O(h^{\frac{2}{\nu^{\s_1}}})
    \ee
    where $\alpha\in\R^*$ and $M_1$ is a symmetric $2\times2$ matrix. Both are $h$-independent and can be computed explicitly.
\end{theorem}

This theorem shows a splitting of eigenvalues, contrary to what one might expect. In the Morse case, there is some interaction as explained above Assumption \ref{ass.triple}. Here we have no such effect just as if \eqref{ass.gener} was true. However, if we relax Assumption \ref{ass.triple} $iii)$ and let $d=2$ or $\nu^{\s_2} = \nu^{\s_1}$, then the exponent of $h$ in front of $M_1$ in \eqref{eq:thmtriple} becomes zero, and we recover an interaction, as in the Morse case.

To work with this potential, we follow the approach of \cite{No24}. We begin by defining several cutoff functions:

Let $\delta>0$, and define $\chi$ such that $\chi\equiv 1$ on $\{V\leq \bsigma + \delta\}$ and $\supp\chi\subset\{V\leq \bsigma + 2\delta\}$. We split $\supp\chi$ into five regions, each containing a critical point, following the geometric constructions from Section 4. Reading Figure \ref{fig:TripleWell} from left to right, we have $\um\in\underline{E},\ \s_1\in\ccc_1,\ \m_1\in E_1,\ \s_2\in\ccc_2$, and $ \m_2\in E_2$, where $\underline{E}$, $E_1$, and $E_2$ denote $E(\um)$, $E(\m_1)$, and $E(\m_2)$ respectively, following the notation of \eqref{def:labeling}. Similarly, $\ccc_1$ and $\ccc_2$ denote $\ccc_{\s_1,3\tau,3\delta}$ and $\ccc_{\s_2,3\tau,3\delta}$, respectively, as defined in \eqref{eq:ccc}. In the coordinates where $\s_j=0$, we define on $\ccc_{\s_j,3\tau,3\delta}$:
\bes
\tilde{\theta}_j^\pm(x) = C_j^{-1}\int_0^{\pm\ell_j(x)}\zeta(r/\tau)e^{-\frac{r^{\nu_1^{\s_j}}}{\nu_1^{\s_j}h}}dr
\ees
where $\ell_j$ is the function constructed in Section \ref{sec:sharp}, $\nu_1^{\s_j}$ denotes the decay rate in the decreasing direction at $\s_j$, and $C_j$ is a normalizing constant satisfying \eqref{eq:ch}. The sign of $\ell$ is chosen (see Remark \ref{rem:signell}) so that the following cutoff funcitons are smooth:
\bes
\theta_{\um} = \left\{\begin{array}{ccc}
0&\mbox{ on }&E_1\cup \ccc_2 \cup E_2\\
2&\mbox{ on }&\underline{E},\\
1+\tilde{\theta}_1^+&\mbox{ on }&\ccc_1.\\
\end{array}\right.\qquad
\theta_{\m_2} = \left\{\begin{array}{ccc}
0&\mbox{ on }&\underline{E}\cup \ccc_1 \cup E_1,\\
2&\mbox{ on }&E_2,\\
1+\tilde{\theta}_2^+&\mbox{ on }&\ccc_2.\\
\end{array}\right.
\ees

\bes
\theta_{\m_1} = \left\{\begin{array}{ccc}
0&\mbox{ on }&\underline{E}\cup E_2,\\
2&\mbox{ on }& E_1,\\
1+\tilde{\theta}_1^-&\mbox{ on }&\ccc_1,\\
1+\tilde{\theta}_2^-&\mbox{ on }&\ccc_2.\\
\end{array}\right.
\ees

We define the elementary functions
\bes
g_\m(x) = h^{-\k(\m)}c(\m)\chi(x)\theta_\m(x)e^{-(V(x)-V_\m)/h},
\ees
for $\m \in \uuu^{(0)}$. Here, $c(\m)$ is such that $\vvert{g_\m} = 1$. Applying the Laplace method, we obtain $c(\m) = c_0(\m)(1+O(h^{\frac2{\overline{\nu}^\m}}))$, where
\bes
c_0(\m)^{-2} = 4\prod_{i=1}^d\frac{2\Gamma(\frac{1}{\nu_i^{\m}})}{\nu_i^{\m}(2t_i^\m)^{\frac{1}{\nu_i^{\m}}}}.
\ees
This definition extends to any critical point $x^* \in \uuu$ via
\bes
c_0(x^*)^{-2} = 4\prod_{i=1}^d\frac{2\Gamma(\frac{1}{\nu_i^{x^*}})}{\nu_i^{x^*}(2|t_i^{x^*}|)^{\frac{1}{\nu_i^{x^*}}}}.
\ees

Let $\fff = \R^{\uuu^{(0)}}\simeq\R^3$ be the finite-dimensional vector space of real-valued functions on $\uuu^{(0)}$, endowed with the following scalar product
\bes
\<\eta,\eta'\>_{\fff} = \sum_{\m\in\uuu^{(0)}}\eta(\m)\eta'(\m).
\ees

We construct the following orthonormal basis of $\fff$, where $\widehat{c},\ \widehat{c^1}$, and $ \widehat{c^2}$ denote normalizing constants:
\bes
\forall \m\in\uuu^{(0)}, \ \ \eta_{\um}(\m) = \frac{\widehat{c}}{c(\m)}.
\ees

\bes
\eta_{\m_1} : \left\{\begin{array}{l}
\eta_{\m_1}(\um) = \frac{\widehat{c^1}}{c(\m_1)},\\
\eta_{\m_1}(\m_1) = -\frac{\widehat{c^1}}{c(\um)},\\
\eta_{\m_1}(\m_2) = 0.\\
\end{array}\right.\qquad
\eta_{\m_2} : \left\{\begin{array}{l}
\eta_{\m_2}(\um) = \frac{\widehat{c^2}}{c(\um)},\\
\eta_{\m_2}(\m_1) = \frac{\widehat{c^2}}{c(\m_1)},\\
\eta_{\m_2}(\m_2) = -\frac{\widehat{c^2}}{c(\m_2)}\gamma,\\
\end{array}\right.
\ees
where $\gamma = \frac{c(\m_2)^2}{c(\um)^2} + \frac{c(\m_2)^2}{c(\m_1)^2}$. We are now ready to construct the corresponding quasimodes:

\begin{defin}\label{def:triplephii} The ground state is defined by
    \bes
    \phii_{\um} = \widehat c\, h^{-\k(\um)}e^{-(V-V_\m)/h}(1+O(h^\beta)),
    \ees
    for some $\beta>0$. For $\m\in\{\m_1,\m_2\}$, we set
    \bes
    \phii_\m = \sum_{\m'\in\uuu^{(0)}}\eta_\m(\m')g_{\m'}.
    \ees
\end{defin}

We have the following proposition:
\begin{proposition}
    The family $(\phii_\m)_{\m\in\uuu^{(0)}}$ is quasi-orthonormal, meaning that there exists a constant $C>0$ such that for all $\m,\m'\in\uuu^{(0)}$,
    \bes
    \<\phii_\m,\phii_{\m'}\> = \delta_{\m,\m'} + O(e^{-C/h}).
    \ees
\end{proposition}
\bp Consider
\bes
\tilde\phii_{\um} = \sum_{\m'\in\uuu^{(0)}}\eta_{\um}(\m')g_{\m'} = h^{-\k(\um)}\widehat c\sum_{\m'\in\uuu^{(0)}}\theta_{\m'}\chi e^{-(V-V_\m)/h}.
\ees
From \cite[Proposition 4.1]{No24}, this holds for the family $(\tilde\phii_{\um},\phii_{\m_1},\phii_{\m_2})$. By Assumption \ref{ass.triple} $i)$, for $\eps>0$ small enough, we have $\tilde\phii_{\um} \equiv \phii_{\um}$ on $\uuu^{(0)} + B(0,\eps)$. Moreover, by the definition of $\uuu^{(0)}$, there exists $c>0$ such that $V - V_\m \geq c$ on $\R^d\setminus\big(\uuu^{(0)} + B(0,\eps)\big)$, which therefore leads to $\tilde\phii_{\um} - \phii_{\um} = O(e^{-c/h})$. When then deduce the announced result.

\ep

For any $\m,\m'\in\uuu^{(0)}\setminus\{\um\}$, we have
\bes
\<\Delta_V\phii_\m,\phii_{\m'}\> = \sum_{\tilde{\m},\tilde{\m}'}\eta_\m(\tilde{\m})\eta_{\m'}(\tilde{\m}')\nnn_{\tilde{\m},\tilde{\m}'},
\ees
where
\bes
\nnn_{\tilde{\m},\tilde{\m}'} = h^{-\k(\tilde{\m})-\k(\tilde{\m}')}c(\tilde{\m})c(\tilde{\m}')e^{-2S/h}\<\Delta_V\chi\theta_{\tilde{\m}}e^{-(V-\bsigma)/h},\chi\theta_{\tilde{\m}'}e^{-(V-\bsigma)/h}\>
\ees
with $S = \bsigma - V_\m$, as defined in \eqref{eq:defS}. Applying with another Laplace method, we obtain
\bes
\<\Delta_V\phii_\m,\phii_{\m'}\> = h^2e^{-2S/h}\sum_{\tilde{\m},\tilde{\m}'}\eta_\m(\tilde{\m})\eta_{\m'}(\tilde{\m}')N_{\tilde{\m},\tilde{\m}'}(1+O(h^{\frac{2}{\overline{\nu}}}))
\ees
where
\bes
N_{\tilde{\m},\tilde{\m}'} = h^{-\k(\tilde{\m})-\k(\tilde{\m}')}c(\tilde{\m})c(\tilde{\m}')(-1)^{1-\delta_{\tilde{\m},\tilde{\m}'}}\sum_{\s\in \j(\tilde{\m})\cap \j(\tilde{\m}')}h^{2\k(\s)-\frac{2}{\nu_1^{\s}}}z_1(\s),
\ees
with
\bes
z_1(\s) = \Bigg(\frac{\nu_1^\s(2|t_1^\s|)^{\frac{1}{\nu_1^\s}}}{2\Gamma(\frac{1}{\nu_1^\s})}\Bigg)^2c_0(\s)^{-2}.
\ees
To ensure the validity of this computation, we must guarantee that sub-principal terms associated with some $\nnn_{\tilde \m,\tilde \m'}$ do not dominate the principal terms of other entries $\nnn_{\hat \m,\hat \m'}$. While this may be difficult to verify in the general case, under Assumption \ref{ass.triple}, the condition becomes straightforward to check.

The leading sub-principal term has the form
\be\label{eq:subprinc}
ah^{2-2\k(\m)+2\k(\s_1)-\frac2{\nu_1^{\s_1}}+\frac{2}{\overline{\nu}}}e^{-2S/h}
\ee
for some $a\neq0$ $h$-independent while the smallest principal one is of the form
\be\label{eq:princ}
bh^{2-2\k(\m)+2\k(\s_2)-\frac2{\nu_1^{\s_2}}}e^{-2S/h}
\ee
(with the same restriction on $b$ as on $a$). But $h^{2\k(\s)-\frac{2}{\nu_1^{\s}}}$ can be simplified to $h^{\frac{d-2}{\nu^{\s}}}$, therefore we can verify that for \eqref{eq:subprinc} to remain smaller than \eqref{eq:princ}, the following inequality must hold
\bes
\frac{d-2}{\nu^{\s_1}} + \frac2{\overline{\nu}} > \frac{d-2}{\nu^{\s_2}} \iff \frac{\nu^{\s_2}}{\nu^{\s_1}} > 1 - \frac{2}{d},
\ees
using the identity $\overline{\nu} = \nu^{\s_1}$, which justifies $iii)$ of Assumption \ref{ass.triple}.

We now have
\bes
\<\Delta_V\phii_\m,\phii_{\m'}\> = h^{2-2\k(\m)}e^{-2S/h}\sum_{\tilde{\m},\tilde{\m}'}\eta_\m(\tilde{\m})\eta_{\m'}(\tilde{\m}')\tilde N_{\tilde{\m},\tilde{\m}'}(1+O(h^{\frac{2}{\overline{\nu}}}))
\ees
where $h^{-2\k(\m)}\tilde N_{\tilde{\m},\tilde{\m}'} = N_{\tilde{\m},\tilde{\m}'}$ (we just factor the common power of $h$ to simplify the following calculus). We are now in a position to explicitly compute the interaction matrix.

Let us start with $\<\Delta_V\phii_{\m_1},\phii_{\m_2}\>$. Computing each term of the sum over $\uuu^{(0)}$, noticing that $\j(\um) = \{\s_1\},\ \j(\m_1) = \{\s_1,\s_2\}$, and $\j(\m_2) = \{\s_2\}$, we have
\bes
\begin{split}
    \eta_{\m_1}(\um)\eta_{\m_2}(\um)\tilde N_{\um,\um} &= h^{\frac{d-2}{\nu^{\s_1}}}\widehat{c^1}\widehat{c^2}z_1(\s_1)\frac{c(\um)}{c(\m_1)},\\
    \eta_{\m_1}(\m_1)\eta_{\m_2}(\um)\tilde N_{\m_1,\um} &= h^{\frac{d-2}{\nu^{\s_1}}}\widehat{c^1}\widehat{c^2}z_1(\s_1)\frac{c(\m_1)}{c(\um)},\\
    \eta_{\m_1}(\um)\eta_{\m_2}(\m_1)\tilde N_{\um,\m_1} &= -h^{\frac{d-2}{\nu^{\s_1}}}\widehat{c^1}\widehat{c^2}z_1(\s_1)\frac{c(\um)}{c(\m_1)},\\
    \eta_{\m_1}(\m_1)\eta_{\m_2}(\m_1)\tilde N_{\m_1,\m_1} &= -h^{\frac{d-2}{\nu^{\s_1}}}\widehat{c^1}\widehat{c^2}z_1(\s_1)\frac{c(\m_1)}{c(\um)}\\
    &\phantom{******} - h^{\frac{d-2}{\nu^{\s_2}}}\widehat{c^1}\widehat{c^2}z_1(\s_2)\frac{c(\m_1)}{c(\um)},\\
    \eta_{\m_1}(\um)\eta_{\m_2}(\m_2)\tilde N_{\um,\m_2} &= 0,\\
    \eta_{\m_1}(\m_1)\eta_{\m_2}(\m_2)\tilde N_{\m_1,\m_2} &= - h^{\frac{d-2}{\nu^{\s_2}}}\widehat{c^1}\widehat{c^2}z_1(\s_2)\frac{c(\m_1)}{c(\um)}\gamma.\\
\end{split}
\ees
Where we did not write the three terms involving $\eta_{\m_1}(\m_2)$ because they are worth zero. We thus obtain
\bes
\<\Delta_V\phii_{\m_1},\phii_{\m_2}\> = - h^{2-2\k(\m)+\frac{d-2}{\nu^{\s_2}}}e^{-2S/h}\widehat{c^1}\widehat{c^2}z_1(\s_2)\frac{c(\m_1)}{c(\um)}(1+\gamma).
\ees
The key observation here is that there is no $\s_1$ appearing in this expression, which may be expected since the interaction was analyzed around $\s_2$, yet $\um$ still appears, and $\s_1$ would have appeared as well if the $\k(\m)$ were not all equal.

Now for $\<\Delta_V\phii_{\m_2},\phii_{\m_2}\>$
\bes
\begin{split}
    \eta_{\m_2}(\um)\eta_{\m_2}(\um)\tilde N_{\um,\um} &= h^{\frac{d-2}{\nu^{\s_1}}}\widehat{c^2}^2z_1(\s_1),\\
    \eta_{\m_2}(\m_1)\eta_{\m_2}(\um)\tilde N_{\m_1,\um} &= - h^{\frac{d-2}{\nu^{\s_1}}}\widehat{c^2}^2z_1(\s_1),\\
    \eta_{\m_2}(\m_2)\eta_{\m_2}(\um)\tilde N_{\m_2,\um} &= 0,\\
    \eta_{\m_2}(\m_1)\eta_{\m_2}(\m_1)\tilde N_{\m_1,\m_1} &= h^{\frac{d-2}{\nu^{\s_1}}}\widehat{c^2}^2z_1(\s_1) + h^{\frac{d-2}{\nu^{\s_2}}}\widehat{c^2}^2z_1(\s_2),\\
    \eta_{\m_2}(\m_1)\eta_{\m_2}(\m_2)\tilde N_{\m_1,\m_2} &= h^{\frac{d-2}{\nu^{\s_2}}}\widehat{c^2}^2z_1(\s_2)\gamma,\\
    \eta_{\m_2}(\m_2)\eta_{\m_2}(\m_2)\tilde N_{\m_2,\m_2} &= h^{\frac{d-2}{\nu^{\s_2}}}\widehat{c^2}^2z_1(\s_2)\gamma^2.\\
\end{split}
\ees
Counting two times the terms that appear twice, we obtain
\bes
\<\Delta_V\phii_{\m_2},\phii_{\m_2}\> = h^{2-2\k(\m)+\frac{d-2}{\nu^{\s_2}}}e^{-2S/h}\widehat{c^2}^2z_1(\s_2)(1 + \gamma)^2.
\ees
Once again we have no $\s_1$ appearing there. Finally, we need to check that the term of order $h^{\frac{d-2}{\nu^{\s_1}}}$ in $\<\Delta_V\phii_{\m_1},\phii_{\m_1}\>$ does not vanish:
\bes
\begin{split}
    \eta_{\m_1}(\um)\eta_{\m_1}(\um)\tilde N_{\um,\um} &= h^{\frac{d-2}{\nu^{\s_1}}}\big(\widehat{c^1}\frac{c(\um)}{c(\m_1)}\big)^2z_1(\s_1),\\
    \eta_{\m_1}(\m_1)\eta_{\m_1}(\um)\tilde N_{\m_1,\um} &= h^{\frac{d-2}{\nu^{\s_1}}}\widehat{c^1}^2z_1(\s_1),\\
    \eta_{\m_1}(\m_1)\eta_{\m_1}(\m_1)\tilde N_{\m_1,\m_1} &= \big(\widehat{c^1}\frac{c(\m_1)}{c(\um)}\big)^2\big(h^{\frac{d-2}{\nu^{\s_1}}}z_1(\s_1) + h^{\frac{d-2}{\nu^{\s_2}}}z_1(\s_2)\big),
\end{split}
\ees
the terms involving $\m_2$ being zero.

This concludes the proof of Theorem \ref{thm.triple}, since $\Delta_V\phii_{\um} = 0$. We also obtain explicit expressions for $\alpha$ and $M_1$, namely:
\be\label{eq:alphatriple}
\alpha = \widehat{c^1}^2z_1(\s_1)\Big(\frac{c(\um)}{c(\m_1)} + \frac{c(\m_1)}{c(\um)}\Big)^2
\ee
and
\bes
M_1 = z_1(\s_2)\begin{pmatrix}
    \Big(\widehat{c^1}\frac{c(\m_1)}{c(\um)}\Big)^2 & -\widehat{c^1}\widehat{c^2}\frac{c(\m_1)}{c(\um)}(1+\gamma)\\ -\widehat{c^1}\widehat{c^2}\frac{c(\m_1)}{c(\um)}(1+\gamma) & \widehat{c^2}^2(1+\gamma)^2
\end{pmatrix}.
\ees

Let $M_1 = \begin{pmatrix}m_1&m_2\\m_2&m_3\end{pmatrix}$. Then, $M = \begin{pmatrix}\alpha +h^\beta m_1&h^\beta m_2\\h^\beta m_2&h^\beta m_3\end{pmatrix}$, where $\beta = (d-2)\big(\frac1{\nu^{\s_2}}-\frac1{\nu^{\s_1}}\big)$. Let $\lambda$ be an eigenvalue of $M$; then it is a root of
\bes
(\lambda - \alpha - h^\beta m_1)(\lambda - h^\beta m_3) - h^{2\beta}m_2^2.
\ees
If $\lambda$ is sufficiently far from $\alpha$, then the first factor is of order $O(1)$. Consequently, the second term must be of order $h^{2\beta}$, which implies that $\lambda = h^\beta m_3(1 + O(h^\beta))$. Hence
\bes
\sigma(\mmm) = \{0,\alpha,m_3h^\beta\}h^\mu e^{-2S/h}(1+O(h^\beta)),
\ees
where $\alpha$ is given by \eqref{eq:alphatriple}, $m_3 = z_1(\s_2)\widehat{c^2}^2(1+\gamma)^2$, and $\mu = 2-2\k(\m)+\frac{d-2}{\nu^{\s_1}}$. These values, up to a small perturbation, correspond to the diagonal entries of the matrix $\mmm$.  As in the proof of Theorem \ref{thm:2}, the small eigenvalues of $\Delta_V$ are given by the eigenvalues of$\mmm$ up to a multiplicative factor of the form $1 + O(e^{-c/h})$. We therefore observe a spectral splitting between the wells in the eigenvalues of $\Delta_V$; as if the generic assumption \eqref{ass.gener} were satisfied.

\newpage
\appendix
\section{Some technical results}

\subsection{Essential spectrum of $\Delta_V$}\label{ssec:essspec}
Under Assumption \ref{ass.confin}, for all points outside the compact set $K$,
\bes
|\nabla V|^2-h\Delta V \geq |\nabla V|^2-hC|\nabla V|^2 \geq \frac{1-hC}{C^2} \geq \frac1{2C^2}
\ees
for $h$ small enough. Hence, because $\liminf (|\nabla V|^2-h\Delta V) \geq \frac1{2C^2} =: C_0$, we have $\sigma_{ess}(\Delta_V)\subset[C_0,+\infty)$ (see, for instance, the examples following \cite[Theorem 4.19]{DiSj99_01}).

\subsection{Proof of Lemma \ref{lem:soussolaffine}}
Without loss of generality, assume that $V\geq 0$. By Assumption \ref{ass.confin}, there exists $R>0$ such that $|\nabla V(x)| \geq \frac1C$ for all $|x|\geq R$

Let $x_0\in\R^d$. We consider the maximal solution to following the Cauchy problem
\be\label{eq:edo}
\left\{
\begin{aligned}
&\dot x(t) = -\nabla V(x(t)),
\\&x(0) = x_0.
\end{aligned}
\right.
\ee
We then see that
\bes
\frac{d}{dt}(V(x(t))) = -|\nabla V(x(t))|^2.
\ees

Let us first show that if $\lambda>0$ satisfies $V(x_0) \leq\lambda$ and $|x_0| > R + \lambda C$, then the solution to \eqref{eq:edo} is defined on $\R_+$. Moreover, for all $s\geq0$, $|x(s)|\geq R$. Let $t$ be the supremum of time for which $|x(s)|\geq R$ holds. Suppose, for contradiction, that $t<\infty$. Then,
\bes
\begin{split}
    |x(t) - x(0)| &= \bigg|\int_0^t-\nabla V(x(s))ds\bigg| \leq \int_0^t|\nabla V(x(s))|ds\\
    &\leq C\int_0^t|\nabla V(x(s))|^2ds = C(V(x_0) - V(x(t)) \leq \lambda C,
\end{split}
\ees
which leads to $|x(t)|\leq |x_0| + \lambda C$. But we also have that
\bes
|x(t)| \geq |x_0| - |x(t) - x_0| > R + \lambda C - \lambda C = R.
\ees
This contradicts the maximality of $t$.

Now assume that such an $x_0$ exists. Then we have
\bes
0\leq V(x(t)) = V(x_0) + \int_0^td(V(x(s))) = V(x_0) - \int_0^t|\nabla V(x(s))|^2ds \leq \lambda - \frac1{C^2}t.
\ees
Since $x(t)$ is defined for all $t\geq0$, the inequality becomes a contradiction for sufficiently large $t$. This leads to
\be\label{eq:sublevel}
\{V\leq\lambda\}\subset \overline B(0,R+\lambda C)
\ee
where the right-hand side denotes the closed ball centered at $0$ with radius $R+\lambda C$.

Now let $x\in\R^d$ be such that $|x| > R$. Then, there exists $\eps>0$ such that $\lambda = \frac1C(|x| - R - \eps) > 0$, hence we can write $|x| = R + \lambda C + \eps$. By \eqref{eq:sublevel}, we know that $V(x) > \lambda = \frac1C|x| + a$, with $a = -\frac{R+\eps}{C}$. Therefore, there exists $\tilde a\in\R$ such that
\be\label{eq:borneinfV}
\forall |x| > R,\ \ V(x) > \frac1C|x| + \tilde a.
\ee
Because $\{|x|\leq R\}$ is compact, we can find $b\in\R$ such that \eqref{eq:borneinfV} holds on the whole space.

\subsection{Proof of Proposition \ref{prop:saddle}}
Before proceeding, we define a family of paths. We will then prove the following lemma, adapted from \cite[Lemma 3.1]{AsBoMi22}. Let $a=(a_1,\ldots,a_d)\in[1,+\infty)^d$ and $x,y\in\R^d$. We define
\bes
\begin{array}{c|ccc}
     \gamma^a_{x,y}:&[0,1]&\to&\R^d\\
     &t&\mapsto&(t^{a_i}(x_i-y_i))_i + y.
\end{array}
\ees
Observe that $\gamma^a_{x,y}(0) = y$ and $\gamma^a_{x,y}(1) = x$. When $a\equiv 1$, this corresponds to the standard linear parametrization of the segment $[y,x]$.
\begin{lemma}\label{lem:path}
    Let $\phii$ be a local smooth diffeomorphism of $\R^d$ defined in a neighborhood of $b\in\R^d$, and $a=(a_1,\ldots,a_d)\in[1,+\infty)^d$. Then there exists $r_b>0$ such that for all $0<r<r_b$,
    \bes
    \forall x\in\phii(B(b,r)),\ \gamma^a_{x,\phii(b)}([0,1])\subset\phii(B(b,r)).
    \ees
\end{lemma}

\bp
In the following, we denote $\underline a = \inf_ia_i$ and $\overline a = \sup_ia_i$. The goal is to show that for all $x\in B(0,r)$ and all $t\in[0,1]$, the point $\gamma^a_{\phii(b+x),\phii(b)}(t)$ belongs to $\phii(B(b,r))$. Define
\bes
\forall t\in[0,1],\ g(t) = |\phii^{-1}(\gamma^a_{\phii(b+x),\phii(b)}(t))-b|^2\in[0,r^2).
\ees
We aim to show that $g(t)<r^2$ for all $t\in[0,1]$. First, we see that $g(0) = 0$ and $g(1) = |x|^2 < r^2$. Then, there exists $E$ a set large enough such that
\bes
\begin{split}
    g(t) &= |\phii^{-1}(\gamma^a_{\phii(b+x),\phii(b)}(t))-b|^2\\
    &\leq \vvert{d\phii^{-1}}_{L^\infty(E)}^2 |\gamma^a_{\phii(b+x),\phii(b)}(t)-\phii(b)|^2\\
    &= \vvert{d\phii^{-1}}_{L^\infty(E)}^2 \sum_i|t^{a_i}(\phii(b+x)_i-\phii(b)_i)|^2\\
    &\leq \vvert{d\phii^{-1}}_{L^\infty(E)}^2\vvert{d\phii}_{L^\infty(B(b,r))}^2 t^{2\underline a}|x|^2.
\end{split}
\ees
Thus there exists $C>0$ such that for all $\eps>0$ and all $t\leq\eps$,
\bes
g(t) \leq C\eps^{2\underline a}|x|^2.
\ees
Choosing $\eps$ such that $C\eps^{2\underline a} \leq 1$, we have $g(t) < r^2$ for all $t\in[0,\eps]$. Moreover, the Taylor formula implies
\bes
\begin{split}
    g'(t) &= 2\< \partial_t(\phii^{-1}\circ\gamma^a_{\phii(b+x),\phii(b)})(t),\phii^{-1}\circ\gamma^a_{\phii(b+x),\phii(b)}(t) - b\>\\
    &= 2\< d_{\gamma^a_{\phii(b+x),\phii(b)}(t)}\phii^{-1}({\gamma^a_{\phii(b+x),\phii(b)}}'(t)),\phii^{-1}\circ\gamma^a_{\phii(b+x),\phii(b)}(t) - b\>\\
    &= 2\< d_{\phii(b)+O(t^{\underline a}x)}\phii^{-1}\big((a_it^{a_i-1}d_b\phii_i(x))_i + O(t^{\underline a-1}x^2)\big),\\&\phantom{************}d_{\phii(b)}\phii^{-1}\big((t^{a_i}d_b\phii_i(x))_i\big) + O(t^{2\underline a}x^2)\>.
\end{split}
\ees
Here we denoted $\phii_i$ the $i$-th component of $\phii$. Denoting $A = d_{\phii(b)}\phii^{-1}$ and $\Gamma_1(t) = \diag(t^{a_i-1})$, $\Gamma_2(t) = \diag(a_it^{a_i-1})$, we observe that
\bes
g'(t) = 2t\<A\Gamma_2(t)A^{-1}x,A\Gamma_1(t)A^{-1}x\> + \vvert{A\Gamma_1(t)A^{-1}x}O(t^{\underline a}x^2) + O(t^{2\underline a-1}x^3).
\ees
But using that $A$ is invertible and thus $A^*A$ is positive definite,
\bes
\begin{split}
    \<A^*A\Gamma_2(t)A^{-1}x,\Gamma_1(t)A^{-1}x\> &\geq c\<\Gamma_2(t)A^{-1}x,\Gamma_1(t)A^{-1}x\>\\
    &\geq c't^{2\overline a-2}\<(AA^*)^{-1}x,x\>\\
    &\geq c''t^{2\overline a-2}|x|^2,
\end{split}
\ees
for some $c,c',c''>0$. Thus we have for $\eps$ defined above
\bes
\exists 0<c_{\eps}<1,\exists r_{\eps}>0,\forall |x|<r_{\eps},\forall t\geq\eps,\ \ g'(t) \geq c_{\eps}|x|^2
\ees
and it leads to
\bes
\forall t\geq\eps,\ g(t) \leq (1-(1-t)c_{\eps})|x|^2 < r_{\eps}^2.
\ees

\ep

We can now prove the proposition.

\bp
Without loss of generality, we may assume that $\s=0$ and $V(\s)=0$. By Assumption \ref{ass.morse}, there exist a smooth diffeomorphism $U$ and two sequences $(t_i)_i\subset\R^*, (\nu_i)_i\subset\N\setminus\{0,1\}$ such that $U$ is defined in a neighborhood of $0$, $U(0) = 0$, and for $x$ in that neighborhood,
\bes
V\circ U (x) = \sum_{i=1}^dt_ix_i^{\nu_i}.
\ees

We denote $X_0 = \{V<0\}$. Let $\overline \nu = \sup\negthinspace_i\nu_i$, and set $a = (\overline\nu/\nu_i)_i\in[1,+\infty)^d$. Let $x\in U^{-1}(X_0\cap B(0,r))$ with $0<r<r_0$, $r_0$ given by Lemma \ref{lem:path}. Then, we observe that for all $s\in[0,1]$, the following holds:
\bes
 V\circ U (\gamma^a_{x,0}(s)) = \sum_{i=1}^dt_i(s^{\frac{\overline\nu}{\nu_i}}x_i)^{\nu_i} = s^{\overline\nu}V\circ U(x),
\ees
hence
\be\label{eq:a}\tag{a}
\gamma^a_{x,0}((0,1])\subset U^{-1}(X_0).
\ee
Moreover, by applying Lemma \ref{lem:path} to $\phii = U^{-1}$ and $b=0$, we obtain
\be\label{eq:b}\tag{b}
\gamma^a_{x,0}((0,1])\subset U^{-1}(B(0,r)),
\ee
thus, combining \eqref{eq:a} and \eqref{eq:b}, we have
\be\label{eq:c}\tag{c}
\gamma^a_{x,0}((0,1])\subset U^{-1}(X_0\cap B(0,r)).
\ee
Next, note that
\be\label{eq:d}\tag{d}
\gamma^a_{x,0}(s)\xrightarrow[s\to0]{}0
\ee
and
\be\label{eq:e}\tag{e}
\forall \eta \ll 1,\ U^{-1}(X_0\cap B(0,r))\cap B(0,\eta) = U^{-1}(X_0)\cap B(0,\eta)
\ee
because $U^{-1}$ is an open map since it is a diffeomorphism. Then, up to proving that
\be\label{eq:f}\tag{f}
U^{-1}(X_0)\cap B(0,\eta) \mbox{ is connected}
\ee
we have
\bes
X_0\cap B(0,r) \mbox{ is connected}
\ees
by combining \eqref{eq:c}, \eqref{eq:d}, \eqref{eq:e} and \eqref{eq:f}, and using the fact that connectedness is preserved under diffeomorphisms Therefore, it just remains to prove \eqref{eq:f}. Now, let us show that if $0\notin\uuu^{(1)}$, then
\bes
U^{-1}(\{V<0\})\cap B(0,\eta) = \{x\in B(0,\eta),\ V\circ U(x)<0\}
\ees
is connected for $\eta>0$ small enough.

We consider the case where $0\in\uuu^{odd}$ which means that there is at least one odd $\nu$. Without any loss of generality, we can consider that $\nu_1$ is odd and $t_1<0$ (in the case $t_1>0$, we just replace $\eta/2$ by $-\eta/2$ in the following).

Start from a point $(x_1,\ldots,x_d)$ in $\{x\in B(0,\eta),\ V\circ U(x)<0\}$, the goal is to connect it to $(\eta/2,0,\ldots,0)$ via a path that remains within the set. In the following, by saying that we link $a$ to $b$ we mean that we create a segment from one to another (so a path of the form $t\mapsto a + t(b-a)$), one can check that each time we do that, the whole path stays in $\{x\in B(0,\eta),\ V\circ U(x)<0\}$.

The first step of the path is to link all the coordinates $x_i$ such that $t_ix_i^{\nu_i}$ is positive to zero. Thus, either $x_1=0$ or $t_1x_1^{\nu_1}<0$ and all other non-zero contributions to $V\circ U$ are negative ones. Then we halve all the $x_i$ (so we link $x_i$ to $x_i/2$), this way, the ending point is in $B(0,\eta/2)$ and we know that its first coordinate is non-negative. Now we link $x_1$ to $\eta/2$, because we were in $B(0,\eta/2)$ and we had $x_1\geq0$, we remain in $B(0,\eta)$. At last, we link all the other coordinates to $0$, this way, we finally reached the desired point.

For $0\in\uuu^{even}$ the setting is extremely similar to the well-known Morse case and the proof is exactly the same.

Recall that for open subsets of a Euclidean space, arc-connectedness and connectedness are equivalent, we then have proven the proposition.

\ep

\subsection{Degenerate Laplace method}
This method generalizes  the standard Laplace method to functions that are not necessarily Morse. We assume only that they satisfy Assumption \ref{ass.morse}. For a more general introduction to the theory of Laplace methods, we refer to \cite{ArGuVa88}, in particular to Sections 7.3.3 and 8.3.2.

\begin{proposition}\label{LaplaceMethod}
    Let $\vvv\subset\R^d$ be bounded, $\phii:\vvv\to\R$, $x_0\in\vvv$ such that $x_0$ is the unique global minimum of $\phii$. Suppose there exist a neighborhood $\www\subset\R^d$ of $0$, and a smooth change of variable $U:\www\to\vvv$ satisfying $U(0)=x_0$, $|\det d_0U|\neq0$. Suppose also there exist sequences $(t_i)_{1\leq i\leq d}\subset\R_+^*$ and $\nu = (\nu_i)_{1\leq i\leq d}\subset2\N^*$ such that
    \be\label{eq:Laplace1}
    \forall x\in\www,\ \ \phii(U(x)) - \phii(x_0) = \sum_{i=1}^d t_ix_i^{\nu_i}.
    \ee
    Then for any $a_h\sim \sum_{j\geq0}h^ja_j\in\ccc^\infty(\overline{\vvv})$,
    \be\label{eq:Laplace2}
    \int_{\vvv}a_h(x)e^{-\frac{\phii(x)-\phii(x_0)}{h}}dx = a_0(x_0)|\det d_0U|h^{\sum_{i=1}^d\frac1{\nu_i}}\prod_{i=1}^d\frac{2\Gamma(\frac{1}{\nu_i})}{\nu_it_i^{\frac{1}{\nu_i}}}(1+O(h^{\min_{1\leq i\leq d}\frac{2}{\nu_i}})).
    \ee
    Additionally, if $a_h = O((x-x_0)^{2\alpha})$ uniformly with respect to $h$, we have
    \be\label{eq:Laplace3}
    \int_{\vvv}a_h(x)e^{-\frac{\phii(x)-\phii(x_0)}{h}}dx = O(h^{\sum_{i=1}^d\frac{2\alpha_i+1}{\nu_i}}).
    \ee
\end{proposition}

\bp
Let $\ccc\subset\www$ be a cubic neighborhood of $0$ where \eqref{eq:Laplace1} holds. We split
\bes
\int_{\vvv}a_h(x)e^{-\frac{\phii(x)-\phii(x_0)}{h}}dx = \int_{U\ccc}a_h(x)e^{-\frac{\phii(x)-\phii(x_0)}{h}}dx + \int_{\vvv\setminus U\ccc}a_h(x)e^{-\frac{\phii(x)-\phii(x_0)}{h}}dx.
\ees
Since $x_0$ is the unique global minimum of $\phii$, it follows that
\bes
\int_{\vvv\setminus U\ccc}a_h(x)e^{-\frac{\phii(x)-\phii(x_0)}{h}}dx = O(e^{-\delta/h})
\ees
for some $\delta>0$ using that $\vvv$ is bounded. For the remaining integral, we apply the change of variables $x\mapsto U(x)$,
\bes
\int_{U\ccc}a_h(x)e^{-\frac{\phii(x)-\phii(x_0)}{h}}dx = \int_{\ccc}a_h(U(x))e^{-\frac1h\sum_{i=1}^d t_ix_i^{\nu_i}}|\det d_xU|dx
\ees
and we denote $b_h(x) = a_h(U(x))|\det d_xU|$,
\bes
T=\begin{pmatrix}\big(\frac{t_1}{h}\big)^{\frac{1}{\nu_1}}&&\\&\ddots&\\&&\big(\frac{t_d}{h}\big)^{\frac{1}{\nu_d}}\end{pmatrix}
\ees
thus $\det T^{-1} = h^{\sum_{i=1}^d\frac1{\nu_i}}\prod_{i=1}^dt_i^{-\frac{1}{\nu_i}}$. We therefore have
\bes
\begin{split}
    \int_\ccc b_h(x)e^{-\frac1h\sum_{i=1}^d t_ix_i^{\nu_i}}dx &= \det T^{-1}\int_{T\ccc} b_h(T^{-1}x)e^{-x^{\nu}}dx\\
    &\sim \det T^{-1}\sum_{j\geq0}h^j\int_{T\ccc}b_j(T^{-1}x)e^{-x^{\nu}}dx,
\end{split}
\ees
where $b_j = a_j\circ U|\det dU|$. Now, using a Taylor expansion, for any $n,j\in\N$,
\bes
\int_{T\ccc} b_j(T^{-1}x)e^{-x^{\nu}}dx = \sum_{|\gamma|\leq n}\frac{\partial^{2\gamma}(b_j\circ T^{-1})(0)}{(2\gamma)!}\int_{T\ccc} x^{2\gamma}e^{-x^{\nu}}dx + O(h^{\min\frac{2n+2}{\nu_i}}),
\ees
noticing that $\int_{T\ccc} x^{\gamma}e^{-x^{\nu}}dx = 0$ if $\gamma\notin 2\N^d$, since $\nu\subset2\N^*$. Using that
\bes
\partial^{2\gamma}(b_j\circ T^{-1})(0) = h^{\sum_{i=1}^d\frac{2\gamma_i}{\nu_i}}\prod_{i=1}^dt_i^{-\frac{2\gamma_i}{\nu_i}}\partial^{2\gamma}b_j(0)
\ees
and choosing $n=0$, we just need to compute $\int_{T\ccc}e^{-x^{\nu}}$. However, using that
\bes
T\ccc = \prod_{i=1}^d[-h^{-\frac{1}{\nu_i}}c_i,h^{-\frac{1}{\nu_i}}c_i]
\ees
for some $c_i>0$, then
\bes
\begin{split}
    \int_{T\ccc}e^{-x^{\nu}}dx &= \int_{\R^d}e^{-x^{\nu}}dx - \prod_{i=1}^d2\int_{x_i\geq c_ih^{-\frac{1}{\nu_i}}}e^{-x_i^{\nu_i}}dx_i\\
    \int_{x_i\geq ch^{-\frac{1}{\nu_i}}}e^{-x_i^{\nu_i}}dx_i &\leq \int_{x_i\geq ch^{-\frac{1}{\nu_i}}}e^{-(ch^{-\frac{1}{\nu_i}})^{\nu_i-1}x_i}dx_i\\
    &= c'h^{1-\frac{1}{\nu_i}}\int_{x_i\geq c'h^{-1}}e^{-x_i}dx_i = O(e^{-\tilde c/h}).
\end{split}
\ees
Hence there exists $c>0$ such that $\int_{T\ccc}e^{-x^{\nu}}dx = \int_{\R^d}e^{-x^{\nu}}dx + O(e^{-c/h})$. Moreover, a straightforward computation shows that
\bes
\begin{split}
    \int_{\R^d}e^{-x^{\nu}}dx &= \prod_{i=1}^d\int_\R e^{-x_i^{\nu_i}}dx_i = \prod_{i=1}^d2\int_{\R_+}e^{-x_i^{\nu_i}}dx_i\\
    &= \prod_{i=1}^d\frac{2}{\nu_i}\int_{\R_+} e^{-x_i}x_i^{\frac{1}{\nu_i}-1}dx_i = \prod_{i=1}^d\frac{2\Gamma(\frac{1}{\nu_i})}{\nu_i}.
\end{split}
\ees
Altogether with $b_0(0) = a_0(U(0))|\det d_0U| = a_0(x_0)|\det d_0U|$ we obtain \eqref{eq:Laplace2}. In the second case ($a_h = O((x-x_0)^{2\alpha})$), we notice that for any $|2\gamma|\leq |2\alpha|-1$, $\partial^{2\gamma}b(0) = 0$, therefore choosing $n=|\alpha|$, we have
\bes
\int_{T\ccc}b_h(T^{-1}x)e^{-x^{\nu}}dx = O(h^{\sum_{i=1}^d\frac{2\alpha_i+1}{\nu_i}})
\ees
and thus we obtain \eqref{eq:Laplace3}.

\ep

\begin{remark}
    If $\nu_i = 1$ for all $1 \leq i \leq d$, then $\phii$ is Morse and we recover the standard Laplace method by noting that $\prod_{i=1}^dt_i = \det\Hess_{x_0}\frac\phii2$.
\end{remark}

\begin{remark}
    In \eqref{eq:Laplace3}, the assumption $a_h = O((x-x_0)^{2\alpha})$ can be replaced with $a_h = O((x-x_0)^{\alpha})$. Due to the parity of the terms involved, the result becomes
    \bes
    \int_{\vvv}a_h(x)e^{-\frac{\phii(x)-\phii(x_0)}{h}}dx = O(h^{\sum_{i=1}^d\frac{\alpha'_i+1}{\nu_i}}),
    \ees
    where $\alpha'_i = \left\{\begin{array}{cc}
        \alpha_i & \mbox{ if }\alpha_i\in2\N \\
        \alpha_i+1 & \mbox{ if }\alpha_i\notin2\N.
    \end{array}\right.$
\end{remark}

\vfill
\bibliographystyle{siam}

\end{document}